\documentclass[11pt,a4paper]{article}
\usepackage[T1]{fontenc}
\usepackage[utf8]{inputenc}
\usepackage{lmodern,amsmath,amssymb,amsthm,mathtools}
\usepackage[margin=28mm]{geometry}
\usepackage{microtype,booktabs,enumitem}
\usepackage[hidelinks]{hyperref}
\hypersetup{pdftitle={Thresholds of the Pólya–Szegő inequality for the fractional logarithmic p-Laplacian},pdfauthor={Hichem Hajaiej}}
\newtheorem{theorem}{Theorem}[section]
\newtheorem{proposition}[theorem]{Proposition}
\newtheorem{lemma}[theorem]{Lemma}
\newtheorem{corollary}[theorem]{Corollary}
\theoremstyle{definition}\newtheorem{definition}[theorem]{Definition}
\theoremstyle{remark}\newtheorem{remark}[theorem]{Remark}
\newcommand{\R}{\mathbb R}
\newcommand{\X}{X^{s,\log,p}(\R^N)}
\newcommand{\Lsp}{\mathcal L_{s,p}}
\newcommand{\diam}{\operatorname{diam}}
\newcommand{\supp}{\operatorname{supp}}
\newcommand{\esssupp}{\operatorname{ess\,supp}}
\newcommand{\Rmon}{R_{\mathrm{mon}}}
\newcommand{\Rpos}{R_{\mathrm{pos}}}
\newcommand{\Dsch}{D_{\mathrm{Sch}}}
\numberwithin{equation}{section}
\makeatletter
\renewcommand{\@seccntformat}[1]{\csname the#1\endcsname.\quad}
\makeatother
\title{Thresholds of the P\'olya--Szeg\H{o} inequality for the fractional logarithmic $p$-Laplacian}
\author{Hichem Hajaiej\\[4pt]
\small Department of Mathematics\\
\small California State University, Los Angeles\\
\small Los Angeles, CA 90032, USA\\
\small\href{mailto:hhajaie@calstatela.edu}{hhajaie@calstatela.edu}}
\date{}
\begin{document}
\maketitle
\begin{abstract}
We establish three surprising features of the P\'olya--Szeg\H{o} principle for the full signed energy $\mathcal L_{s,p}$ of the fractional logarithmic $p$-Laplacian. Let $N\ge1$, $1<p<\infty$, and let $0<s<1$ be the fractional order. For a nonnegative function $u$ in the energy space, $u^H$ denotes its polarization with respect to a half-space $H$, and $u^*$ its Schwarz symmetric decreasing rearrangement. Set $d(u)=\diam(\esssupp u)$, and let $\Rmon$ be the radius up to which the interaction kernel is strictly decreasing.

\emph{First, Schwarz symmetrization decreases the energy on small supports but can increase it on larger supports.} In every dimension we prove
\[
 \begin{aligned}
 d(u)\le\Rmon&\quad\Longrightarrow\quad\Lsp(u^*)\le\Lsp(u),\\
 \Lsp(v^*)&>\Lsp(v)\quad\text{for some }0\le v\in C_c^\infty(\R^N).
 \end{aligned}
\]
For every $\ell>2\Rmon$, such a counterexample can be chosen with $2\Rmon<d(v)<\ell$.

\emph{Second, polarization and Schwarz symmetrization have different sharp thresholds on the line.} Their sharp uniform diameter bounds are $\Rmon$ and $2\Rmon$, respectively. We prove
\[
 \begin{aligned}
 d(u)\le\Rmon&\ \Longrightarrow\ \Lsp(u^H)\le\Lsp(u)
       &&(N\ge1),\\
 d(u)\le2\Rmon&\ \Longrightarrow\ \Lsp(u^*)\le\Lsp(u)
       &&(N=1).
 \end{aligned}
\]
The polarization bound holds for every $H$ and is sharp in every dimension; the Schwarz bound is sharp on the line. For every $\Rmon<\ell\le2\Rmon$, there are a smooth nonnegative function $u$ on the line and a half-line $H$ such that
\[
 \Rmon<d(u)<\ell,\qquad
 \Lsp(u^*)<\Lsp(u)<\Lsp(u^H).
\]
Thus an individual polarization can strictly increase the energy while the final Schwarz rearrangement strictly decreases it.

\emph{Third, the sharp Steiner threshold changes with dimension.} It is $2\Rmon$ for $N=1$, where Steiner and Schwarz symmetrization coincide, but only $\Rmon$ for every $N\ge2$.

We also prove the sharp doubled Schwarz bound in every dimension for functions with ball superlevel sets whose centers may vary with the level, with equality precisely for translates of the Schwarz rearrangement. An exact cumulative Riesz-deficit formula characterizes the unrestricted higher-dimensional Schwarz threshold, whose explicit value remains open within $\Rmon\le \Dsch\le2\Rmon$. Complete equality classifications accompany the main comparison theorems.
\end{abstract}


\section{Introduction}
The classical P\'olya--Szeg\H{o} inequality states that symmetric decreasing rearrangement does not increase the Dirichlet energy of a nonnegative Sobolev function. More precisely, if $u\in W^{1,p}(\R^N)$, $1<p<\infty$, and $u^*$ denotes its Schwarz symmetric decreasing rearrangement, then
\[
 \int_{\R^N}|\nabla u^*|^p\,dx\le \int_{\R^N}|\nabla u|^p\,dx.
\]
This energy-decreasing property is a cornerstone of rearrangement methods in variational analysis; see \cite{HLP,LL,Baernstein}. Related inequalities for nonlocal energies follow from the Riesz rearrangement principle \cite{BLL,Burchard}. The fractional Sobolev comparison is recalled in Corollary~\ref{cor:fractionalinput}; see also \cite[Theorem~A.1]{FS} and \cite{DNPV}. Polarization, a two-point rearrangement used to approximate Steiner and Schwarz symmetrizations, is discussed in \cite{BS,HApprox}.

Three phenomena distinguish the full signed energy studied here. First, Schwarz rearrangement decreases this energy under a suitable support-diameter bound (Theorem~\ref{thm:schwarz}), while smooth functions with larger supports can give the reverse inequality (Proposition~\ref{prop:schwarzcore}). Support size and superlevel-set geometry therefore enter the comparison. The essential-support diameter is defined in \eqref{eq:esssupportdiameter}.

Second, a single polarization and full Schwarz rearrangement have different optimal diameter bounds in dimension one: $\Rmon$ and $2\Rmon$, respectively. Corollary~\ref{cor:step-endpoint-separation} exhibits the resulting strict separation: a polarization increases the energy while Schwarz rearrangement decreases it.

Third, the sharp Steiner diameter changes with dimension. It is $2\Rmon$ in dimension one, where Steiner and Schwarz symmetrization coincide (Theorem~\ref{thm:schwarzonedimensional}), and $\Rmon$ in every dimension $N\ge2$ (Proposition~\ref{prop:sharpsteiner}). These thresholds guarantee comparison for every admissible function below the bound; exceeding a bound allows counterexamples but does not force every function to violate the inequality. The radii are defined in \eqref{eq:introradii}.

Fix $N\ge1$, $1<p<\infty$, and $0<s<1$.
The fractional logarithmic $p$-Laplacian considered here is obtained by differentiating the fractional $p$-Laplacian with respect to its order at a fixed $s\in(0,1)$; see \cite{BGHO}. To specify the operator and its normalization, write $g_p(a)=|a|^{p-2}a$, with $g_p(0)=0$, and let $C_t=C(N,t,p)>0$ be differentiable near $s$. For sufficiently regular functions, the fractional $p$-Laplacian is
\begin{equation}\label{eq:introfractionalp}
 (-\Delta)^t_pu(x)=C_t\,\operatorname{PV}\!\int_{\R^N}
 \frac{g_p(u(x)-u(y))}{|x-y|^{N+tp}}\,dy.
\end{equation}
We define its fixed-positive-order logarithmic derivative by
\begin{equation}\label{eq:intrologoperator}
 \begin{aligned}
 \mathfrak L_{s,p}u
 &:=\left.\partial_t\bigl[(-\Delta)^t_pu\bigr]\right|_{t=s},\\
 \mathfrak L_{s,p}u(x)
 &=C_s\,\operatorname{PV}\!\int_{\R^N}
 g_p(u(x)-u(y))\,
 \frac{b_s-p\log|x-y|}{|x-y|^{N+sp}}\,dy,
 \qquad b_s=\frac{C_s'}{C_s}.
 \end{aligned}
\end{equation}
Here $\operatorname{PV}$ means the limit of the integrals over
$\R^N\setminus B_\varepsilon(x)$ as $\varepsilon\downarrow0$.
The pointwise formula is used only where this limit exists and order
differentiation is justified. In the energy space introduced below, the
operator is understood weakly through
\begin{equation}\label{eq:introweakoperator}
 \langle\mathfrak L_{s,p}u,\varphi\rangle
 =\frac{C_s}{2}\iint_{\R^{2N}}
 g_p(u(x)-u(y))(\varphi(x)-\varphi(y))
 k_s(|x-y|)\,dx\,dy,
 \qquad\varphi\in C_c^\infty(\R^N),
\end{equation}
where
\begin{equation}\label{eq:introkernel}
 k_s(r)=(b_s-p\log r)r^{-N-sp}.
\end{equation}
The corresponding full-space scalar energy is
\begin{equation}\label{eq:introenergy}
 \Lsp(u)=\frac{C_s}{2p}\iint_{\R^{2N}}
 |u(x)-u(y)|^p k_s(|x-y|)\,dx\,dy.
\end{equation}
This distinguishes the operator $\mathfrak L_{s,p}$ from its scalar energy
$\mathcal L_{s,p}$: their relation is
$D\mathcal L_{s,p}(u)[\varphi]=\langle\mathfrak L_{s,p}u,\varphi\rangle$.
The concrete normalization and its differentiability restrictions are given
in Section~\ref{sec:definitions}.

Write $q=N+sp$. Two distinct radii govern this problem:
\begin{equation}\label{eq:introradii}
 \Rpos=e^{b_s/p},\qquad
 \Rmon=e^{b_s/p+1/q}.
\end{equation}
The kernel changes sign at $\Rpos$ and stops decreasing at $\Rmon>\Rpos$.
For nonnegative functions it is the second radius that controls
polarization. Comparison therefore persists beyond the sign change of
the kernel. The negative tail is retained in every theorem; replacing
the kernel by its positive part would give a different functional.

The new theory concerns this full signed energy at a fixed positive
order. The classical positive-kernel argument does not apply because
the kernel eventually increases, and differentiating the fractional
P\'olya--Szeg\H{o} inequality does not preserve its direction. We instead
identify the occupied pairs that contribute to the sorting defect,
control their distances, and pass through the logarithmic singularity
while keeping the signed remainder. This yields comparison, sharp
counterexamples, and rigidity within one framework.

\par
The resulting sharp diameter bounds distinguish a single sorting step from the fully symmetrized function. Theorem~\ref{thm:polar} and Proposition~\ref{prop:sharpdiameter} prove that
$\Rmon$ is the largest diameter guaranteeing every polarization
comparison. Theorem~\ref{thm:schwarzonedimensional} proves the exact
Schwarz threshold $2\Rmon$ in dimension one. {For every diameter allowance between these two radii, Corollary~\ref{cor:step-endpoint-separation} constructs a smooth function for which a polarization strictly increases the energy while Schwarz rearrangement strictly decreases it.} This separates the final
rearrangement inequality from monotonicity of its approximating sorting
steps. The new ingredient is a correction that relates the line
interaction to an interaction on a circle and admits its own
rearrangement comparison; the circle Riesz inequality itself is the
classical input \cite[Theorem~2]{BH}. For every $\ell>2\Rmon$,
Proposition~\ref{prop:schwarzcore} constructs a smooth translated-core
counterexample with support diameter less than $\ell$, without a restriction
$sp<1$. For Steiner symmetrization, the transverse directions available when $N\ge2$ allow a different counterexample, giving the sharp bound $\Rmon$ (Proposition~\ref{prop:sharpsteiner}).

\par
Additional comparison ranges emerge when the geometry of the superlevel sets is prescribed. Theorem~\ref{rem:ballsuperlevelthreshold} establishes the sharp
$2\Rmon$ bound in every dimension for functions whose positive
superlevel sets are balls. Their centers may vary with the level;
radial symmetry of the function is not assumed. Equality forces one
common center, even at the critical diameter. For a different geometric class, Theorem~\ref{rem:affinethreshold} proves the larger sufficient bound
\[
 \diam(\esssupp u)\le2e^{\delta_p}\Rmon,
 \qquad \delta_p>0,
\]
for common linear deformations of smooth radial decreasing profiles.
{Here $\delta_p$ is the positive series in \eqref{eq:translationmarginseries}; its proof and numerical evaluation are in Appendix~\ref{sec:schwarzthreshold}.}
This affine class consists of common linear deformations; the ball-superlevel class permits centers to vary with the level. The two bounds concern different geometric families.
The underlying uniform moment gap excludes the entire linear-deformation
family from unrestricted threshold-minimizing sequences, even when both
the profile and the linear map vary (Corollary~\ref{cor:affineexcluded}).
Optimality of the affine bound remains open.

\par
To understand the scale dependence of these comparisons, we use the cumulative form of the Riesz rearrangement inequality \cite{BLL,Burchard}. Theorem~\ref{thm:schwarzdeficit} expresses the signed Schwarz energy gap in terms of a nonnegative cumulative Riesz deficit. Its sign is determined exactly by
whether a geometric mean $g_u$ of the deficit-weighted interaction
distances lies below or above $\Rmon$. Every bounded-support nonnegative
admissible nonsymmetric profile consequently has one critical dilation: the usual
comparison is strict below that scale and strictly reversed above it.
At the critical scale equality holds for a nonsymmetric function, which
explains why small-support rigidity cannot extend to arbitrary supports.
Write $\Dsch$ for the largest support diameter guaranteeing the Schwarz comparison for every nonnegative function in the energy space. Optimizing the profile criterion gives
\[
 \Dsch=\Rmon e^{\theta_{N,s,p}},
 \qquad 0\le\theta_{N,s,p}\le\log2.
\]
Here $\theta_{N,s,p}$ is an explicit logarithmic moment infimum
(Corollary~\ref{cor:schwarzvariational}). The parameter
$\theta_{N,s,p}$ and the ratio $\Dsch/\Rmon$ are independent
of the chosen positive differentiable normalization; the unnormalized
diameter need not be.

\par
The logarithmic-moment formulation also clarifies equality at the critical diameter. Theorem~\ref{thm:ordermoment} proves that the logarithmic moment of every
fixed nonsymmetric profile increases strictly with $s$ over its admissible orders, while its
infimum over all profiles is nondecreasing. For smooth profiles the
derivative is a strictly positive variance. In dimension one the
infimum is nevertheless the constant $\log2$. It is therefore not
attained by any nonsymmetric profile at any positive order.
This gives the complete equality characterization at diameter
$2\Rmon$ (Corollary~\ref{cor:onedimensionalendpoint}), beyond the
strict-kernel-decrease regime of the polarization theorem.
The comparison theorems also give common-center rigidity without an
assumption excluding plateaus or critical sets. Such assumptions can
play a role for local gradient energies \cite{BZ,CF,CCF}.
Theorem~\ref{thm:quantdeficit} further bounds the logarithmic gap by the
ordinary fractional rearrangement deficit. Theorem~\ref{thm:twoorderdeficit} strengthens this estimate using fractional deficits at two orders. In dimension one and for ball superlevel sets in every dimension, it gives an explicit strictly positive, profile-dependent bound at diameter $2\Rmon$, where the one-order estimate has a zero right-hand side. The quantitative proof is contained in the appendix.

\par
Beyond the power-difference energy, the sorting argument applies to certain more general interactions. In the language of submodular rearrangement interactions \cite{BH,Hajaiej10}, Theorems~\ref{thm:generalG} and~\ref{thm:Phi} extend the method to
submodular value interactions and to integrands depending jointly on
distance and the two function values. Under the stated growth and integrability assumptions, the latter theorem
uses decrease in distance of a
four-value sorting gap. Neither separation into a kernel and a value
interaction nor positivity of that gap is required. Strictness gives
the corresponding equality classifications. This identifies which part
of the classical rearrangement mechanism survives for signed
interactions, rather than imposing positivity on the whole integrand.

\par
These extensions require structural assumptions on the interaction; monotonicity and convexity in a local energy-density variable alone are insufficient. Theorem~\ref{thm:nonlinearcounter} constructs a smooth compactly
supported $u\ge0$, with support diameter at most $\Rpos/2$, and a smooth
integrand $j$ satisfying $\partial_tj>0$ and $\partial_{tt}j>0$, such that
\[
 \int_{\R^N}j(u^*(x),e_{u^*}(x))\,dx
 >\int_{\R^N}j(u(x),e_u(x))\,dx.
\]
Here $e_u$ is the pointwise density of the full energy, defined in
\eqref{eq:localdensity}. The construction applies when $sp<1$.
The scalar energy itself obeys the small-support comparison in this
regime. Thus separate monotonicity and strict convexity in the density
variable do not supply the coupling needed for a generalized
nonlocal P\'olya--Szeg\H{o} theorem. Theorem~\ref{thm:affinedensity}
provides a positive result under a suitable convex-weight condition.

\begin{samepage}
The uniform diameter thresholds for unrestricted nonnegative functions are:
\begin{center}
\begin{tabular}{lll}
\toprule
Operation & Dimension & Uniform diameter threshold\\
\midrule
Polarization & $N\ge1$ & $\Rmon$ (sharp)\\
Steiner & $N=1$ & $2\Rmon$ (sharp)\\
Steiner & $N\ge2$ & $\Rmon$ (sharp)\\
Schwarz & $N=1$ & $2\Rmon$ (sharp)\\
Schwarz & $N\ge2$ & $\Rmon\le \Dsch\le2\Rmon$\\
\bottomrule
\end{tabular}
\end{center}
\end{samepage}
For $N\ge2$, determining the exact unrestricted Schwarz threshold remains open. The sharp doubled bound is proved in dimension one and for ball superlevel sets in every dimension. Proposition~\ref{prop:doubledcriterion} identifies the additional estimate needed to extend it to arbitrary higher-dimensional profiles.

For local gradient energies, polarization preserves the relevant
value--gradient pairs; generalized local inequalities are discussed in
\cite{Hajaiej10}. For standard fractional energies, a nonnegative
decreasing kernel orders the reflected interactions at every scale.
Our small-support proof replaces that global ordering by control of the
active distances in the defect identity. The constructive approximation
by successive polarizations is based on Hajaiej
\cite[unnumbered Theorem, p.~410]{HApprox}; its slice extension for Steiner
symmetrization is proved in Theorem~\ref{thm:polarapprox}. The new analytic passage
combines lower semicontinuity of the singular positive part with
$L^p$ continuity of the integrable signed remainder. This also retains
the full zero exterior interaction on bounded domains.

The logarithmic operator at $s=0$ has a different, renormalized energy
\cite{CW,DJS}. The unrestricted endpoint comparison and rigidity of
Hajaiej, Li and Zhan \cite[Theorems~1.1--1.2]{HLZ} are recalled precisely
in Theorem~\ref{thm:endpointliterature}. They coexist with the
positive-order counterexamples proved here. For regular compactly
supported functions, a justified limit $s\downarrow0$ recovers the
endpoint comparison because $\Rmon(s)\to\infty$
(Remark~\ref{rem:main-to-endpoint}); equality does not transfer through
that limiting argument alone. At an interior positive order, even the
comparison cannot be obtained by differentiating the fractional
inequality: a positive fractional gap can have a negative order
derivative. The dilation counterexamples make this obstruction explicit.

In view of these differences, the determination of the exact universal Schwarz threshold for $N\ge2$ is
considerably more delicate than the corresponding polarization and Steiner
problems. The kernel changes sign and is not globally decreasing
(Lemma~\ref{lem:kernel}). Corollary~\ref{cor:step-endpoint-separation}
proves in dimension one that an individual polarization can strictly
increase the energy even when the final Schwarz rearrangement strictly
decreases it. Theorem~\ref{thm:schwarzdeficit}, specifically
\eqref{eq:schwarzdeficitformula}, expresses the signed energy gap in terms
of a nonnegative cumulative Riesz deficit weighted by
$\log(\Rmon/r)$, which changes sign. Corollary~\ref{cor:schwarzvariational}
then characterizes the universal threshold by a logarithmic-moment
infimum over admissible profiles. The remaining difficulty is to control
the balance between the two distance ranges uniformly over all profiles;
Proposition~\ref{prop:doubledcriterion} gives the precise balance estimate
equivalent to the unrestricted doubled threshold. The ball-superlevel
and affine comparison theorems
(Theorems~\ref{rem:ballsuperlevelthreshold} and~\ref{rem:affinethreshold})
apply to their stated geometric classes and do not establish that estimate
for arbitrary higher-dimensional functions. These criteria and supporting
estimates are presented in Appendix~\ref{sec:schwarzthreshold}.

The preliminary geometric and approximation results, including their precise hypotheses, are collected in Section~\ref{sec:classicalinputs}; Lemma~\ref{lem:positivekernel} supplies the positive-kernel comparison used in the proofs. The main text establishes the comparison theorems, equality cases, and sharp diameter bounds. Appendix~\ref{sec:schwarzthreshold} contains the cumulative deficit representation and its consequences. Corollary~\ref{cor:boundedPS} gives the bounded-domain comparison, and Appendix~\ref{sec:appendix} contains the extensions and counterexamples for more general interactions.

\section{Definitions and preliminary results}\label{sec:definitions}
\subsection{Order differentiation and normalization}
Fix $N\ge1$, $1<p<\infty$, and $0<s<1$. Following \cite{BGHO}, for $0<t<1$ we define the fractional $p$-Laplacian normalization by
\begin{equation}\label{eq:normalization}
 C(N,t,p)=
 \begin{cases}
 \displaystyle\frac{tp\,2^{2(t-1)}\Gamma\!\left(\frac{N+tp}{2}\right)}
 {\pi^{(N-1)/2}\Gamma(1-t)\Gamma\!\left(\frac{p+1}{2}\right)},
 &\displaystyle\frac12<t<1,\\[10pt]
 \displaystyle\frac{tp\,2^{2t-1}\Gamma\!\left(\frac{N+tp}{2}\right)}
 {\pi^{N/2}\Gamma(1-t)},
 &\displaystyle 0<t\le\frac12,
 \end{cases}
\end{equation}
where $\Gamma$ is the Gamma function. Set $C_t=C(N,t,p)$. The two branches differ at $t=\frac12$ by the factor
$\sqrt\pi/[2\Gamma((p+1)/2)]$; hence this piecewise normalization is not generally continuous there. Throughout the order-differentiation statements, we assume $s\ne\frac12$ unless $p=2$, in which case the branches coincide and
\[
 C(N,t,2)=\frac{t\,2^{2t}\Gamma(N/2+t)}{\pi^{N/2}\Gamma(1-t)}.
\]
On either open branch, logarithmic differentiation gives
\begin{equation}\label{eq:normalizationderivative}
 \frac{C_s'}{C_s}=\frac1s+2\log2
 +\frac p2\psi\!\left(\frac{N+sp}{2}\right)+\psi(1-s),
\end{equation}
where $\psi=\Gamma'/\Gamma$ is the digamma function. More generally, the comparison results below apply to any positive normalization continuously differentiable near the fixed order $s$. Write
\begin{equation}\label{eq:radii}
 b_s=\frac{C_s'}{C_s},\qquad q=N+sp,\qquad
 \Rpos=e^{b_s/p},\qquad \Rmon=e^{b_s/p+1/q}.
\end{equation}
\begin{remark}[Explicit evaluation of the monotonicity radius]
For the normalization \eqref{eq:normalization}, substitution of \eqref{eq:normalizationderivative} into \eqref{eq:radii} gives
\[
 \displaystyle
 \Rmon=\exp\!\left[
 \frac1{ps}+\frac{2\log2}{p}
 +\frac12\psi\!\left(\frac{N+sp}{2}\right)
 +\frac{\psi(1-s)}p+\frac1{N+sp}\right].
\]
This formula is used at the differentiability points specified above. The radius depends on $N,s,p$; numerical values for $p=2$, rounded to four decimal places, are
\begin{center}
\begin{tabular}{ccrrr}
\toprule
$N$ & $p$ & $s=\frac14$ & $s=\frac12$ & $s=\frac34$\\
\midrule
1 & 2 & 9.7177 & 2.5163 & 0.6265\\
2 & 2 & 11.4328 & 2.8950 & 0.7086\\
3 & 2 & 12.9315 & 3.2310 & 0.7824\\
\bottomrule
\end{tabular}
\end{center}
For example, using $\psi(1/2)=-\gamma-2\log2$ and $\psi(3/2)=\psi(1/2)+2$, where $\gamma$ is Euler's constant, we obtain
\[
 N=p=2,\quad s=\frac12:
 \qquad \Rmon={\frac12 e^{7/3-\gamma}}\approx2.8950.
\]
Thus the support-diameter hypothesis in this case reads
{$\diam(\esssupp u)\le\tfrac12 e^{7/3-\gamma}$.}

For fixed $N,p$, the expansions $\psi(1-s)=-\gamma+O(s)$ as $s\downarrow0$ and $\psi(1-s)=-(1-s)^{-1}-\gamma+O(1-s)$ as $s\uparrow1$ give
\[
 \begin{aligned}
 \log\Rmon&=\frac1{ps}+O(1),&&s\downarrow0,\\
 \log\Rmon&=-\frac1{p(1-s)}+O(1),&&s\uparrow1.
 \end{aligned}
\]
Consequently $\Rmon\to\infty$ as $s\downarrow0$, whereas $\Rmon\to0$ as $s\uparrow1$. These endpoint statements describe the geometric radius; they do not assert convergence of the signed energies or their rearrangement inequalities at the endpoints.
\end{remark}
The normalized fractional energy and its order derivative are
\begin{align}
 F_t(u)&=\frac{C_t}{2p}\iint_{\R^{2N}}
 \frac{|u(x)-u(y)|^p}{|x-y|^{N+tp}}\,dx\,dy,\label{eq:F}\\
 \Lsp(u)&=\frac{C_s}{2p}\iint_{\R^{2N}}
 |u(x)-u(y)|^p k_s(|x-y|)\,dx\,dy.\label{eq:L}
\end{align}
For smooth compactly supported functions, differentiation under the integral gives $\Lsp(u)=\partial_t F_t(u)|_{t=s}$. The first variation of \eqref{eq:L} corresponds to the order derivative of the fractional $p$-Laplacian, with the normalization in \eqref{eq:F}. We fix a differentiability point of $C_t$ throughout.

\subsection{Why the fractional inequality cannot simply be differentiated}\label{sec:derivative-obstruction}
The classical fractional P\'olya--Szeg\H{o} inequality asserts
\[
 F_t(u^*)\le F_t(u),\qquad 0<t<1.
\]
For a fixed smooth nonnegative $u$ whose rearrangement is also smooth, define the differentiable gap
\[
 \Delta_u(t)=F_t(u)-F_t(u^*)\ge0.
\]
The desired logarithmic inequality at a positive order $s$ would instead assert
\begin{equation}\label{eq:gapderivative}
 \Delta_u'(s)=\Lsp(u)-\Lsp(u^*)\ge0.
\end{equation}
Nonnegativity of a differentiable function does not imply nonnegativity of its derivative: for example, $e^{-t}>0$ has a negative derivative. Thus the fractional inequality does not imply \eqref{eq:gapderivative}.

More explicitly, for $h>0$ the inequality at $s+h$ gives only
\[
 \frac{F_{s+h}(u^*)-F_s(u^*)}{h}
 \le \frac{F_{s+h}(u)-F_s(u)}{h}
       +\frac{\Delta_u(s)}h.
\]
When $\Delta_u(s)>0$, the last term diverges as $h\downarrow0$, so the fractional comparison gives no ordering of the two derivatives. If exceptionally $\Delta_u(s)=0$ and the gap is differentiable near $s$, then $s$ is an interior minimum and $\Delta_u'(s)=0$. This special case supplies no derivative comparison for a general positive gap.

At the kernel level, the positive kernels $C_t r^{-N-tp}$ have derivative
\[
 \left.\partial_t(C_t r^{-N-tp})\right|_{t=s}
 =C_s(b_s-p\log r)r^{-N-sp}=C_s k_s(r),
\]
which is signed and eventually increasing in $r$. Their order derivatives therefore do not retain the positivity and radial decrease used in the classical proof. Differentiation under the integral may be entirely justified while differentiation of the inequality remains invalid.

There is also a direct obstruction within the fractional family. For $u_A(x)=u(x/A)$, dilation commutes with Schwarz rearrangement and
\begin{align*}
 \Delta_{u_A}(t)&=A^{N-tp}\Delta_u(t)\ge0,\\
 \Delta_{u_A}'(s)
 &=A^{N-sp}\bigl[\Delta_u'(s)-p\log A\,\Delta_u(s)\bigr].
\end{align*}
If $\Delta_u(s)>0$, the second expression is negative for all sufficiently large $A$, even though the fractional gap stays nonnegative at every order. Proposition~\ref{prop:dilation} supplies smooth examples. Proposition~\ref{prop:schwarzcore} below gives another smooth family with an exactly known rearrangement. Hence an unrestricted logarithmic inequality is actually false, rather than merely inaccessible to this limiting argument. The small-support result requires the separate kernel-ordering argument proved below. This discussion concerns differentiation at a fixed positive order; it does not rule out endpoint limits with a common leading term and separately justified first-order asymptotics, which concern different limiting functionals.

\subsection{Recent endpoint P\'olya--Szeg\H{o} and rigidity results}\label{sec:endpointcomparison}
For the order-zero logarithmic $p$-Laplacian of \cite{DJS}, use the renormalized form
\begin{align}
 \mathcal E_{L,p}(u)
 &=\frac{C_{N,p}}2\iint_{|x-y|<1}
 \frac{|u(x)-u(y)|^p}{|x-y|^N}\,dx\,dy\notag\\
 &\quad+\frac{C_{N,p}}2\iint_{|x-y|\ge1}
 \frac{|u(x)-u(y)|^p-|u(x)|^p-|u(y)|^p}{|x-y|^N}\,dx\,dy\notag\\
 &\quad+\rho_{N,p}\|u\|_p^p,\label{eq:endpointenergy}\\
 C_{N,p}&=\frac{p}{|\mathbb S^{N-1}|},\qquad
 \rho_{N,p}=2\log2-\gamma+\frac p2\psi(N/2).\notag
\end{align}
Set
\[
 X^p=\left\{u\in L^p(\R^N):
 \iint_{|x-y|<1}\frac{|u(x)-u(y)|^p}{|x-y|^N}\,dx\,dy<\infty\right\},
 \qquad
 \mathsf S_p(u)=\int_{\R^N}|u|^p\log^-|u|\,dx,
\]
where $\log^-t=\max\{0,-\log t\}$ and the integrand is zero at $u=0$.

\begin{theorem}[Endpoint comparison and rigidity: Hajaiej--Li--Zhan]\label{thm:endpointliterature}
Hajaiej, Li and Zhan prove the comparison and rigidity in
\cite[Theorem~1.1]{HLZ}, the spectral comparison and its equality case in
\cite[Theorem~1.2]{HLZ}, and the finiteness criterion in
\cite[Proposition~2.3 and Lemma~3.2]{HLZ}.
Let $N\ge1$, $1<p<\infty$, and let $u\in X^p$ be real-valued with $\mathsf S_p(u)<\infty$. Then $u^\#=(|u|)^*\in X^p$ and
\[
 \mathcal E_{L,p}(u^\#)\le\mathcal E_{L,p}(u).
\]
For $u\ne0$, equality holds exactly when
\[
 u=\sigma u^\#(\cdot-x_0)\quad\text{a.e.},
 \qquad \sigma\in\{-1,1\},\quad x_0\in\R^N.
\]
No support-size restriction is imposed. For $u\in X^p$, finiteness of $\mathcal E_{L,p}(u^\#)$ is equivalent to $\mathsf S_p(u)<\infty$. Furthermore, for a bounded nonempty open set $\Omega\subset\R^N$, define
\[
 X_0^p(\Omega)=\{u\in X^p:u=0\text{ a.e. on }\R^N\setminus\Omega\},
\]
and set
\[
 \lambda_{1,L,p}(\Omega)
 =\inf_{0\ne u\in X_0^p(\Omega)}
 \frac{\mathcal E_{L,p}(u)}{\|u\|_p^p}.
\]
For bounded nonempty open $\Omega$ and a ball $B$ of equal volume,
\[
 \lambda_{1,L,p}(B)\le\lambda_{1,L,p}(\Omega),
\]
with equality precisely when $\Omega$ agrees almost everywhere with a translate of $B$.
\end{theorem}
This theorem concerns the endpoint functional \eqref{eq:endpointenergy}; it is not a theorem about \eqref{eq:L}. Its equality statement requires a direct rigidity argument in addition to any first-order limiting comparison.

\begin{remark}[Recovering the endpoint inequality from our main result]\label{rem:main-to-endpoint}
Our fixed-positive-order Schwarz theorem recovers the inequality in Theorem~\ref{thm:endpointliterature} for regular compactly supported functions by letting $s\downarrow0$. Two facts make this possible.

First, for the normalization \eqref{eq:normalization},
\[
 \log R_{\mathrm{mon}}(s)=\frac1{ps}+O(1),
 \qquad R_{\mathrm{mon}}(s)\longrightarrow\infty.
\]
Hence every fixed compact support satisfies the diameter hypothesis of Theorem~\ref{thm:schwarz} for all sufficiently small $s>0$. If $u\ge0$ has compact support and $u,u^*\in W^{\sigma,p}(\R^N)$ for some $\sigma\in(0,1)$, that theorem gives
\begin{equation}\label{eq:main-endpoint-comparison}
 \mathcal L_{s,p}(u^*)\le\mathcal L_{s,p}(u)
 \qquad\text{for all sufficiently small }s>0.
\end{equation}
Indeed, positive fractional regularity places both functions in the logarithmic energy space at each smaller order. This class includes, in particular, smooth compactly supported functions with smooth rearrangement.

Second, the renormalized energy satisfies
\begin{equation}\label{eq:main-endpoint-limit}
 \lim_{s\downarrow0}p\mathcal L_{s,p}(f)=\mathcal E_{L,p}(f)
\end{equation}
for every compactly supported $f\in W^{\sigma,p}(\R^N)$. To verify this limit, put $M=\|f\|_p^p$ and write
\begin{align*}
 Q_s(f)
 &=\iint_{|x-y|<1}\frac{|f(x)-f(y)|^p}{|x-y|^{N+sp}}\,dx\,dy\\
 &\quad+\iint_{|x-y|\ge1}
 \frac{|f(x)-f(y)|^p-|f(x)|^p-|f(y)|^p}{|x-y|^{N+sp}}\,dx\,dy.
\end{align*}
Separating the mass term in the tail gives the exact identity
\[
 pF_s(f)=A(s)M+\frac{C_s}{2}Q_s(f),
 \qquad A(s)=\frac{C_s|\mathbb S^{N-1}|}{ps}.
\]
The explicit normalization extends $A$ smoothly to $s=0$, with
\[
 A(0)=1,\qquad A'(0)=\rho_{N,p},\qquad
 C_s\to0,\qquad C_s'\to C_{N,p}.
\]
Furthermore, $Q_s(f)\to Q_0(f)$ and $Q_s'(f)$ stays bounded as $s\downarrow0$. For the short-range term, the factor $|\log|x-y||$ introduced by differentiation is controlled by the extra fractional regularity, taking $s<\sigma/2$. For the renormalized tail, the numerator vanishes unless both function values are nonzero; compact support bounds all active distances. Differentiating the exact identity therefore yields
\[
 p\mathcal L_{s,p}(f)
 =A'(s)M+\frac{C_s'}2Q_s(f)+\frac{C_s}2Q_s'(f)
 \longrightarrow\rho_{N,p}M+\frac{C_{N,p}}2Q_0(f)
 =\mathcal E_{L,p}(f).
\]
Applying this limit to \eqref{eq:main-endpoint-comparison} proves
\[
 \mathcal E_{L,p}(u^*)\le\mathcal E_{L,p}(u).
\]
For signed functions in the same regular compact-support class, the modulus comparison applies eventually as well, since $R_{\mathrm{pos}}(s)\to\infty$, and gives the inequality with $u^\#=(|u|)^*$. This argument establishes the endpoint inequality in the stated class. It does not by itself extend it to all of $X^p$ or transfer the equality classification through the limit.

To explain the latter point, define the nonnegative deficits
\[
 \delta_s(u)=p\bigl[\mathcal L_{s,p}(u)-\mathcal L_{s,p}(u^*)\bigr],
 \qquad
 \delta_0(u)=\mathcal E_{L,p}(u)-\mathcal E_{L,p}(u^*).
\]
The preceding argument gives $\delta_s(u)\to \delta_0(u)$. For nonsymmetric $u$, the positive-order rigidity theorem gives $\delta_s(u)>0$, but the limit only implies $\delta_0(u)\ge0$. Strict positivity can disappear, as the scalar example $d(s)=s\to0$ shows. Endpoint equality therefore requires additional information, such as a positive limiting deficit for each nonsymmetric profile or a strict endpoint polarization argument. The same issue arises for the constant-sign condition in real-valued equality cases. Theorem~\ref{thm:endpointliterature} supplies these endpoint rigidity conclusions directly.
\end{remark}

\begin{remark}[Why an endpoint limit can work while an interior derivative does not]
For the normalization \eqref{eq:normalization}, and for smooth compactly supported functions for which the expansions are justified,
\begin{equation}\label{eq:endpointexpansion}
 pF_t(u)=\|u\|_p^p+t\mathcal E_{L,p}(u)+o(t)
 \qquad(t\downarrow0).
\end{equation}
Indeed, writing $M=\|u\|_p^p$ and subtracting $|u(x)|^p+|u(y)|^p$ on $|x-y|\ge1$, one obtains
\[
 \iint\frac{|u(x)-u(y)|^p}{|x-y|^{N+tp}}\,dx\,dy
 =\frac{2|\mathbb S^{N-1}|}{tp}M+Q_t(u),
\]
where $Q_t(u)\to Q_0(u)$ is the sum of the two renormalized integrals in \eqref{eq:endpointenergy}. Also
\[
 C_t=tC_{N,p}\bigl(1+\rho_{N,p}t+o(t)\bigr).
\]
These identities give \eqref{eq:endpointexpansion}. For compactly supported smooth functions, the renormalized tail has bounded interaction support and the short-range part is controlled by smoothness, justifying the limit.

If both $u$ and $u^\#$ admit this expansion, equimeasurability cancels their common leading mass term:
\[
 0\le\frac{p[F_t(u)-F_t(u^\#)]}{t}
 \longrightarrow\mathcal E_{L,p}(u)-\mathcal E_{L,p}(u^\#).
\]
Thus a first-order endpoint comparison can follow from the fractional inequality in this class. Extending it to the full space in Theorem~\ref{thm:endpointliterature} requires the endpoint integrability and approximation theory. At a fixed positive $s$, by contrast, the gap $F_s(u)-F_s(u^\#)$ generally does not vanish, so the cancellation needed for an ordered difference quotient is absent, as Section~\ref{sec:derivative-obstruction} explains. Finally, positivity of a gap at every $t>0$ does not imply strict positivity of its first-order limit: the gap could be $o(t)$. Consequently, equality classifications cannot be inferred solely from fractional strictness.
\end{remark}

\subsection{Energy space, supports, and zero extension}
\begin{definition}
Set
\[
 [u]_{s,\log,p}^p=\iint_{|x-y|<1}
 \frac{|u(x)-u(y)|^p\log(1/|x-y|)}{|x-y|^{N+sp}}\,dx\,dy,
\]
and
\[
 \X=\{u\in L^p(\R^N):[u]_{s,\log,p}<\infty\}.
\]
For an open set $\Omega$, define
\[
 X_0^{s,\log,p}(\Omega)=\{u\in\X:u=0\text{ a.e. on }\R^N\setminus\Omega\}.
\]
This definition uses zero extension and does not impose a separate density convention.
\end{definition}
The essential support is the complement of the largest open set on which $u=0$ almost everywhere. For a bounded nonempty essential support, its diameter is
\begin{equation}\label{eq:esssupportdiameter}
\diam(\esssupp u)
=\sup\bigl\{|x-y|:x,y\in\esssupp u\bigr\}.
\end{equation}
A support-diameter hypothesis means that there is a measurable set $S$, of the stated diameter, outside which $u$ vanishes. Changing representatives on null sets has no effect.

For a radial kernel $K$, abbreviate
\[
 E_K(u)=\iint_{\R^{2N}}|u(x)-u(y)|^pK(|x-y|)\,dx\,dy.
\]
When $K(|\cdot|)\in L^1(\R^N)$,
\begin{equation}\label{eq:tailbound}
 |E_K(u)|\le 2^p\|K(|\cdot|)\|_1\|u\|_p^p.
\end{equation}
Indeed, $|a-b|^p\le2^{p-1}(|a|^p+|b|^p)$. Splitting distances at a fixed number smaller than one shows that every $u\in\X$ has
\begin{equation}\label{eq:absolute}
 \iint_{\R^{2N}}\frac{|u(x)-u(y)|^p(1+|\log|x-y||)}{|x-y|^{N+sp}}\,dx\,dy<\infty.
\end{equation}
To verify \eqref{eq:absolute} explicitly, denote its left-hand side by $I(u)$ and split
\[
 I(u)=I_{<}(u)+I_{\ge}(u),
\]
where the two terms are the integrals over $|x-y|<e^{-1}$ and $|x-y|\ge e^{-1}$, respectively. For $0<r<e^{-1}$, we have $\log(1/r)>1$, and therefore
\[
 1+|\log r|=1+\log(1/r)\le 2\log(1/r).
\]
It follows directly from the defining seminorm that
\[
 \begin{aligned}
 I_{<}(u)
 &\le 2\iint_{|x-y|<e^{-1}}
 \frac{|u(x)-u(y)|^p\log(1/|x-y|)}{|x-y|^{N+sp}}\,dx\,dy\\
 &\le 2[u]_{s,\log,p}^{\,p}<\infty.
 \end{aligned}
\]
For the remaining distances, set
\[
 K(z)=\mathbf1_{\{|z|\ge e^{-1}\}}
 \frac{1+|\log|z||}{|z|^{N+sp}}.
\]
Writing $\sigma_{N-1}=|\mathbb S^{N-1}|$, polar coordinates give
\[
 \begin{aligned}
 \|K\|_{L^1(\R^N)}
 &=\sigma_{N-1}\left(
 \int_{e^{-1}}^1(1-\log r)r^{-1-sp}\,dr
 +\int_1^\infty(1+\log r)r^{-1-sp}\,dr\right)\\
 &=\sigma_{N-1}\left(
 \int_{e^{-1}}^1(1-\log r)r^{-1-sp}\,dr
 +\frac1{sp}+\frac1{(sp)^2}\right)<\infty.
 \end{aligned}
\]
Here the first integral is finite because its interval is bounded away from zero, and the tail integrals converge because $sp>0$. Using Tonelli's theorem and the elementary $p$-power inequality yields
\[
 \begin{aligned}
 I_{\ge}(u)
 &=\iint_{\R^{2N}}|u(x)-u(y)|^pK(x-y)\,dx\,dy\\
 &\le 2^{p-1}\iint_{\R^{2N}}
 (|u(x)|^p+|u(y)|^p)K(x-y)\,dx\,dy\\
 &=2^p\|K\|_{L^1(\R^N)}\|u\|_p^p<\infty.
 \end{aligned}
\]
Combining the two estimates proves \eqref{eq:absolute}, with the quantitative bound
\[
 I(u)\le 2[u]_{s,\log,p}^{\,p}
 +2^p\|K\|_{L^1(\R^N)}\|u\|_p^p.
\]
The split is taken strictly below $1$ because $\log(1/r)$ tends to zero as $r\uparrow1$ and hence cannot control the added constant $1$ uniformly on $(0,1)$.

\subsection{Polarization}
Let $H$ be a half-space and let $\sigma_H$ denote reflection in its boundary. Boundary values are irrelevant. For $u\ge0$, its polarization is
\begin{equation}\label{eq:polarization}
 u^H(x)=\begin{cases}
 \max\{u(x),u(\sigma_Hx)\},&x\in H,\\
 \min\{u(x),u(\sigma_Hx)\},&x\notin H.
 \end{cases}
\end{equation}
Polarization preserves distributions, $L^p$ norms, and integrals of integrable functions of the values. If $x,y\in H$, then
\begin{equation}\label{eq:geometry}
 |x-y|\le |x-\sigma_Hy|,
\end{equation}
strictly for interior points. Reflection itself preserves every radial-kernel energy. For the definition of polarization (also called two-point rearrangement), its distribution-preserving property, and its $L^p$ continuity, see Hajaiej \cite[Sections~1--2 and Corollary~3.4]{HApprox}; further background is available in Brock and Solynin \cite{BS}.

\subsection{Steiner and Schwarz symmetrization}
Fix a unit vector $e$. Write $x=z+te$, with $z\in e^\perp$. The Steiner symmetrization $u^{\star_e}$ is obtained by symmetrically decreasingly rearranging $t\mapsto u(z+te)$ on $\R$, for almost every $z$. For every $\alpha>0$, the slice of $\{u^{\star_e}>\alpha\}$ at $z$ is the centered interval of length $|\{t:u(z+te)>\alpha\}|$. The Schwarz rearrangement $u^*$ has centered ball superlevel sets of the same measure as those of $u$.

Both operations preserve $L^p$ norms and global distributions. Steiner symmetrization also preserves the distribution in each slice. If $A:[0,\infty)\to[0,\infty)$ is continuous and strictly increasing with $A(0)=0$, then
\begin{equation}\label{eq:commute}
 (A(u))^H=A(u^H),\qquad
 (A(u))^{\star_e}=A(u^{\star_e}),\qquad (A(u))^*=A(u^*).
\end{equation}
For definitions through superlevel sets and the basic properties of Steiner and Schwarz symmetrization, see \cite[Section~2]{VS06}. Further background on equimeasurability and rearrangement inequalities is provided by Baernstein \cite{Baernstein}; for symmetric decreasing rearrangement, see also \cite[Chapter~3]{LL}.

\subsection{Classical results used in the proofs}\label{sec:classicalinputs}
We state the exact forms of the external geometric and approximation results
that are used. No equality case of the set Riesz inequality is assumed in our
proofs. The strict positive-kernel energy comparisons needed later are proved
separately in Lemma~\ref{lem:positivekernel}.

\begin{theorem}[Set form of the Riesz rearrangement inequality]
\label{thm:rieszinput}
Let $N\ge1$, let $E,F\subset\R^N$ be Lebesgue measurable sets of finite
measure, and let $r>0$. Denote by $E^*,F^*$ the centered balls of measures
$|E|,|F|$, respectively. Then
\begin{equation}\label{eq:rieszinput}
 \int_E\int_F\mathbf1_{\{|x-y|<r\}}\,dy\,dx
 \le\int_{E^*}\int_{F^*}\mathbf1_{\{|x-y|<r\}}\,dy\,dx.
\end{equation}
Both integrals are finite. This follows from Burchard and Hajaiej
\cite[Theorem~2]{BH} with the product integrand, the functions
$\mathbf1_E,\mathbf1_F$, and the kernel
$K(t)=\mathbf1_{\{t<r\}}$.
\end{theorem}

\begin{theorem}[Supermodular Riesz inequality on the line and circle]
\label{thm:supermodularinput}
Let $X$ be $\R$ with Lebesgue measure and Euclidean distance, or the circle
$\mathbb T_d$ of circumference $d>0$ with arc-length measure and geodesic
distance. Rearrangement is about one fixed origin in $X$. Let
$F:[0,\infty)^2\to[0,\infty)$ be Borel measurable, vanish on the coordinate
axes, and be supermodular, meaning that for all $a,b,h,\ell\ge0$,
\[
 F(a+h,b+\ell)+F(a,b)\ge F(a+h,b)+F(a,b+\ell).
\]
Let $K:[0,\infty)\to[0,\infty]$ be measurable and nonincreasing. If
$f,g:X\to[0,\infty)$ are measurable and, on $\R$, all their positive
superlevel sets have finite measure, then
\[
 \begin{aligned}
 &\iint_{X^2}F(f(x),g(y))K(\operatorname{dist}_X(x,y))\,dx\,dy\\
 &\quad\le\iint_{X^2}F(f^*(x),g^*(y))
 K(\operatorname{dist}_X(x,y))\,dx\,dy.
 \end{aligned}
\]
The integrals are nonnegative extended integrals. This is the $m=2$
specialization of Burchard and Hajaiej \cite[Theorem~2]{BH}; its negative-part condition is
automatic because $F\ge0$. The theorem on the unit circle gives the
stated circle version by scaling distance and arc-length measure.
No equality classification from that theorem is used here.
\end{theorem}

\begin{theorem}[Euclidean isodiametric inequality]
\label{thm:isodiametricinput}
If $E\subset\R^N$ is Lebesgue measurable and has finite diameter $d$, then
\[
 |E|\le\omega_N(d/2)^N.
\]
For $d>0$, equality holds if and only if $E$ agrees almost everywhere with a
ball of radius $d/2$. We use only the volume bound and its consequence
$\diam(E^*)\le d$. This is the isodiametric inequality in
Maggi, Ponsiglione and Pratelli \cite[inequality~(1.1)]{MPP}.
\end{theorem}

\begin{theorem}[Approximation by successive polarizations]
\label{thm:polarapprox}
Let $N\ge1$, $1<p<\infty$, and let $u\in L^p(\R^N)$ be nonnegative.
Choose either $Tu=u^*$, or $Tu=u^{\star_e}$ for one fixed unit vector $e$.
There is a sequence of half-spaces $(H_n)_{n\ge1}$ such that the successive
iterates
\[
 u_0=u,\qquad u_n=(u_{n-1})^{H_n}
\]
satisfy $\|u_n-Tu\|_p\to0$. For Schwarz rearrangement the half-spaces may be
chosen to contain the origin. For Steiner symmetrization, each boundary is
perpendicular to $e$ and its selected half-space contains $e^\perp$.
In either case $T(u_n)=Tu$ for every $n$. If $u=0$ almost everywhere outside
a centered ball $B_M(0)$, every $u_n$ also vanishes outside that ball.
The Schwarz approximation is Hajaiej's constructive theorem
\cite[unnumbered Theorem, p.~410]{HApprox}. The Steiner approximation
follows from its one-dimensional version by the slice argument below.
\end{theorem}
\begin{proof}
For Schwarz rearrangement, fix a dense sequence $(G_j)_{j\ge1}$ of
closed half-spaces containing the origin. Hajaiej
\cite[unnumbered Theorem, p.~410]{HApprox} proves that the block iterates
\[
 v_0=u,\qquad v_{m+1}=v_m^{G_1\cdots G_{m+1}}
\]
converge to $u^*$ strongly in $L^p$. Concatenate the finite blocks
$(G_1),(G_1,G_2),(G_1,G_2,G_3),\ldots$ into a sequence of individual
polarizations. Its block endpoints converge to $u^*$. At every
individual step the distance to $u^*$ is nonincreasing, since $u^*$ is
fixed by these polarizations and polarization is nonexpansive in $L^p$
by Hajaiej \cite[Corollary~3.4]{HApprox}. Thus the full sequence converges.

For Steiner symmetrization in direction $e$, write $x=z+te$, with
$z\in e^\perp$. Fix one dense sequence $(I_j)$ of closed half-lines in
$\R$ containing zero, and apply the same block construction with
\[
 G_j=\{z+te:z\in e^\perp,\ t\in I_j\}.
\]
For almost every $z$, Fubini's theorem gives
$u(z+\cdot e)\in L^p(\R)$. The one-dimensional instance of Hajaiej's
theorem, with this same sequence $(I_j)$, gives convergence on that
slice to its symmetric decreasing rearrangement. The slice error is
bounded by
\[
 \|v_m(z+\cdot e)-u^{\star_e}(z+\cdot e)\|_{L^p(\R)}^p
 \le2^p\|u(z+\cdot e)\|_{L^p(\R)}^p.
\]
The right-hand side is integrable in $z$. Dominated convergence
therefore gives strong $L^p(\R^N)$ convergence of the block endpoints
to $u^{\star_e}$. Nonexpansiveness again extends it to every individual
polarization in the concatenated sequence.

Equimeasurability gives $T(u_n)=Tu$; in the Steiner case this holds
slice by slice. The selected half-spaces contain the origin in the
Schwarz case and $e^\perp$ in the Steiner case. Their polarizations
preserve any enclosing centered ball, as follows directly from the
reflection distance inequality. These observations prove the remaining
assertions.
\end{proof}

For comparison with local energies, the classical local
P\'olya--Szeg\H{o} inequality in the form mentioned in the introduction is:
if $1<p<\infty$ and $u\ge0$ belongs to $W^{1,p}(\R^N)$, then
$u^*\in W^{1,p}(\R^N)$ and
\[
 \int_{\R^N}|\nabla u^*|^p\,dx\le\int_{\R^N}|\nabla u|^p\,dx.
\]
This is \cite[Corollary~6.3]{VS06} with the ambient domain $\R^N$ and
Schwarz symmetrization. It is recalled for context, not used to deduce a
comparison for the signed nonlocal energy.

\section{Analytic and geometric lemmas}
\begin{lemma}[Kernel geometry]\label{lem:kernel}
For $N\ge1$, $1<p<\infty$, $0<s<1$ and the radii \eqref{eq:radii}, the kernel $k_s$ is positive on $(0,\Rpos)$, negative on $(\Rpos,\infty)$, strictly decreasing on $(0,\Rmon]$, and strictly increasing on $[\Rmon,\infty)$.
\end{lemma}
\begin{proof}
The sign assertion follows from its definition. Differentiation gives
\[
 k_s'(r)=r^{-q-1}(pq\log r-qb_s-p),
\]
whose bracket vanishes exactly at $r=\Rmon$.
\end{proof}

\begin{lemma}[Signed-tail continuity and lower semicontinuity]\label{lem:lsc}
Extend $\Lsp$ to $L^p$ by assigning $+\infty$ outside $\X$. It is lower semicontinuous for strong $L^p$ convergence.
\end{lemma}
\begin{proof}
Choose $0<\delta<\min\{1,\Rpos\}$. Set $P=k_s\mathbf1_{(0,\delta)}\ge0$ and $G=k_s\mathbf1_{[\delta,\infty)}$, so $G(|\cdot|)\in L^1$. The map $v\mapsto v(x)-v(y)$ is bounded from $L^p(\R^N)$ to $L^p(\R^{2N},|G(|x-y|)|\,dx\,dy)$, by \eqref{eq:tailbound}. Hence $E_G$ is continuous in $L^p$: the $p$th powers converge in $L^1$ of this weighted measure. Fatou's lemma gives lower semicontinuity of $E_P$, after taking an almost-everywhere convergent subsequence realizing the lower limit. Finite $P$-energy is equivalent to membership in $\X$ for $L^p$ functions, since $P(r)$ is comparable to $r^{-q}\log(1/r)$ sufficiently near zero. The conclusion follows.
\end{proof}

\begin{lemma}[Polarization does not increase support diameter]\label{lem:diameter}
If $u\ge0$ vanishes outside $S$, then $u^H$ vanishes outside a set of diameter at most $\diam S$.
\end{lemma}
\begin{proof}
Let $S^H$ be the polarization of $S$: in $H$ it is $(S\cup\sigma_HS)\cap H$, and outside $H$ it is $(S\cap\sigma_HS)\setminus H$. The positivity set of $u^H$ is contained in $S^H$ up to null sets. If two points of $S^H$ lie in $H$, choose for each an original point in $S$ from its reflection pair. Their distance is unchanged if the choices are on the same side, and can only be larger than the new distance if they are on opposite sides, by \eqref{eq:geometry}. If one new point is outside $H$, both points in its reflection pair belong to $S$, and one can choose the representative matching the side chosen for the other point; the distance is then unchanged. If both new points are outside $H$, both already belong to $S$. Thus every new distance is bounded by an original distance.
\end{proof}

\begin{lemma}[Two-point defect]\label{lem:Q}
For real $a,b,c,d$, put $a_+=\max(a,b)$, $a_-=\min(a,b)$, and similarly $c_\pm$. Define
\[
 Q_p(a,b,c,d)=|a-c|^p+|b-d|^p-|a_+-c_+|^p-|a_--c_-|^p.
\]
Then $Q_p\ge0$, and $Q_p>0$ if and only if $(a-b)(c-d)<0$.
\end{lemma}
\begin{proof}
For equally ordered pairs $Q_p=0$. If $a>b$ and $d>c$, direct integration gives
\begin{equation}\label{eq:strictQ}
 Q_p=p(p-1)\int_b^a\int_c^d|\alpha-\beta|^{p-2}\,d\beta\,d\alpha>0.
\end{equation}
The diagonal singularity is locally integrable when $1<p<2$. The reverse ordering follows by interchanging the two pairs.
\end{proof}

\begin{lemma}[Auxiliary positive-kernel comparison and rigidity]
\label{lem:positivekernel}
Let $N\ge1$, $1<p<\infty$, and let $K:(0,\infty)\to[0,\infty)$ be finite,
measurable and nonincreasing, with
$K(|\cdot|)\mathbf1_{\{|\cdot|>\delta\}}\in L^1(\R^N)$ for every
$\delta>0$. Suppose that $u\ge0$ belongs to $L^p(\R^N)$, has bounded
support and $E_K(u)<\infty$. Then for every half-space $H$ and every unit
vector $e$,
\[
 E_K(u^H)\le E_K(u),\qquad
 E_K(u^{\star_e})\le E_K(u),\qquad E_K(u^*)\le E_K(u),
\]
and all three transformed energies are finite. There is no support-diameter
restriction. If $K$ is strictly decreasing on $(0,\infty)$ and $u\ne0$,
equality holds in the three comparisons, respectively, exactly when
\[
 u^H=u\ \text{or}\ u^H=u\circ\sigma_H;\qquad
 u=u^{\star_e}(\cdot-ae)\ (a\in\R);\qquad
 u=u^*(\cdot-b)\ (b\in\R^N),
\]
where all identities are almost everywhere.
\end{lemma}
\begin{proof}
Cap the kernel at the origin by $K_\varepsilon(r)=K(\max\{r,\varepsilon\})$.
It is integrable on $\R^N$. With $r=|x-y|$,
$\rho=|x-\sigma_Hy|$ and the four values from Lemma~\ref{lem:Q}, splitting
the double integral into its four reflected quadrants gives
\[
 E_{K_\varepsilon}(u)-E_{K_\varepsilon}(u^H)
 =\int_{H\times H}[K_\varepsilon(r)-K_\varepsilon(\rho)]Q_p\,dx\,dy\ge0.
\]
Here $r\le\rho$ by \eqref{eq:geometry} and $Q_p\ge0$ by
Lemma~\ref{lem:Q}. Monotone convergence as $\varepsilon\downarrow0$
proves the polarization comparison and finiteness. The kernel differences
also increase to $K(r)-K(\rho)$, yielding the corresponding exact defect
identity. If $K$ is strictly decreasing and both ordering sets
\[
 A_H=\{x\in H:u(x)>u(\sigma_Hx)\},\quad
 B_H=\{x\in H:u(x)<u(\sigma_Hx)\}
\]
have positive measure, their product has $Q_p>0$ and $r<\rho$ almost
everywhere. The defect is then strictly positive. Equality is therefore
precisely the global ordering alternative stated in the lemma.

For $T$ equal to Steiner or Schwarz symmetrization, take the successive
polarizations from Theorem~\ref{thm:polarapprox}. Their energies decrease;
Fatou's lemma along an almost-everywhere convergent subsequence gives
$E_K(Tu)\le\liminf_nE_K(u_n)\le E_K(u)$. If the endpoint energies are equal,
every polarization step has equality and hence is an identity or a reflection
of the entire function. For Steiner symmetrization, the compositions have
form $u_n(z+te)=u(z+(\epsilon_nt+c_n)e)$, with $\epsilon_n\in\{-1,1\}$.
Their supports stay in a fixed centered ball; since $u\ne0$, $(c_n)$ is bounded.
A subsequence has $\epsilon_n=\epsilon$ and $c_n\to c$. Strong continuity
of translations in $L^p$ and evenness of $Tu$ in $t$ give
$u=Tu(\cdot-ce)$. In the Schwarz case the compositions have form
$u_n(x)=u(O_nx+c_n)$, $O_n\in O(N)$. The same enclosing-ball argument
bounds $c_n$; compactness of $O(N)$ supplies convergent subsequences.
Strong continuity of Euclidean motions and radiality of $Tu$ give
$u=Tu(\cdot-b)$ for one $b\in\R^N$. Reflection and translation invariance
prove the converses.
\end{proof}

\begin{corollary}[Fractional P\'olya--Szeg\H{o} comparison used in this paper]
\label{cor:fractionalinput}
Let $N\ge1$, $1<p<\infty$ and $0<t<1$. Suppose that $u\ge0$ belongs to
$L^p(\R^N)$, has bounded support, and
\[
 [u]_{t,p}^p:=\iint_{\R^{2N}}
 \frac{|u(x)-u(y)|^p}{|x-y|^{N+tp}}\,dx\,dy<\infty.
\]
Then the fractional energies of $u^H,u^{\star_e},u^*$ are finite and do not
exceed $[u]_{t,p}^p$. For $u\ne0$, equality has exactly the three configurations
listed in Lemma~\ref{lem:positivekernel}. No restriction $tp<N$ is needed
in this bounded-support $L^p$ formulation.
\end{corollary}
\begin{proof}
Apply Lemma~\ref{lem:positivekernel} to the strictly decreasing kernel
$K(r)=r^{-N-tp}$. Its integral outside every positive radius is finite
because $tp>0$. This proves both the comparisons and the equality cases.
\end{proof}
\noindent
For bibliographic precision, \cite[Theorem~A.1]{FS} states a corresponding
comparison and equality theorem in the homogeneous fractional Sobolev space
with the range $1\le p<N/t$. Corollary~\ref{cor:fractionalinput}, whose proof
is given here, is the bounded-support $L^p$ result invoked in our arguments
for all $p>1$; we do not apply the homogeneous-space theorem outside its
stated range.

\section{Polarization inequality and complete equality cases}
The comparison radius is set by decrease of the kernel, rather than its
sign. The following theorem retains the full negative tail and gives
the complete global-ordering alternative for equality. Its diameter
bound is proved sharp in Proposition~\ref{prop:sharpdiameter}.

\begin{theorem}[Signed-kernel polarization comparison and rigidity]\label{thm:polar}
Fix $N\ge1$, $1<p<\infty$, and $0<s<1$ at a differentiability point of $C_t$. Let $u\ge0$ belong to $\X$ and vanish almost everywhere outside a measurable set $S$ with $\diam S\le\Rmon$. For every half-space $H$, $u^H\in\X$ and
\begin{equation}\label{eq:polarineq}
 \Lsp(u^H)\le\Lsp(u).
\end{equation}
Equality holds if and only if $u^H=u$ a.e. or $u^H=u\circ\sigma_H$ a.e. Equivalently, one of the two orderings
\[
 u(x)\ge u(\sigma_Hx)\ \text{for a.e. }x\in H,
 \qquad\text{or}\qquad
 u(x)\le u(\sigma_Hx)\ \text{for a.e. }x\in H
\]
holds throughout $H$.
\end{theorem}
\begin{proof}
Abbreviate $\sigma=\sigma_H$. For $x,y\in H$, write
\[
 a=u(x),\quad b=u(\sigma x),\quad c=u(y),\quad d=u(\sigma y),
 \qquad r=|x-y|,\quad \rho=|x-\sigma y|.
\]
For an integrable radial kernel $K$, splitting $\R^{2N}$ into its four reflected quadrants gives
\begin{equation}\label{eq:defectK}
 E_K(u)-E_K(u^H)=\int_{H\times H}[K(r)-K(\rho)]Q_p(a,b,c,d)\,dx\,dy.
\end{equation}
Indeed, the same-side terms have coefficient $K(r)$, the cross-side terms coefficient $K(\rho)$, and the sum of the four value differences is unchanged by sorting both pairs.

If $Q_p>0$ and $a>b,c<d$, nonnegativity gives $a,d>0$, whence $x,\sigma y\in S$ up to null sets. In the reverse ordering, $b,c>0$, so $\sigma x,y\in S$. In either case
\begin{equation}\label{eq:activedistance}
 r\le\rho\le\diam S\le\Rmon.
\end{equation}
Thus decrease of the kernel is needed only on active distances.

For sufficiently small $\varepsilon>0$, define $k_\varepsilon(r)=k_s(\max\{r,\varepsilon\})$. This is integrable as a radial function on $\R^N$ and nonincreasing on $(0,\Rmon]$. Equations \eqref{eq:defectK}--\eqref{eq:activedistance} yield $E_{k_\varepsilon}(u^H)\le E_{k_\varepsilon}(u)$.

Fix $\delta$ as in Lemma~\ref{lem:lsc}, and take $\varepsilon<\delta$. Then $k_\varepsilon=P_\varepsilon+G$, where $P_\varepsilon=k_\varepsilon\mathbf1_{(0,\delta)}\ge0$ increases to $P$ as $\varepsilon\downarrow0$, and $G$ is independent of $\varepsilon$. Since $\|u^H\|_p=\|u\|_p$, \eqref{eq:tailbound} gives a uniform bound for $E_{P_\varepsilon}(u^H)$. Monotone convergence implies finite $P$-energy of $u^H$, hence $u^H\in\X$. Taking limits proves \eqref{eq:polarineq}. On the set where $Q_p>0$, the differences $k_\varepsilon(r)-k_\varepsilon(\rho)$ also increase to $k_s(r)-k_s(\rho)$, by \eqref{eq:activedistance}. A further monotone-convergence argument gives the exact identity
\begin{equation}\label{eq:defect}
 \Lsp(u)-\Lsp(u^H)=\frac{C_s}{2p}\int_{H\times H}
 [k_s(r)-k_s(\rho)]Q_p(a,b,c,d)\,dx\,dy.
\end{equation}

For equality, let $A=\{x\in H:u(x)>u(\sigma x)\}$ and $B=\{x\in H:u(x)<u(\sigma x)\}$. If both have positive measure, then on $A\times B$ almost everywhere, $Q_p>0$ and $r<\rho\le\Rmon$. Strict decrease makes \eqref{eq:defect} strictly positive. Hence equality requires $|A|=0$ or $|B|=0$. These alternatives give $u^H=u\circ\sigma$ and $u^H=u$, respectively. Conversely, both identities imply equality, the former by reflection invariance.
\end{proof}

\begin{remark}
Equality for a single polarization does not, by itself, imply radial symmetry or Steiner symmetry. In particular, a globally reversed ordering merely reflects the function and leaves its energy unchanged. Neither strict positivity of $k_s$ nor the stronger bound $\diam S\le\Rpos$ is needed in Theorem~\ref{thm:polar}.
\end{remark}

\section{Steiner symmetrization and common-center rigidity}
\begin{theorem}[Steiner comparison and common-translation rigidity]\label{thm:steiner}
Fix $N\ge1$, $1<p<\infty$ and $0<s<1$ at a differentiability point of
$C_t$. Let $u\ge0$ belong to $\X$ and vanish almost everywhere outside a
measurable set $S$ with $\diam S\le\Rmon$, and let $e$ be a unit vector.
Then $u^{\star_e}\in\X$ and
\begin{equation}\label{eq:steiner}
 \Lsp(u^{\star_e})\le\Lsp(u).
\end{equation}
For nonzero $u$, equality holds if and only if there is a single $a\in\R$ such that
\begin{equation}\label{eq:steinerequality}
 u(z+te)=u^{\star_e}(z+(t-a)e)\quad\text{for a.e. }(z,t)\in e^\perp\times\R.
\end{equation}
\end{theorem}
\begin{proof}
Theorem~\ref{thm:polarapprox}, applied to the nonnegative $L^p$ function $u$
with $T=\star_e$, supplies successive polarizations $u_0=u$,
$u_n=(u_{n-1})^{H_n}$ converging strongly in $L^p$ to $u^{\star_e}$,
with reflecting hyperplanes perpendicular to $e$. Lemma~\ref{lem:diameter} keeps every support diameter at most $\diam S$. Therefore Theorem~\ref{thm:polar} gives
\[
 \Lsp(u_{n+1})\le\Lsp(u_n)\le\Lsp(u).
\]
Lemma~\ref{lem:lsc} yields
\begin{equation}\label{eq:chain}
 \Lsp(u^{\star_e})\le\liminf_n\Lsp(u_n)\le\Lsp(u),
\end{equation}
including membership of the limit in $\X$.

Assume $u\ne0$ and equality in \eqref{eq:steiner}. The decreasing energies in \eqref{eq:chain} must all equal $\Lsp(u)$. The equality part of Theorem~\ref{thm:polar} shows that every step either does nothing or reflects the whole current function. Since all reflections act only in the $e$-coordinate, compositions have the form
\begin{equation}\label{eq:orbit}
 u_n(z+te)=u(z+(\epsilon_n t+c_n)e),\qquad \epsilon_n\in\{-1,1\}.
\end{equation}
The original support is bounded. If a subsequence had $|c_n|\to\infty$, its supports would escape every bounded slab in $t$. Strong $L^p$ convergence would then force the limit to vanish on every bounded slab, contrary to $\|u^{\star_e}\|_p=\|u\|_p>0$. Thus $c_n$ is bounded. Along a subsequence, $\epsilon_n=\epsilon$ and $c_n\to c$. Continuity of translations in $L^p$ gives
\[
 u^{\star_e}(z+te)=u(z+(\epsilon t+c)e).
\]
Since $u^{\star_e}$ is even in $t$, solving this relation yields \eqref{eq:steinerequality} with $a=c$. Conversely, \eqref{eq:steinerequality} gives equality by translation invariance.
\end{proof}

\begin{remark}[Meaning of the common center]
Equality means that, for almost every slice and almost every positive level, the nonempty superlevel set in that slice is an interval centered at the same number $a$. A center depending on $z$ or on the level is not allowed, except for immaterial empty slices or null-set changes. No assumption on the measure of a critical set or on flat portions of the symmetrized function is imposed. In dimension one this theorem is the corresponding Schwarz statement.
\end{remark}

\section{P\'olya--Szeg\H{o} inequality for Schwarz symmetrization}
For Schwarz rearrangement, the sufficient diameter bound $\Rmon$ holds in every dimension. We first prove the corresponding comparison and equality cases, then establish the optimal bound $2\Rmon$ in dimension one. In higher dimensions, we give additional sufficient bounds for specified classes of superlevel-set geometries; these do not identify the unrestricted sharp diameter.

\subsection{Definition and precise statement}
For $u\ge0$, let $\mu_u(\alpha)=|\{u>\alpha\}|$ be its distribution function and denote the volume of the unit ball by $\omega_N$. For the distribution function, equimeasurability, and the associated symmetric decreasing rearrangement, see \cite[Chapter~3]{LL}. The Schwarz symmetrization $u^*$ is characterized by
\[
 \{u^*>\alpha\}=B_{(\mu_u(\alpha)/\omega_N)^{1/N}}(0)
 \quad\text{up to null sets, for every }\alpha>0.
\]
Thus $u^*$ is radial, nonincreasing in the radius, and equimeasurable with $u$. Unlike Steiner symmetrization, which rearranges one-dimensional slices, Schwarz symmetrization rearranges each whole superlevel set into a ball. For this superlevel-set definition of Schwarz symmetrization, see \cite[Section~1]{VS09}; for the definition of Steiner symmetrization and its relation to Schwarz symmetrization, see \cite[Section~2]{VS06}. A systematic treatment of both operations and their basic properties is given in \cite{Baernstein}.

\begin{theorem}[Schwarz P\'olya--Szeg\H{o} inequality and rigidity]\label{thm:schwarz}
Fix $N\ge1$, $1<p<\infty$, and $0<s<1$ {at a differentiability point of the chosen positive normalization}. Let $u\ge0$ belong to $\X$, and suppose that $u=0$ almost everywhere outside a measurable set $S$ satisfying
\[
 \diam S\le\Rmon=\exp\!\left(\frac{b_s}{p}+\frac{1}{N+sp}\right).
\]
Then $u^*\in\X$ and
\begin{equation}\label{eq:schwarzfull}
 \begin{aligned}
 &\iint_{\R^{2N}}|u^*(x)-u^*(y)|^p k_s(|x-y|)\,dx\,dy\\
 &\qquad\le
 \iint_{\R^{2N}}|u(x)-u(y)|^p k_s(|x-y|)\,dx\,dy.
 \end{aligned}
\end{equation}
Equivalently, $\Lsp(u^*)\le\Lsp(u)$. If $u\ne0$, equality holds if and only if there is $a\in\R^N$ with
\begin{equation}\label{eq:schwarzequality}
 u(x)=u^*(x-a)\quad\text{for a.e. }x\in\R^N.
\end{equation}
In terms of superlevel sets, all positive-measure superlevel sets are balls with one common center, up to null sets and almost every level.
\end{theorem}
\begin{proof}
Apply Theorem~\ref{thm:polarapprox} to $u\ge0$ with $T=*$: it supplies
successive polarizations $u_0=u$, $u_n=(u_{n-1})^{H_n}$ converging strongly
in $L^p(\R^N)$ to $u^*$, with every $H_n$ containing the origin. Lemma~\ref{lem:diameter} shows that all positivity sets have diameter at most $\diam S$. Hence Theorem~\ref{thm:polar} applies at every step and gives
\begin{equation}\label{eq:schwarzchain}
 \Lsp(u_{n+1})\le\Lsp(u_n)\le\Lsp(u).
\end{equation}
By Lemma~\ref{lem:lsc},
\[
 \Lsp(u^*)\le\liminf_{n\to\infty}\Lsp(u_n)\le\Lsp(u).
\]
This also implies $u^*\in\X$, and proves \eqref{eq:schwarzfull}. In particular, the signed tail has been retained throughout; the argument does not replace $k_s$ by its positive part.

Now suppose $u\ne0$ and the endpoint energies are equal. The decreasing energies in \eqref{eq:schwarzchain} must all equal $\Lsp(u)$. Equality at each polarization step, by Theorem~\ref{thm:polar}, means that the current function is either unchanged or entirely reflected. Composing these operations gives
\[
 u_n(x)=u(O_nx+c_n),\qquad O_n\in O(N),\quad c_n\in\R^N.
\]
Choose $M$ such that $u$ is supported in $B_M(0)$. Polarizations toward half-spaces containing the origin preserve this support bound. Since $u_n\ne0$, a positive-measure set of points $x\in B_M(0)$ has $O_nx+c_n\in B_M(0)$; therefore $|c_n|\le2M$. Compactness of $O(N)$ and boundedness of $(c_n)$ yield a subsequence for which $O_n\to O$ and $c_n\to c$. The action of Euclidean motions is strongly continuous on $L^p$, as follows by approximating with continuous compactly supported functions. Thus
\[
 u^*(x)=u(Ox+c)\quad\text{a.e.}
\]
Solving for $u$ gives $u(y)=u^*(O^{-1}(y-c))=u^*(y-c)$, because $u^*$ is radial. This is \eqref{eq:schwarzequality}. Conversely, a translate has exactly the same energy by translation invariance. The superlevel-set formulation follows from
$u(x)=\int_0^\infty\mathbf1_{\{u>t\}}(x)\,dt$ for nonnegative measurable
$u$, together with the defining superlevel sets of $u^*$; the identity is
Tonelli's layer-cake formula, recalled in \cite[Chapter~1]{LL}.
\end{proof}

\begin{remark}
The zero function gives equality trivially. For nonzero functions, the equality statement does not require excluding plateaus or critical sets. The diameter condition $\diam S\le\Rmon$ is sufficient in every dimension. It is not claimed to be the optimal universal Schwarz condition. Proposition~\ref{prop:explicitschwarz} constructs a smooth counterexample at a larger diameter, but does not establish failure at every diameter above $\Rmon$.
\end{remark}

\subsection{The sharp one-dimensional Schwarz diameter}
\begin{theorem}[Schwarz comparison up to twice the monotonicity radius in dimension one]
\label{thm:schwarzonedimensional}
Let $N=1$, $1<p<\infty$, $0<s<1$, and let $C_t>0$ be continuously
differentiable near $s$. Set $R=\Rmon$ as in \eqref{eq:radii}. If
$u\ge0$ belongs to $X^{s,\log,p}(\R)$ and
\[
 \diam(\esssupp u)\le2R,
\]
then $u^*\in X^{s,\log,p}(\R)$ and
\begin{equation}\label{eq:schwarzonedimensional}
 \Lsp(u^*)\le\Lsp(u).
\end{equation}
The bound $2R$ is the optimal uniform support-diameter bound: for every
$\ell>2R$ there exists a nonnegative $u\in C_c^\infty(\R)$, with
$u^*\in C_c^\infty(\R)$ and $\diam(\esssupp u)<\ell$, such that
$\Lsp(u^*)>\Lsp(u)$. Consequently the optimal one-dimensional Schwarz
threshold, defined in Corollary~\ref{cor:schwarzthreshold}, is $2R$.
The same threshold applies to one-dimensional Steiner symmetrization,
which coincides with Schwarz rearrangement. For $u\not\equiv0$, equality in \eqref{eq:schwarzonedimensional}, including at diameter $2\Rmon$, holds if and only if $u=u^*(\cdot-a)$ almost everywhere for one $a\in\R$. No restriction $sp<1$ is imposed.
\end{theorem}
\begin{proof}
The zero function is immediate. Otherwise, put $d=\diam(\esssupp u)>0$
and $m=d/2\le R$. By translation we may assume that $u$ vanishes outside
$I=(-m,m)$. Write $k=k_s$ and $q=1+sp$.

\emph{Step 1: a monotone correction kernel.}
For $0<t<m\le R$,
\begin{equation}\label{eq:circlepairedderivative}
 k'(m-t)+k'(m+t)<0.
\end{equation}
Indeed, $k'(r)=pq\,r^{-q-1}\log(r/R)$. On putting $z=t/m$, the sum
divided by $pq\,m^{-q-1}$ is
\[
 \begin{aligned}
 &\log(m/R)\big[(1-z)^{-q-1}+(1+z)^{-q-1}\big]\\
 &\qquad+\frac{\log(1-z)}{(1-z)^{q+1}}
       +\frac{\log(1+z)}{(1+z)^{q+1}}<0.
 \end{aligned}
\]
The first term is nonpositive, while the remaining sum is negative since
$-\log(1-z)>\log(1+z)>0$ and $(1-z)^{-q-1}>(1+z)^{-q-1}$.
Define
\begin{equation}\label{eq:circlecorrection}
 B(r)=\begin{cases}
 0,&0\le r\le m,\\
 k(d-r)-k(r),&m<r<d.
 \end{cases}
\end{equation}
It is nonnegative and nondecreasing, because $B(m)=0$ and
$B'(r)=-k'(d-r)-k'(r)>0$ on $(m,d)$ by
\eqref{eq:circlepairedderivative}. Moreover $B(r)\to\infty$ as $r\uparrow d$.

\emph{Step 2: finiteness of the correction interaction.}
Let $h_p(a,b)=a^p+b^p-|a-b|^p$. It vanishes on the axes and satisfies
$0\le h_p(a,b)\le a^p+b^p$. It is supermodular, as follows from
\[
 h_p(a,b)=p(p-1)\int_0^a\int_0^b
 |\alpha-\beta|^{p-2}\,d\beta\,d\alpha.
\]
The diagonal singularity is locally integrable for $1<p<2$.
Lemma~\ref{lem:positivekernel}, applied to
$K(r)=r^{-q}(1+|\log r|)$, shows that $u^*\in X^{s,\log,p}(\R)$:
this kernel is decreasing for $q>1$, has integrable tail, and has finite
energy by \eqref{eq:absolute}. Hence both $v=u,u^*$ have finite
$E_{k^+}(v)$, where $k^+=\max\{k,0\}$.

Set $C_d=\sup_{m\le r\le d}|k(r)|<\infty$ and
\[
 J_B(v)=\iint_{I^2}h_p(v(x),v(y))B(|x-y|)\,dx\,dy.
\]
For $m<r<d$, $B(r)\le k^+(d-r)+C_d$. The term with $C_d$ is at most
$2dC_d\|v\|_p^p$. To bound the other term, translate coordinates to
$I=(0,d)$ and use symmetry to restrict to $x<y$. Since
\[
 d-y+x=x-(y-d)=(x+d)-y,
\]
the terms with $v(x)^p$ and $v(y)^p$ are bounded by the left and right
zero-exterior interactions, respectively. Explicitly,
\[
 \begin{aligned}
 &2\int_{\substack{0<x<y<d\\y-x>m}}
 [v(x)^p+v(y)^p]k^+(d-y+x)\,dx\,dy\\
 &\quad\le2\int_I v(x)^p\int_{\R\setminus I}
 k^+(|x-z|)\,dz\,dx\le E_{k^+}(v).
 \end{aligned}
\]
Thus
\begin{equation}\label{eq:circlecorrectionfinite}
 0\le J_B(v)\le E_{k^+}(v)+2dC_d\|v\|_p^p<\infty
 \qquad(v=u,u^*).
\end{equation}

\emph{Step 3: reverse comparison for the correction.}
For $T>0$, extend $B_T=\min\{B,T\}$ to $[0,\infty)$ by $B_T(r)=T$
for $r\ge d$. The kernel $T-B_T$ is bounded, nonnegative and
nonincreasing, with compact support. Apply
Theorem~\ref{thm:supermodularinput} on $\R$ with $F=h_p$ and $f=g=u$.
Equimeasurability preserves the finite quantity
$\iint_{I^2}h_p(v(x),v(y))\,dx\,dy$. Cancelling its multiple $T$ gives
$J_{B_T}(u^*)\le J_{B_T}(u)$. Monotone convergence yields
\begin{equation}\label{eq:circlereverse}
 J_B(u^*)\le J_B(u).
\end{equation}

\emph{Step 4: comparison on a circle and the singular limit.}
Identify the endpoints of $I$ to form the circle $\mathbb T_d$ of
circumference $d$, with arc-length measure and geodesic distance
\[
 \delta(x,y)=\min\{|x-y|,d-|x-y|\}.
\]
Symmetric decreasing rearrangement on this circle about $0$ agrees with
ordinary line rearrangement on $I$: their superlevel sets are the same
centered intervals, of lengths at most $d$.
For $0<\varepsilon<\min\{m,\Rpos\}$ put
$k_\varepsilon(r)=k(\max\{r,\varepsilon\})$. On $[0,m]$ the kernel
$k_\varepsilon-k(m)$ is bounded, nonnegative and nonincreasing.
Theorem~\ref{thm:supermodularinput} on $\mathbb T_d$ with $F=h_p$
therefore gives
\begin{equation}\label{eq:circlecomparison}
 \Delta J_{\mathrm{circ},\varepsilon}
 :=\iint_{I^2}[h_p(u^*(x),u^*(y))-h_p(u(x),u(y))]
 k_\varepsilon(\delta(x,y))\,dx\,dy\ge0.
\end{equation}
Here the added constant $-k(m)$ cancels by equimeasurability.

Define $B_\varepsilon(r)=0$ for $r\le m$ and
$B_\varepsilon(r)=k_\varepsilon(d-r)-k(r)$ for $m<r<d$.
Then $0\le B_\varepsilon\uparrow B$ as $\varepsilon\downarrow0$, and
\[
 k_\varepsilon(|x-y|)
 =k_\varepsilon(\delta(x,y))-B_\varepsilon(|x-y|).
\]
Importantly, $B_\varepsilon$ itself need not be nondecreasing;
\eqref{eq:circlereverse} was proved using $B_T$, not $B_\varepsilon$.
Writing $\Delta J_{\mathrm{line},\varepsilon}$ for the integral in
\eqref{eq:circlecomparison} with $k_\varepsilon(|x-y|)$ in place of
$k_\varepsilon(\delta(x,y))$, we obtain
\[
 \Delta J_{\mathrm{line},\varepsilon}
 =\Delta J_{\mathrm{circ},\varepsilon}
   +J_{B_\varepsilon}(u)-J_{B_\varepsilon}(u^*).
\]
By \eqref{eq:circlecorrectionfinite}, monotone convergence and
\eqref{eq:circlereverse}, its lower limit is nonnegative.

The full-space kernel $k_\varepsilon(|\cdot|)$ is integrable.
Equimeasurability cancels the separate power terms in
$|a-b|^p=a^p+b^p-h_p(a,b)$, giving the exact identity
\[
 \Lsp^{\varepsilon}(u)-\Lsp^{\varepsilon}(u^*)
 =\frac{C_s}{2p}\Delta J_{\mathrm{line},\varepsilon},
\]
where $\Lsp^\varepsilon=(C_s/(2p))E_{k_\varepsilon}$.
For $r<\varepsilon<\Rpos$, $0<k_\varepsilon(r)\le k(r)$, and elsewhere
the kernels agree. Absolute convergence from \eqref{eq:absolute} thus
permits dominated convergence in both full energies. This proves
\eqref{eq:schwarzonedimensional}.

Finally, Proposition~\ref{prop:schwarzcore}, specialized to $N=1$,
constructs the stated smooth counterexamples for every $\ell>2R$.
Combined with the inequality at $d=2R$, this proves exact sharpness.
The equality assertion follows from Corollary~\ref{cor:onedimensionalendpoint} below. Its proof uses only the comparison part established here at smaller positive orders, so there is no circular use of rigidity.
\end{proof}

\subsection{Sharp comparisons for ball superlevel sets and localized cores}
\begin{theorem}[Sharp Schwarz threshold and rigidity for ball superlevel sets]
\label{rem:ballsuperlevelthreshold}
For a measurable function $u:\R^N\to[0,\infty)$ and a level $t>0$, its
\emph{superlevel set} is
\begin{equation}\label{eq:superleveldefinition}
 E_t=\{x\in\R^N:u(x)>t\}.
\end{equation}
Assume $N\ge1$, $1<p<\infty$, $0<s<1$, and the positive normalization
continuously differentiable near $s$ fixed in Theorem~\ref{thm:schwarzonedimensional}.
Suppose $u\ge0$ belongs to $\X$, has bounded essential support, and, for
almost every $t>0$, $E_t$ is either null or agrees up to a null set with
a ball $B_{b_t}(c_t)$. The centers $c_t$ may depend on $t$; a common center
is not assumed. Then
\begin{equation}\label{eq:ballsuperlevelthreshold}
 \diam(\esssupp u)\le2\Rmon
 \quad\Longrightarrow\quad u^*\in\X,\qquad \Lsp(u^*)\le\Lsp(u).
\end{equation}
For $u\not\equiv0$, equality holds exactly when $u=u^*(\cdot-a)$ almost everywhere for one $a\in\R^N$, including at diameter $2\Rmon$.
The bound $2\Rmon$ is sharp within this class in every dimension:
for each $\ell>2\Rmon$ there is such a nonnegative smooth function, with
smooth rearrangement and support diameter less than $\ell$, for which the
inequality is strictly reversed. There is no restriction $sp<1$.
For $N\ge2$ this is a restricted-class result; it does not determine the
Schwarz threshold for arbitrary superlevel sets.

\end{theorem}
\begin{proof}
\emph{Comparison.}
Put $q=N+sp$. Lemma~\ref{lem:positivekernel}, applied to
the decreasing kernel $r^{-q}(1+|\log r|)$, and \eqref{eq:absolute} show
that $u^*\in\X$. Both full signed energies therefore converge absolutely.
{For $0<b\leq\Rmon$, let
\[
 T_b(x)=\int_{|y|>b}k_s(|x-y|)\,dy,\qquad |x|<b.
\]
Lemma~\ref{lem:ballpotential} proves that this exterior-ball potential is strictly radially increasing, also at $b=\Rmon$. The explicit series and its derivation are in Appendix~\ref{sec:ballpotentialappendix}.}

Let $A\subset B=B_b(c)$ be measurable, with $b\le \Rmon$, and let $A_c^*$
be the ball centered at $c$ of measure $|A|$. Radial monotonicity gives
\[
 \int_A T_b(x-c)\,dx\ge\int_{A_c^*}T_b(x-c)\,dx.
\]
Indeed, the part of $A$ outside $A_c^*$ has the same measure as the
missing part inside $A_c^*$, and the potential is no smaller on the
former. Translating both centered sets to the origin yields
\begin{equation}\label{eq:nestedballexterior}
 \int_A\int_{B^c}k_s(|x-y|)\,dy\,dx
 \ge\int_{A^*}\int_{(B^*)^c}k_s(|x-y|)\,dy\,dx.
\end{equation}
These integrals have finite negative parts: $k_s^-(|\cdot|)$ is
integrable, so the negative part is bounded by
$|A|\|k_s^-(|\cdot|)\|_1$. Possible positive boundary divergences are
allowed in this intermediate comparison.

For $a,b\ge0$, direct integration gives the scalar identity
\[
 \begin{aligned}
 |a-b|^p={}&p(p-1)\int_{0<\alpha<\beta<\infty}
 (\beta-\alpha)^{p-2}\\
 &\times\bigl[\mathbf1_{\{a>\beta,\ b\le\alpha\}}
             +\mathbf1_{\{b>\beta,\ a\le\alpha\}}\bigr]
 \,d\alpha\,d\beta.
 \end{aligned}
\]
Apply it first with $|k_s|$ and use Tonelli's theorem. Absolute energy
convergence then justifies Fubini's theorem for $k_s$ and proves
\begin{equation}\label{eq:nestedlevelenergy}
 \Lsp(u)=C_s(p-1)\int_{0<\alpha<\beta<\infty}
 (\beta-\alpha)^{p-2}
 \int_{E_\beta}\int_{E_\alpha^c}k_s(|x-y|)\,dy\,dx
 \,d\alpha\,d\beta.
\end{equation}
For almost every pair $\alpha<\beta$ with nonnull levels, the outer set
$E_\alpha$ is a ball of radius
$b_\alpha\le\diam(\esssupp u)/2\le \Rmon$, and $E_\beta\subset E_\alpha$.
Apply \eqref{eq:nestedballexterior} to this pair and integrate against
the nonnegative weight in \eqref{eq:nestedlevelenergy}. Null levels give
zero contributions. The result is exactly
$\Lsp(u)\ge\Lsp(u^*)$, proving \eqref{eq:ballsuperlevelthreshold}.

\emph{Rigidity, including the endpoint.}
{The strict radial increase of $T_b$, including at $b=\Rmon$, is Lemma~\ref{lem:ballpotential}.} For a measurable $A\subset B_b(c)$ of prescribed measure, the potential integral is consequently minimized uniquely, up to null sets, by the ball centered at $c$. To verify uniqueness, remove the common part of $A$ and this ball. The remaining sets have equal measure and lie strictly outside and inside its bounding sphere; strict radial increase orders their integrals strictly unless they are null.

All the paired deficits in \eqref{eq:nestedlevelenergy} are nonnegative. Finite absolute energy makes the relevant integrals finite for almost every pair of levels. Equality therefore forces $E_\beta$ to be centered at $c_\alpha$ for almost every nonnull pair $\alpha<\beta$. Since $E_\beta$ is a ball, $c_\beta=c_\alpha$ for almost every such pair. Fubini's theorem then gives a single center for almost all nonempty levels; equal-volume nested balls already coincide almost everywhere. Layer cake yields $u=u^*(\cdot-a)$. Translation invariance gives the converse.

\emph{Sharpness within the class.}
For each $\ell>2\Rmon$, Proposition~\ref{prop:schwarzcore} chooses
$\Rmon<r_0<r_1<\ell/2$, a radial nonincreasing smooth background $v$, equal
to $1$ on $B_{r_0}$ and supported in $\overline B_{r_1}$, and a small
nonnegative radial nonincreasing smooth core $w$. For a sufficiently
small nonzero translation $a$, its function $u_a=v+w(\cdot-a)$ satisfies
\[
 \diam(\esssupp u_a)=2r_1<\ell,\qquad
 u_a^*=v+w,\qquad \Lsp(u_a^*)>\Lsp(u_a).
\]
For levels $0<t<1$, its superlevel sets are those of $v$; for $t\ge1$,
they are the translated superlevel sets of the core, when nonempty.
Thus all these sets are balls, and the counterexamples belong to the
stated class. This proves sharpness at $2\Rmon$ in every dimension.
\end{proof}

\begin{proposition}[A radial background with an arbitrary localized core]
\label{prop:arbitrarycore}

The exterior-potential argument also applies to a radial background with
an arbitrary sufficiently localized core. Fix $N\ge1$, $1<p<\infty$,
$0<s<1$, and the same differentiable normalization. Let
$v\in C_c^\infty(\R^N)$ be nonnegative, radial and nonincreasing, with
\[
 v=1\ \text{on }B_{r_0},\qquad
 \supp v=\overline B_{r_1},\qquad 0<r_0<r_1.
\]
Let $w\in C_c^\infty(B_{r_0})$ be nonnegative and satisfy
$\diam(\supp w)\le\Rmon$. No symmetry or assumption on the superlevel
sets of $w$ is imposed. For $u=v+w$,
\begin{equation}\label{eq:arbitrarycorecomparison}
 r_1\le\Rmon\quad\Longrightarrow\quad
 u^*=v+w^*,\qquad \Lsp(u^*)\le\Lsp(u).
\end{equation}
Here $\diam(\esssupp u)=2r_1$. Thus this is another sharp
$2\Rmon$ uniform diameter bound for a specified family, now allowing
nonball superlevel sets at the higher levels. There is no restriction $sp<1$.
\end{proposition}
\begin{proof}
The zero core is immediate. Otherwise choose $b_0<r_0$ with
$\supp w\subset B_{b_0}$. Equimeasurability puts $\supp w^*$ in
$\overline B_{b_0}$. At levels $0<t<1$ the sets $\{u>t\}$ are those
of $v$, whereas at levels $t\ge1$ they are $\{w>t-1\}$. Consequently
$u^*=v+w^*$. Theorem~\ref{thm:schwarz}, applied to the core alone, gives
\begin{equation}\label{eq:arbitrarycoreinternal}
 \Lsp(w^*)\le\Lsp(w),
\end{equation}
because $\diam(\supp w)\le\Rmon$.

For $a,d\ge0$ define
\[
 G_p(a,d)=(a+d)^p-a^p-d^p,\qquad
 V_a(x)=\int_{\R^N}G_p(a,1-v(y))k_s(|x-y|)\,dy.
\]
For $x\in B_{b_0}$ these potentials are finite: the integrand vanishes
on the plateau, so its nonzero part stays a positive distance from $x$,
and the kernel has integrable tail. Direct expansion of the full energy
gives, for $f=w,w^*$,
\begin{equation}\label{eq:arbitrarycoredecomposition}
 \Lsp(v+f)=\Lsp(v)+\Lsp(f)
          +\frac{C_s}{p}\int_{B_{b_0}}V_{f(x)}(x)\,dx.
\end{equation}
Indeed, on pairs inside the plateau the background cancels. On a pair
with only one point in the core the excess integrand is
$G_p(f(x),1-v(y))$; the two orientations give the factor $C_s/p$.
All energies converge absolutely, including those with $w^*$, by
Lemma~\ref{lem:positivekernel} and \eqref{eq:absolute}.

For $\tau\ge0$ let
\[
 \begin{aligned}
 H_\tau(x)&=\left.\partial_aV_a(x)\right|_{a=\tau},\\
 g_\tau(\rho)&=
 p\big[(\tau+1-v(\rho))^{p-1}-\tau^{p-1}\big],
 \end{aligned}
\]
where $v(\rho)$ denotes the radial profile and the derivative at $0$ is
from the right. Differentiation is justified on $B_{b_0}$ by the positive
separation from the background transition. The function $g_\tau$ is
nondecreasing, vanishes on $[0,r_0]$, and is constant on $[r_1,\infty)$.
Expressing it through its nonnegative Stieltjes measure gives
\begin{equation}\label{eq:arbitrarycorepotential}
 H_\tau(x)=\int_{r_0}^{r_1}T_\rho(x)\,dg_\tau(\rho),
 \qquad
 T_\rho(x)=\int_{|y|>\rho}k_s(|x-y|)\,dy.
\end{equation}
For $r_1\le\Rmon$, formula \eqref{eq:ballpotentialseries} shows that
each $T_\rho$ in this representation is radially nondecreasing on
$B_{b_0}$. Therefore $H_\tau$ is also radially nondecreasing there.
Centering the superlevel set of the core minimizes its potential integral:
\[
 \int_{\{w>\tau\}}H_\tau(x)\,dx
 \ge\int_{\{w^*>\tau\}}H_\tau(x)\,dx.
\]
The possibly negative constant $H_\tau(0)$ cancels because these sets
have equal measure. Integrate in $\tau$ from $0$ to $\|w\|_\infty$.
Since $V_0=0$, Fubini's theorem gives
\begin{equation}\label{eq:arbitrarycorecrosscomparison}
 \int_{B_{b_0}}V_{w(x)}(x)\,dx
 \ge\int_{B_{b_0}}V_{w^*(x)}(x)\,dx.
\end{equation}
Fubini is valid here because the core amplitudes are bounded and all
interactions with the background transition have a fixed positive
separation. Equations \eqref{eq:arbitrarycoreinternal},
\eqref{eq:arbitrarycoredecomposition}, and
\eqref{eq:arbitrarycorecrosscomparison} prove
\eqref{eq:arbitrarycorecomparison}.

\emph{Sharpness for this family.}
For each $\ell>2\Rmon$, Proposition~\ref{prop:schwarzcore} chooses
$\Rmon<r_0<r_1<\ell/2$ and a sufficiently small translated smooth radial
core. The core can be chosen with diameter at most $\Rmon$ and remains
compactly inside $B_{r_0}$. Thus these functions belong to the same
background-and-core family and satisfy
\[
 \diam(\esssupp u)<\ell,\qquad \Lsp(u^*)>\Lsp(u).
\]
This proves that $2\Rmon$ is the optimal uniform diameter threshold for
this family. Its geometric assumptions are still essential to the proof;
no unrestricted higher-dimensional Schwarz threshold is inferred.
\end{proof}

\begin{remark}[Explicit examples with moving centers]
\label{rem:explicitballlevels}
The following three families satisfy the superlevel-set assumption in every
dimension. Fix $L>0$, write $e_1=(1,0,\ldots,0)$, and set $z=x/L$.
The notation $a_+=\max\{a,0\}$ denotes the positive part.
\begin{enumerate}[label=(\roman*),leftmargin=*]
\item \emph{An elementary family involving a square root.}
Define
\begin{equation}\label{eq:explicitballgauge}
 \rho(z)=\frac23\left(\sqrt{z_1^2+3|z|^2}-z_1\right),
 \qquad u_m(x)=\bigl(1-\rho(x/L)\bigr)_+^m,
 \qquad m\ge1.
\end{equation}
Squaring the inequality defining $\rho(z)<\lambda$, or completing the
square in the resulting quadratic expression, gives
\[
 \{z:\rho(z)<\lambda\}
 =B_\lambda\!\left(\frac\lambda2e_1\right),\qquad \lambda>0.
\]
Hence, for $0<t<1$,
\begin{equation}\label{eq:explicitgaugesuperlevels}
 \{u_m>t\}
 =B_{L(1-t^{1/m})}\!\left(\frac L2(1-t^{1/m})e_1\right).
\end{equation}
For $t\ge1$ the superlevel set is empty. The centers vary with $t$, and
\[
 \esssupp u_m=\overline B_L\!\left(\frac L2e_1\right),
 \qquad \diam(\esssupp u_m)=2L.
\]

\item \emph{A rational family.}
For $m\ge1$, put
\begin{equation}\label{eq:explicitrationalfamily}
 v_m(x)=\begin{cases}
 \displaystyle\left(\frac{1-|x|^2/L^2}{3+2x_1/L}\right)^m,
       &|x|<L,\\[6pt]
 0,&|x|\ge L.
 \end{cases}
\end{equation}
The denominator is at least $1$ on $\overline B_L$. Write
$\tau=t^{1/m}$ and $\tau_*=(3-\sqrt5)/2$. For $0<t<\tau_*^m$,
completing the square yields
\begin{equation}\label{eq:explicitrationalsuperlevels}
 \{v_m>t\}
 =B_{L\sqrt{1-3\tau+\tau^2}}(-L\tau e_1).
\end{equation}
This ball lies in $B_L$, since
$\tau+\sqrt{1-3\tau+\tau^2}\le1$ for $0<\tau<\tau_*$.
The maximum of $v_m$ is $\tau_*^m$, attained at $-L\tau_*e_1$,
so the superlevel sets are empty for $t\ge\tau_*^m$. In particular,
$\esssupp v_m=\overline B_L$ and its diameter is $2L$.

\item \emph{A smooth compactly supported family.}
Define the explicit smooth cutoff
\begin{equation}\label{eq:explicitsmoothcutoff}
 H(a)=\begin{cases}e^{-1/a},&a>0,\\0,&a\le0,\end{cases}
 \qquad
 \eta(r)=\frac{H(1-r)}{H(1-r)+H(r-\tfrac12)},\quad r\ge0,
\end{equation}
and set
\begin{equation}\label{eq:explicitsmoothballfamily}
 w(x)=\eta\bigl(\rho(x/L)\bigr),
\end{equation}
with $\rho$ from \eqref{eq:explicitballgauge}. The denominator in
\eqref{eq:explicitsmoothcutoff} never vanishes. The function $\eta$ is
$1$ on $[0,1/2]$, strictly decreasing on $(1/2,1)$, and $0$ on
$[1,\infty)$. Its plateau removes the possible nonsmoothness of $\rho$
at the origin, while flatness of the cutoff at its endpoints gives
$w\in C_c^\infty(\R^N)$.
For $0<t<1$, define
\[
 \ell_t=\log\frac{1-t}{t},\qquad
 \lambda_t=\frac34+
 \frac{\ell_t}{4\bigl(\sqrt{16+\ell_t^2}+4\bigr)}.
\]
Solving $\eta(\lambda_t)=t$ gives $1/2<\lambda_t<1$ and
\begin{equation}\label{eq:explicitsmoothsuperlevels}
 \{w>t\}=B_{L\lambda_t}\!\left(\frac{L\lambda_t}{2}e_1\right).
\end{equation}
For $t\ge1$ these sets are empty. Its essential support is
$\overline B_L((L/2)e_1)$ and has diameter $2L$.
\end{enumerate}
The first two families are nonnegative, compactly supported and Lipschitz.
For such a function $f$, its near-diagonal logarithmic energy is bounded
by a constant times
\[
 \int_0^1 r^{p(1-s)-1}\log(1/r)\,dr<\infty,
\]
using $|f(x+h)-f(x)|\le C|h|$ and the bounded measure of the set on which
this difference can be nonzero for $|h|<1$. Thus all three families belong
to $\X$ for every $p>1$ and $0<s<1$. They have no common center for their
nonempty positive superlevel sets. Choosing $L\le\Rmon$ makes
\eqref{eq:ballsuperlevelthreshold} applicable to each family.

\end{remark}

\begin{remark}[Different thresholds for polarization and full symmetrization]
\label{rem:distinctthresholds}
The distinction is surprising: in dimension one, the exact universal
diameter bound is $\Rmon$ for a polarization, but $2\Rmon$ for Schwarz
rearrangement. Thus the kernel's loss of decrease at $\Rmon$ is the exact
obstruction to a comparison valid for every polarization, but it does not
determine the threshold for full rearrangement. Above $\Rmon$, one
polarization can increase the energy while the final Schwarz rearrangement
still decreases it. There is no contradiction with approximation by
polarizations: monotonicity of each intermediate energy is a sufficient
proof mechanism, not a necessary condition for the endpoint comparison.
The circle argument controls the endpoint by an additional comparison.
For Steiner symmetrization the dimension matters: its sharp threshold is
$\Rmon$ for $N\ge2$ by Proposition~\ref{prop:sharpsteiner}, whereas in
dimension one Steiner and Schwarz coincide and both have threshold
$2\Rmon$. No difference between the Steiner and Schwarz thresholds in
the same dimension $N\ge2$ is asserted here, because the exact Schwarz
threshold in those dimensions is still unknown.
\end{remark}

\begin{remark}[Exact threshold for interval superlevel sets in dimension one]
\label{rem:intervalthreshold}
Theorem~\ref{thm:schwarzonedimensional} proves the $2\Rmon$ bound for all
nonnegative admissible functions in dimension one. For functions whose
positive superlevel sets are intervals, up to null sets and almost every
level, there is also the following direct argument.
More precisely, if $u\ge0$ belongs to $X^{s,\log,p}(\R)$, has bounded support,
and has this interval property, then
\begin{equation}\label{eq:intervalthreshold}
 \diam(\esssupp u)\le2\Rmon
 \quad\Longrightarrow\quad \Lsp(u^*)\le\Lsp(u).
\end{equation}
The value $2\Rmon$ is the optimal uniform support-diameter bound within this
class, and by Theorem~\ref{thm:schwarzonedimensional} it is also the
unrestricted one-dimensional threshold. No restriction $sp<1$ is needed.

\emph{Justification.}
Put $R=\Rmon$ and $q=1+sp$. The logarithmic kernel has the midpoint ordering
\begin{equation}\label{eq:midpointordering}
 k_s(m-t)\ge k_s(m+t),\qquad 0<t<m\le R.
\end{equation}
Indeed, $k_s'(r)=pq\,r^{-q-1}\log(r/R)$. For $0<z<m$, the sum of its
values at $m-z$ and $m+z$ is negative: the terms containing $\log(m/R)$
are nonpositive, and the remaining terms have negative sum because
\[
 -\log(1-z/m)>\log(1+z/m)>0,\qquad
 (1-z/m)^{-q-1}>(1+z/m)^{-q-1}.
\]
Integrating this paired derivative over $0<z<t$ proves
\eqref{eq:midpointordering}. For $0<\varepsilon<\Rpos$, the integrable
capped kernel $k_\varepsilon(r)=k_s(\max\{r,\varepsilon\})$ retains that
ordering. To see the only additional case, if $a<\varepsilon<b$, then
$k_\varepsilon(a)=k_s(\varepsilon)\ge k_s(b)$: the positive value
$k_s(\varepsilon)$ is a maximum of $k_s$ on $[\varepsilon,\infty)$.

Consider two nested intervals, after translation written as
\[
 A_c=[c-a,c+a]\subset B=[-b,b],\qquad |c|+a\le b\le R.
\]
Their exterior potential
\[
 T_{b,\varepsilon}(x)=\int_{\R\setminus B}
 k_\varepsilon(|x-y|)\,dy
\]
is even and, for $0<x<b$, satisfies
\[
 T_{b,\varepsilon}'(x)
 =k_\varepsilon(b-x)-k_\varepsilon(b+x)\ge0.
\]
Thus integration of $T_{b,\varepsilon}$ over an interval of length $2a$
contained in $B$ is minimized by $A_0$. Subtracting these potential integrals
from the translation-independent quantity
$|A_c|\int_{\R}k_\varepsilon(|z|)\,dz$ gives
\begin{equation}\label{eq:nestedintervalkernel}
 \int_{A_c}\int_B k_\varepsilon(|x-y|)\,dy\,dx
 \le\int_{A_0}\int_B k_\varepsilon(|x-y|)\,dy\,dx.
\end{equation}
For any two nonempty positive superlevel sets of $u$, the larger interval
has half-length at most $R$ under \eqref{eq:intervalthreshold}. Apply
\eqref{eq:nestedintervalkernel} to each nested pair and integrate against
the nonnegative double-layer weight for
$h_p(a,b)=a^p+b^p-|a-b|^p$, displayed in \eqref{eq:hpdoublelayer} below.
Fubini's theorem is applicable because $k_\varepsilon\in L^1$ and
$h_p(a,b)\le a^p+b^p$. Equimeasurability cancels the separate power terms,
proving the energy comparison for $k_\varepsilon$. The positive singular
part and integrable signed tail argument from Theorem~\ref{thm:polar}
then permits $\varepsilon\downarrow0$, proving
\eqref{eq:intervalthreshold}.

For optimality, the smooth translated-core functions in
Proposition~\ref{prop:schwarzcore} have interval superlevel sets when $N=1$
and violate the comparison with support diameter less than every prescribed
$\ell>2R$. This also shows, with the normalized moment introduced in
Corollary~\ref{cor:schwarzvariational}, that
\[
 \inf\{\Theta(v):v\in\mathcal A_{1,s,p},\
       \text{$v$ has interval positive superlevel sets}\}=\log2,
\]
where $\mathcal A_{1,s,p}$ denotes the diameter-one nonsymmetric class used
there. Corollary~\ref{cor:schwarzvariational} also identifies the
unrestricted one-dimensional infimum as $\log2$.
\end{remark}

\subsection{Comparison beyond the doubled radius for common linear deformations}
The next result gives a geometric distinction complementary to the
moving-center ball examples. Common linear deformations of radial
profiles remain comparable beyond $2\Rmon$ and are uniformly separated
from the unrestricted threshold infimum. The quantitative input is a
translation-moment estimate that is independent of dimension.

{The estimate and its proof are in Appendix~\ref{sec:schwarzthreshold}. The constant \eqref{eq:translationmarginseries} uses all integer translation scales; for $p=2$ it gives a sufficient affine-class diameter of approximately $2.0680\Rmon$.}

{
\begin{theorem}[A uniform moment gap for common linear deformations]
\label{rem:affinethreshold}
Assume $N\ge2$. Let $f\in C_c^\infty(\R^N)$ be nonzero, nonnegative,
radial and nonincreasing, with $\esssupp f=\overline B_{1/2}$. Let
$A\in SL(N,\R)$ be nonorthogonal, put $\Lambda=\|A\|_{\mathrm{op}}$, and define
\[
 v(x)=f(A^{-1}\Lambda x),\qquad \diam(\esssupp v)=1.
\]
For the moment $\Theta$ defined below in
Corollary~\ref{cor:schwarzvariational}, one has the uniform bound
\begin{equation}\label{eq:affinemomentbound}
 \Theta(v)>\log2+\delta_p,\qquad
 \delta_p\text{ is given by \eqref{eq:translationmarginseries}}.
\end{equation}
Consequently every translate and dilation $u$ of such a profile satisfies
\begin{equation}\label{eq:affinethreshold}
 \diam(\esssupp u)\le2e^{\delta_p}\Rmon
 \quad\Longrightarrow\quad \Lsp(u^*)<\Lsp(u).
\end{equation}
Orthogonal $A$ give radial profiles and equality at every scale. The bound
\eqref{eq:affinethreshold} is sufficient; its optimality within this class
is not asserted. The class uses one common linear map for all level sets,
and the bound does not enlarge the sufficient diameter for arbitrary
nonnegative functions.

\end{theorem}
}

\begin{proof}
For $0<\alpha<p$, write
\[
 E_\alpha(w)=\iint_{\R^{2N}}
 \frac{|w(x)-w(y)|^p}{|x-y|^{N+\alpha}}\,dx\,dy,
 \qquad G_\alpha(w)=E_\alpha(w)-E_\alpha(w^*).
\]
In this proof derivatives in $\alpha$ keep $p$ and the function fixed, and
are evaluated at $\alpha=sp$. Differentiating the positive cumulative
deficit representation for $G_\alpha$, proved in
Theorem~\ref{thm:schwarzdeficit}, gives for diameter-one nonsymmetric
profiles
\begin{equation}\label{eq:momentlogderivative}
 \Theta(w)=\partial_\alpha\log G_\alpha(w)-\frac1{N+\alpha}.
\end{equation}
All derivatives here are legitimate for the smooth compactly supported
functions involved. Set $q=N+\alpha$ and
\[
 J(t)=\frac1{|\mathbb S^{N-1}|}\int_{\mathbb S^{N-1}}|A\omega|^{-t}\,d\omega.
\]
A change of variables and the radiality of $f$ give
\begin{equation}\label{eq:affinefactorization}
 G_\alpha(v)=\Lambda^{\alpha-N}E_\alpha(f)[J(q)-1].
\end{equation}
The determinant condition implies $J(0)=J(N)=1$. Since $A$ is nonorthogonal,
strict convexity gives $J(q)>1$ for $q>N$. Put
\[
 P(q)=\frac{J(q)-1}{q(q-N)}>0,\qquad Z(\omega)=-\log|A\omega|.
\]
The second divided difference formula, applied to the exponential, gives
\begin{equation}\label{eq:affinedivdiff}
 P(q)=\frac1{|\mathbb S^{N-1}|}\int_{\mathbb S^{N-1}}
 Z(\omega)^2\int_{\substack{a,b\ge0\\a+b\le1}}
 e^{(Na+qb)Z(\omega)}\,da\,db\,d\omega.
\end{equation}
For completeness, the divided difference of $J$ at $0,N,q$ is exactly
$(J(q)-1)/(q(q-N))$, since its first two values equal $1$; the integral
formula follows by two applications of the fundamental theorem of calculus.
As $Z\ge-\log \Lambda$, $\Lambda\ge1$, and $0\le b\le1$, the logarithmic derivative
of \eqref{eq:affinedivdiff}, a weighted average of $bZ$, satisfies
\[
 \frac{P'(q)}{P(q)}\ge-\log \Lambda.
\]
Combining this with \eqref{eq:momentlogderivative} and
\eqref{eq:affinefactorization} yields
\begin{equation}\label{eq:affineradialreduction}
 \Theta(v)=\partial_\alpha\log(\alpha E_\alpha(f))
       +\log \Lambda+\frac{P'(q)}{P(q)}
 \ge\partial_\alpha\log(\alpha E_\alpha(f)).
\end{equation}

Apply Lemma~\ref{lem:radialtranslationmoment} with $\alpha=sp$ to
\eqref{eq:affineradialreduction}. It gives
$\Theta(v)>\log2+\delta_p$, proving \eqref{eq:affinemomentbound}.
The exact dilation criterion in Corollary~\ref{cor:schwarzvariational}
then yields \eqref{eq:affinethreshold}.
\end{proof}

\begin{corollary}[Exclusion of common linear deformations from threshold sequences]
\label{cor:affineexcluded}
For every $N\ge2$, $s\in(0,1)$ and $p>1$, each diameter-one nonsymmetric
profile $v$ of Theorem~\ref{rem:affinethreshold} satisfies
\[
 \Theta_s(v)-\theta_{N,s,p}\ge\delta_p>0.
\]
Thus no sequence of these profiles can minimize the unrestricted Schwarz
threshold, even when both the radial profile and the linear map vary.
\end{corollary}
\begin{proof}
The theorem gives $\Theta_s(v)\ge\log2+\delta_p$, whereas the
translated-core construction and Corollary~\ref{cor:schwarzvariational}
give $\theta_{N,s,p}\le\log2$.
\end{proof}

\section{Optimality of our results}
This section distinguishes universal optimality from counterexamples that apply only to particular geometries. Throughout, the order $s$ is a differentiability point of the fixed positive normalization, and
\[
 0<\Rmon=\exp\!\left(\frac{b_s}{p}+\frac1{N+sp}\right)<\infty.
\]
No additional upper or lower bound on this radius is required: the counterexample scales are chosen relative to its value. Every strict counterexample below has original support diameter larger than $\Rmon$. The polarization constructions in Propositions~\ref{prop:counterpolar} and
\ref{prop:sharpdiameter} place their active distances in the increasing
range $(\Rmon,\infty)$. Proposition~\ref{prop:sharpsteiner} uses separated
transverse bumps to prove failure at every enlarged Steiner diameter when
$N\ge2$. Proposition~\ref{prop:dilation} gives a separate exact scaling
argument for both symmetrizations. For Schwarz rearrangement,
Proposition~\ref{prop:schwarzcore} uses exterior-potential concavity and proves
failure for every prescribed bound $\ell>2\Rmon$; it does not prove failure for
every $\ell>\Rmon$. In dimension one the sufficient circle comparison in
Theorem~\ref{thm:schwarzonedimensional} matches this upper bound and gives
the exact Schwarz and Steiner threshold $2\Rmon$.
\subsection{A smooth counterexample for one polarization}
\begin{proposition}[Smooth polarization counterexample]\label{prop:counterpolar}
Fix $N\ge1$, $1<p<\infty$, $0<s<1$, and a positive normalization $C_t$ differentiable near $s$, with $\Rmon$ defined by \eqref{eq:radii}. For every $T>\Rmon$ and $0<\varepsilon<T/4$, the two-bump construction below gives a nonnegative $u\in C_c^\infty(\R^N)$ and a half-space $H$ such that
\[
 \diam(\esssupp u)>\Rmon,\qquad \Lsp(u^H)>\Lsp(u).
\]
All active inter-bump distances used in the defect calculation lie above $\Rmon$.
\end{proposition}
\begin{proof}
Let $H=\{x_1>0\}$, $e=e_1$, $T>\Rmon$, and $0<\varepsilon<T/4$. Choose a nonzero nonnegative radial $\varphi\in C_c^\infty(\R^N)$ supported in $\overline B_\varepsilon$ and put
\[
 u(x)=\varphi(x-3Te)+\varphi(x+Te).
\]
The reflection pairs of the bumps are disjoint, and radial symmetry gives
\[
 u^H(x)=\varphi(x-3Te)+\varphi(x-Te).
\]
For $\xi,\eta\in B_\varepsilon$, let $x=3Te+\xi$ and $y=Te+\eta$. Then
\[
 r=|x-y|>2T-2\varepsilon>\Rmon,\qquad
 \rho=|x-\sigma_Hy|>r.
\]
Both distances lie in the strictly increasing region of $k_s$. Define $h_p(a,b)=a^p+b^p-|a-b|^p$, which is strictly positive for $a,b>0$. Self-energies of the bumps cancel, leaving
\[
 \Lsp(u^H)-\Lsp(u)
 =\frac{C_s}{p}\int_{B_\varepsilon}\int_{B_\varepsilon}
 h_p(\varphi(\xi),\varphi(\eta))
 [k_s(\rho)-k_s(r)]\,d\eta\,d\xi>0.
\]
Every integral is finite because the bumps are smooth and their mutual distances are positive. Choosing an occupied point in each bump shows that the original support diameter is at least $4T-2\varepsilon>\Rmon$, so the small-support hypothesis is violated.
\end{proof}

\begin{remark}[An explicit choice and a positive lower bound]\label{rem:explicitpolar}
The construction can be made with no unspecified bump or separation parameter. Set
\[
 T=2\Rmon,\qquad \varepsilon=T/8,\qquad
 \varphi(x)=\begin{cases}
 \exp\!\left(-\dfrac{1}{1-|x|^2/\varepsilon^2}\right),&|x|<\varepsilon,\\
 0,&|x|\ge\varepsilon.
 \end{cases}
\]
Use the two functions in Proposition~\ref{prop:counterpolar}. For every pair of local variables in $B_\varepsilon$,
\[
 r\le 2T+2\varepsilon=\frac{9T}{4},\qquad
 \rho\ge4T-2\varepsilon=\frac{15T}{4}.
\]
Since both ranges lie above $\Rmon$,
\[
 k_s(\rho)-k_s(r)\ge
 \gamma_s:=k_s(15T/4)-k_s(9T/4)>0.
\]
On $B_{\varepsilon/2}$, $\varphi\ge m:=e^{-4/3}$. The function $h_p$ is increasing in each variable on $[0,\infty)^2$, as is seen by differentiating in either variable; thus $h_p(a,b)\ge h_p(m,m)=2m^p$ when $a,b\ge m$. Consequently the strict failure has the quantitative bound
\begin{equation}\label{eq:polarbound}
 \displaystyle
 \Lsp(u^H)-\Lsp(u)
 \ge\frac{2C_s}{p}e^{-4p/3}
 |B_{\varepsilon/2}|^2\,
 [k_s(15T/4)-k_s(9T/4)]>0.
\end{equation}
This counterexample works for every $N\ge1$, $1<p<\infty$, and $0<s<1$, with the stated differentiable normalization.
\end{remark}

\subsection{Sharpness of the diameter bound for polarization}
\begin{proposition}[Optimal uniform polarization radius]\label{prop:sharpdiameter}
Fix $N\ge1$, $1<p<\infty$, $0<s<1$, and the differentiable normalization. For every $\ell>\Rmon$, there are a nonnegative $u\in C_c^\infty(\R^N)$ and a half-space $H$ such that
\begin{equation}\label{eq:sharpdiameter}
 \diam(\esssupp u)<\ell,\qquad \Lsp(u^H)>\Lsp(u).
\end{equation}
Consequently, among diameter-only hypotheses, the largest bound guaranteeing $\Lsp(u^H)\le\Lsp(u)$ for every nonnegative $u\in\X$ and every $H$ is exactly $\Rmon$.
\end{proposition}
\begin{proof}
Set $R=\Rmon$ and choose
\[
 \delta=\frac{\ell-R}{4},\qquad \varepsilon=\frac\delta8,
 \qquad a=R+2\delta,\qquad b=\delta.
\]
Let $H=\{x_1>0\}$, with reflection $\sigma$, and let $\varphi\in C_c^\infty(\R^N)$ be nonzero, radial, nonnegative, and supported in $\overline B_\varepsilon$. Define
\[
 u(x)=\varphi(x-ae_1)+\varphi(x+be_1).
\]
All the relevant bump supports and their reflection pairs are disjoint, since $b>\varepsilon$ and $a-b=R+\delta>2\varepsilon$. Hence
\[
 u^H(x)=\varphi(x-ae_1)+\varphi(x-be_1).
\]
The original support diameter satisfies
\[
 \diam(\esssupp u)\le a+b+2\varepsilon
 =R+\frac{13}{4}\delta<R+4\delta=\ell.
\]
For $\xi,\eta\in B_\varepsilon$, put $x=ae_1+\xi$, $y=be_1+\eta$, and
\[
 r=|x-y|,\qquad \rho=|x-\sigma y|.
\]
Both points lie in the interior of $H$, so $\rho>r$. Moreover,
\[
 r\ge a-b-2\varepsilon=R+\frac34\delta>R,
 \qquad
 \rho\le a+b+2\varepsilon<\ell.
\]
Thus both distances lie in the strictly increasing region of $k_s$, and $k_s(\rho)>k_s(r)$. The same self-energy cancellation as in Proposition~\ref{prop:counterpolar} gives
\[
 \Lsp(u^H)-\Lsp(u)
 =\frac{C_s}{p}\int_{B_\varepsilon}\int_{B_\varepsilon}
 h_p(\varphi(\xi),\varphi(\eta))
 [k_s(\rho)-k_s(r)]\,d\eta\,d\xi>0.
\]
Indeed, $h_p(\alpha,\beta)=\alpha^p+\beta^p-|\alpha-\beta|^p>0$ when both arguments are positive, which occurs on a product set of positive measure. All integrals are finite for smooth compactly supported functions. This proves \eqref{eq:sharpdiameter}. Theorem~\ref{thm:polar} supplies sufficiency at $\ell=\Rmon$, completing the optimality statement.
\end{proof}

{
\begin{corollary}[A sorting step increases energy while its symmetric endpoint decreases it]\label{cor:step-endpoint-separation}
Let $N=1$, $p>1$, $0<s<1$, with the differentiable normalization fixed above. For every $\ell$ with $\Rmon<\ell\leq2\Rmon$, there exist a nonnegative $u\in C_c^\infty(\R)$ and a half-line $H$ such that
\[
 \Rmon<\diam(\esssupp u)<\ell,\qquad
 \Lsp(u^*)<\Lsp(u)<\Lsp(u^H).
\]
In particular $(u^H)^*=u^*$, so the energy increase of this intermediate step coexists with a strictly lower energy at the fully symmetrized endpoint.
\end{corollary}
\begin{proof}
Use the two-bump construction in Proposition~\ref{prop:sharpdiameter} with the given $\ell$, taking the radial bump strictly positive inside its supporting ball. It gives $\diam(\esssupp u)<\ell$ and $\Lsp(u^H)>\Lsp(u)$. The diameter must exceed $\Rmon$ by Theorem~\ref{thm:polar}. The two separated nonzero bumps cannot be a translate of a symmetric decreasing function: their positivity set has two separated components, whereas a nonzero symmetric decreasing function has interval positive superlevel sets and an interval positivity set. Since $\ell\leq2\Rmon$, Theorem~\ref{thm:schwarzonedimensional}, including its rigidity assertion, gives $\Lsp(u^*)<\Lsp(u)$. Equimeasurability under polarization gives $(u^H)^*=u^*$.
\end{proof}
}

\begin{remark}[Uniform sharpness versus individual necessity]
The proposition asserts failure for at least one function and one half-space whenever the allowed diameter bound is enlarged. It does not assert failure for every function whose support diameter exceeds $\Rmon$. It also does not imply failure of a Steiner or Schwarz comparison at every diameter above that radius: a single polarization comparison and the final rearrangement comparison are different statements. The separate construction in Proposition~\ref{prop:sharpsteiner} settles Steiner sharpness for $N\ge2$; Corollary~\ref{cor:schwarzthreshold} only bounds the Schwarz threshold.
\end{remark}

\subsection{Sharpness for Steiner symmetrization in dimensions at least two}
\begin{proposition}[Optimal uniform Steiner radius]\label{prop:sharpsteiner}
Assume $N\ge2$ and fix a unit vector $e$. For every $\ell>\Rmon$, there is a nonnegative $u\in C_c^\infty(\R^N)$ such that
\[
 \diam(\esssupp u)<\ell,\qquad \Lsp(u^{\star_e})>\Lsp(u).
\]
Consequently, $\Rmon$ is the optimal uniform support-diameter bound for Steiner symmetrization in direction $e$.
\end{proposition}
\begin{proof}
Put $R=\Rmon$, choose a unit vector $f\perp e$, and set
\[
 \delta=\frac{\ell-R}{4},\quad a=R+\delta,\quad b=R+2\delta,
 \quad h=\sqrt{b^2-a^2}.
\]
Choose $0<\varepsilon<\delta/8$ and a nonzero nonnegative radial decreasing $\varphi\in C_c^\infty(B_\varepsilon)$. Define
\[
 u(x)=\varphi(x)+\varphi(x-af-he).
\]
The projections of the two bump supports onto $e^\perp$ are disjoint, since $a>2\varepsilon$. Each line parallel to $e$ therefore meets at most one bump. Centering its slice gives exactly
\[
 u^{\star_e}(x)=\varphi(x)+\varphi(x-af).
\]
The original support has diameter at most $b+2\varepsilon<\ell$.

Write $h_p(\alpha,\beta)=\alpha^p+\beta^p-|\alpha-\beta|^p$. Expansion of the energy of two disjoint nonnegative bumps, followed by cancellation of their translation-invariant individual energies, yields
\begin{align*}
 &\Lsp(u^{\star_e})-\Lsp(u)\\
 &=\frac{C_s}{p}\iint_{B_\varepsilon\times B_\varepsilon}
 h_p(\varphi(\xi),\varphi(\eta))
 \bigl[k_s(|af+he+\eta-\xi|)-k_s(|af+\eta-\xi|)\bigr]\,d\xi\,d\eta.
\end{align*}
All distances here exceed $R$, and
\[
 |af+he+\eta-\xi|\ge b-2\varepsilon
 >a+2\varepsilon\ge |af+\eta-\xi|.
\]
Since $k_s$ is strictly increasing on $(R,\infty)$,
\begin{align*}
 \Lsp(u^{\star_e})-\Lsp(u)
 &\ge \frac{C_s}{p}
 [k_s(b-2\varepsilon)-k_s(a+2\varepsilon)]\\
 &\quad\times\iint_{B_\varepsilon\times B_\varepsilon}
 h_p(\varphi(\xi),\varphi(\eta))\,d\xi\,d\eta>0.
\end{align*}
Strictness follows because $h_p(\alpha,\beta)>0$ for $\alpha,\beta>0$. Theorem~\ref{thm:steiner} supplies sufficiency at $\ell=R$.
\end{proof}
\begin{remark}
The transverse vector $f$ is essential to this construction. In dimension one, Steiner symmetrization is Schwarz rearrangement and Proposition~\ref{prop:sharpsteiner} makes no assertion. Nor does this two-bump proof establish Schwarz sharpness in higher dimensions: Schwarz rearrangement merges the bumps and changes their individual energies.
\end{remark}

\subsection{Dilation obstruction for Steiner and Schwarz rearrangement}
\begin{proposition}[Dilation counterexamples and additional equality scales]\label{prop:dilation}
Let $T$ denote either fixed-direction Steiner symmetrization or Schwarz rearrangement. There is a smooth nonnegative compactly supported $v$, with $Tv$ also smooth, for which all sufficiently large dilations $v_R(x)=v(x/R)$ satisfy
\[
 \Lsp(v_R)<\Lsp(Tv_R).
\]
Writing $d_v=\diam(\esssupp v)$ and $\Delta(t)=F_t(v)-F_t(Tv)$, a sufficient explicit condition for strict failure is
\[
 R>\max\left\{\exp\!\left(\frac{\Delta'(s)}{p\Delta(s)}\right),
 \frac{\Rmon}{d_v}\right\}.
\]
In particular, $\diam(\esssupp v_R)=Rd_v>\Rmon$. There is also a dilation giving equality although $v_R$ is not an equality configuration of the corresponding small-support theorem.
\end{proposition}
\begin{proof}
Choose a nonzero smooth radial decreasing bump $\varphi$ supported in $B_1$, positive there, and set $v(x)=\varphi(x-ae)+\varphi(x+ae)$ with $a>1$. For Steiner symmetrization in direction $e$, writing $x=z+te$ gives $Tv(z+te)=\varphi(z+(t/2)e)$. For Schwarz rearrangement, $Tv(x)=\varphi(x/2^{1/N})$. Thus both transforms are smooth. The two-bump function is neither a translate of its Steiner symmetrization in direction $e$ nor a translate of its Schwarz rearrangement.

The ordinary fractional rearrangement gap
\[
 \Delta(t)=F_t(v)-F_t(Tv)
\]
satisfies $\Delta(s)>0$ by Corollary~\ref{cor:fractionalinput} with $t=s$.
The functions $v$ and $Tv$ are smooth and compactly supported, so their
fractional seminorms are finite. For Schwarz rearrangement, equality in that
corollary would require $v$ to be a translate of $v^*$; for Steiner
symmetrization it would require one translation in direction $e$ of
$v^{\star_e}$. The separated two-bump positivity set rules out both
configurations. Since $C_s>0$, strictness also holds for the normalized
fractional gap $\Delta(s)$.

Smoothness implies differentiability of $\Delta$ near $s$, with differentiation justified uniformly on a compact subinterval of $(0,1)$. Dilation commutes with $T$, and
\[
 F_t(v_R)-F_t(Tv_R)=R^{N-tp}\Delta(t).
\]
Differentiating this identity, not an inequality, gives
\begin{equation}\label{eq:scalinggap}
 \Lsp(v_R)-\Lsp(Tv_R)
 =R^{N-sp}\bigl(\Delta'(s)-p\log R\,\Delta(s)\bigr).
\end{equation}
It is negative for $R>\exp(\Delta'(s)/(p\Delta(s)))$ and zero at that threshold. The additional displayed condition $R>\Rmon/d_v$ makes the support-size violation explicit. Dilation cannot turn the two bumps into the corresponding symmetry configuration. The equality scale must also have support diameter larger than $\Rmon$, by the rigidity statements in Theorems~\ref{thm:steiner} and~\ref{thm:schwarz}.
\end{proof}

\begin{remark}
Proposition~\ref{prop:dilation} shows that the equality classifications in the small-support theorems cannot be extended to arbitrary supports: outside the monotonicity regime, exact cancellation can produce additional equality examples. Likewise, differentiating the fractional P\'olya--Szeg\H{o} inequality is invalid; a nonnegative gap need not have a nonnegative derivative. Each counterexample is compactly supported, so no unrestricted version valid on every bounded zero exterior domain can hold.
\end{remark}

\subsection{An explicit smooth counterexample for Steiner symmetrization}
The next construction gives an explicit smooth function, identifies its Steiner symmetrization, and specifies a dilation at which the energy strictly increases. Thus the counterexample does not rely on an unspecified large scale.

\begin{proposition}[Explicit smooth Steiner counterexample]\label{prop:explicitsteiner}
Fix a unit vector $e$ and write $x=z+te$. Define
\[
 \psi(x)=\begin{cases}
 \exp\!\left(-\dfrac{1}{1-|x|^2}\right),&|x|<1,\\
 0,&|x|\ge1,
 \end{cases}
 \qquad
 v(x)=\psi(x-2e)+\psi(x+2e).
\]
Set $w(z+te)=\psi(z+(t/2)e)$ and define the two finite constants
\begin{equation}\label{eq:DGexplicit}
 \Delta_s=F_s(v)-F_s(w)>0,\qquad
 G_s=\Lsp(v)-\Lsp(w).
\end{equation}
In particular, these constants are specified by the convergent integrals \eqref{eq:F} and \eqref{eq:L}. Choose
\begin{equation}\label{eq:explicitR}
 R=\max\left\{\exp\!\left(1+\frac{|G_s|}{p\Delta_s}\right),\Rmon\right\},
 \qquad u_R(x)=v(x/R).
\end{equation}
Then $u_R$ and its Steiner symmetrization are nonnegative smooth compactly supported functions, and
\begin{align}
 u_R^{\star_e}(z+te)&=\psi\!\left(\frac zR+\frac{t}{2R}e\right),\label{eq:explicitSteiner}\\
 \Lsp(u_R^{\star_e})-\Lsp(u_R)
 &=R^{N-sp}\bigl(p\log R\,\Delta_s-G_s\bigr)\notag\\
 &\ge p\Delta_sR^{N-sp}>0.\label{eq:steinerbound}
\end{align}
Here $\diam(\esssupp u_R)=6R>\Rmon$. Thus Steiner symmetrization strictly increases the full logarithmic energy outside the small-support regime.
\end{proposition}
\begin{proof}
For a fixed $z$ with $|z|<1$, $t\mapsto\psi(z+te)$ is even and decreasing for $t>0$. The two translated slices of $v$ are disjoint. For every positive level, their superlevel lengths add, so rearrangement replaces two identical intervals centered at $\pm2$ by one centered interval of twice the length. This proves $v^{\star_e}=w$, including the zero slices with $|z|\ge1$; dilation gives \eqref{eq:explicitSteiner}. In particular both $v$ and $w$ are smooth.

Corollary~\ref{cor:fractionalinput}, with $t=s$ and Steiner direction $e$,
gives $\Delta_s>0$: smoothness ensures finite fractional energy, and the two
separated slice components exclude the common-translation equality case. Absolute convergence gives $G_s\in\R$.

For every smooth compactly supported $f$, order differentiation of the exact fractional scaling identity gives
\[
 \Lsp(f(\cdot/R))
 =R^{N-sp}\bigl(\Lsp(f)-p\log R\,F_s(f)\bigr).
\]
Subtract this identity for $f=v$ and $f=w$. Since dilation commutes with Steiner symmetrization, we obtain
\[
 \Lsp(u_R^{\star_e})-\Lsp(u_R)
 =R^{N-sp}\bigl(p\log R\,\Delta_s-G_s\bigr).
\]
The choice \eqref{eq:explicitR} gives $p\log R\,\Delta_s\ge p\Delta_s+|G_s|$, which proves \eqref{eq:steinerbound}. The two unit-ball supports of $v$, centered at $\pm2e$, have union diameter $6$, so dilation gives the stated support diameter $6R$.
\end{proof}

\begin{remark}[Interpretation and bounded domains]
The scale in \eqref{eq:explicitR} is explicit in terms of two specified finite energy integrals; no elementary closed form for these integrals is asserted. Both counterexamples concern the full signed energy, and do not contradict the small-support inequalities. For the polarization example in Remark~\ref{rem:explicitpolar}, the ball $B_{4T}(0)$ contains both supports and is invariant under the chosen polarization. For the Steiner example, $B_{4R}(0)$ contains the supports of both $u_R$ and $u_R^{\star_e}$ and is Steiner symmetric. Thus both failures also occur for zero exterior energies on sufficiently large balls, with the original and symmetrized function admissible in the same domain.
\end{remark}

\subsection{An explicit Schwarz counterexample}
\begin{proposition}[Explicit smooth Schwarz counterexample]\label{prop:explicitschwarz}
Let $\psi$ and $v(x)=\psi(x-2e)+\psi(x+2e)$ be the smooth bumps from Proposition~\ref{prop:explicitsteiner}. Put
\[
 w_*(x)=\psi(x/2^{1/N}),\qquad
 \Delta_*=F_s(v)-F_s(w_*)>0,\qquad
 G_*=\Lsp(v)-\Lsp(w_*),
\]
and choose
\[
 R_*=\max\left\{\exp\!\left(1+\frac{|G_*|}{p\Delta_*}\right),\Rmon\right\}.
\]
For $f(x)=v(x/R_*)$, the support diameter is $6R_*>\Rmon$, and one has
\[
 f^*(x)=\psi\!\left(\frac{x}{2^{1/N}R_*}\right),
 \qquad
 \Lsp(f^*)-\Lsp(f)
 \ge p\Delta_*R_*^{N-sp}>0.
\]
Both functions are nonnegative, smooth, and compactly supported.
\end{proposition}
\begin{proof}
At every positive level, $v$ has two disjoint equal ball superlevel sets. Their combined volume is that of a ball whose radius is larger by the factor $2^{1/N}$. Thus $v^*=w_*$. Corollary~\ref{cor:fractionalinput}, with $t=s$ and Schwarz rearrangement,
gives $\Delta_*>0$: $v$ has finite fractional energy and its two-component
positivity set excludes the translated-radial equality case. Smoothness gives finite $G_*$. The exact dilation identity \eqref{eq:scalinggap}, with $T=*$ and scale $R_*$, yields
\[
 \Lsp(f^*)-\Lsp(f)
 =R_*^{N-sp}(p\log R_*\,\Delta_*-G_*)
 \ge R_*^{N-sp}(p\Delta_*+|G_*|-G_*),
\]
which gives the asserted bound. Taking a centered ball of radius larger than $\max\{3,2^{1/N}\}R_*$ makes both functions admissible in the same zero exterior space. This does not contradict Theorem~\ref{thm:schwarz}, whose support-diameter hypothesis is essential.
\end{proof}

\subsection{A translated-core Schwarz counterexample at twice the monotonicity radius}
\begin{proposition}[Translated-core Schwarz counterexamples above $2\Rmon$]\label{prop:schwarzcore}
Fix $N\ge1$, $1<p<\infty$, and $0<s<1$, with the differentiable normalization used throughout. For every $\ell>2\Rmon$, there exists a nonnegative $u\in C_c^\infty(\R^N)$ with $u^*\in C_c^\infty(\R^N)$ such that
\[
 \diam(\esssupp u)<\ell,\qquad \Lsp(u^*)>\Lsp(u).
\]
The supports of $u$ and $u^*$ can be the same centered ball. No restriction $sp<1$ is required.
\end{proposition}
\begin{proof}
\textit{Step 1: the background and its tail potentials.}
Put $R=\Rmon$ and choose $R<r_0<r_1<\ell/2$. Take a radial nonincreasing $v\in C_c^\infty(\R^N)$, equal to $1$ on $B_{r_0}$, strictly decreasing from $1$ to $0$ on $r_0<|x|<r_1$, and zero outside $B_{r_1}$. For $\rho\in[r_0,r_1]$ define
\[
 T_\rho(x)=\int_{|y|>\rho}k_s(|x-y|)\,dy.
\]
These functions are smooth for $|x|<r_0/2$: the integration stays away from the singularity and all differentiated tails are integrable. Integration by parts over the exterior of $B_\rho$, with vanishing boundary term at infinity, gives
\begin{equation}\label{eq:tailhessian}
 D^2T_\rho(0)
 =-\frac{|\mathbb S^{N-1}|}{N}\rho^{N-1}k_s'(\rho)I.
\end{equation}
Indeed, the inner boundary has exterior-domain normal $-y/|y|$, and the spherical integral of $\theta_i\theta_j$ is $|\mathbb S^{N-1}|\delta_{ij}/N$. Since $k_s'(\rho)>0$ on this compact interval, continuity yields $\eta\in(0,r_0/2)$ and $\kappa>0$ such that
\begin{equation}\label{eq:tailconcavity}
 D^2T_\rho(x)\le-\kappa I
 \quad (|x|<\eta,\ r_0\le\rho\le r_1).
\end{equation}
Matrix inequalities here mean inequalities of quadratic forms.

\textit{Step 2: the exact rearrangement.}
Choose a nonzero nonnegative radial nonincreasing $w\in C_c^\infty(B_\varepsilon)$ with $0<\varepsilon<\eta/4$. For $|a|<\eta/4$ put
\[
 u_a(x)=v(x)+w(x-a).
\]
The core remains inside the plateau of $v$. For levels below $1$, the superlevel sets are those of $v$; for levels at least $1$, they are the translated ball superlevel sets of the core. Consequently
\begin{equation}\label{eq:corestar}
 u_a^*=v+w=u_0,\qquad
 \supp u_a=\supp u_0=\overline B_{r_1}.
\end{equation}
In particular the support diameter is $2r_1<\ell$.

\textit{Step 3: cancellation and concavity.}
For $\alpha,d\ge0$ set
\[
 G_p(\alpha,d)=(\alpha+d)^p-\alpha^p-d^p,
 \qquad
 V_\alpha(x)=\int_{\R^N}G_p(\alpha,1-v(y))k_s(|x-y|)\,dy.
\]
These potentials are finite and smooth for $|x|<r_0/2$. Writing $v(y)=\widehat v(|y|)$ and $g_\alpha(r)=G_p(\alpha,1-\widehat v(r))$, we have
\begin{align*}
 g_\alpha'(r)
 &=p\bigl[(\alpha+1-\widehat v(r))^{p-1}
 -(1-\widehat v(r))^{p-1}\bigr](-\widehat v'(r))\ge0,\\
 V_\alpha(x)&=\int_{r_0}^{r_1}g_\alpha'(\rho)T_\rho(x)\,d\rho.
\end{align*}
The second identity follows by expressing $g_\alpha(|y|)$ as an integral of exterior-ball indicators; absolute integrability justifies interchanging integrals. Moreover
\[
 \int_{r_0}^{r_1}g_\alpha'(\rho)\,d\rho=G_p(\alpha,1),
\]
which is positive for $\alpha>0$. Thus \eqref{eq:tailconcavity} implies
\[
 D^2V_\alpha(x)\le-\kappa G_p(\alpha,1)I\qquad(|x|<\eta).
\]

Because $v=1$ on the core, direct expansion of the full energy gives the exact identity
\begin{equation}\label{eq:coreenergy}
 \Lsp(u_a)=\Lsp(v)+\Lsp(w)+\frac{C_s}{p}F(a),
 \qquad
 F(a)=\int_{B_\varepsilon}V_{w(\xi)}(a+\xi)\,d\xi.
\end{equation}
For pairs of points in the core the background cancels; for pairs with only one point in the core the extra integrand is precisely $G_p(w,1-v)$. The two orientations of such pairs account for the coefficient $C_s/p$.
Radial symmetry gives $F(-a)=F(a)$, so $\nabla F(0)=0$. The preceding Hessian estimate gives
\[
 D^2F(a)\le-\kappa M_w I,
 \qquad M_w=\int_{B_\varepsilon}G_p(w(\xi),1)\,d\xi>0,
 \qquad |a|<\eta/4.
\]
Therefore, for $0<|a|<\eta/4$, Taylor's formula and \eqref{eq:coreenergy} imply
\[
 \Lsp(u_a^*)-\Lsp(u_a)
 =\frac{C_s}{p}[F(0)-F(a)]
 \ge\frac{C_s\kappa M_w}{2p}|a|^2>0.
\]
All original energies converge absolutely because the functions are smooth and compactly supported. This proves the proposition.
\end{proof}

\begin{corollary}[Bounds for the optimal Schwarz diameter]\label{cor:schwarzthreshold}
Let $\Dsch$ be the supremum of all $d>0$ such that $\Lsp(u^*)\le\Lsp(u)$ for every nonnegative $u\in\X$ with $\diam(\esssupp u)\le d$. Then
\[
 \begin{cases}
 \Dsch=2\Rmon,&N=1,\\
 \Rmon\le \Dsch\le2\Rmon,&N\ge2.
 \end{cases}
\]
\end{corollary}
\begin{proof}
Theorem~\ref{thm:schwarz} proves the lower bound. If some admissible diameter bound $d$ exceeded $2\Rmon$, Proposition~\ref{prop:schwarzcore} with $\ell=d$ would contradict its defining universal inequality. For $N=1$, Theorem~\ref{thm:schwarzonedimensional} gives the matching lower bound $2\Rmon$, including validity at that endpoint.
\end{proof}
\begin{remark}[What remains unresolved in dimensions at least two]
The exact one-dimensional threshold is $2\Rmon$ by
Theorem~\ref{thm:schwarzonedimensional}. For $N\ge2$, the translated-core
construction gives counterexamples only above $2\Rmon$, while the
unrestricted sufficient theorem gives $\Rmon$. Neither this construction
nor the one-dimensional circle argument settles the Schwarz threshold in
those dimensions. In particular, sharpness of $\Rmon$ for Steiner
symmetrization in $N\ge2$ does not settle sharpness for Schwarz rearrangement.
\end{remark}

\section{Concluding discussion}\label{sec:conclusion}
The sharp polarization and Schwarz diameters in dimension one are $\Rmon$ and $2\Rmon$, respectively. Corollary~\ref{cor:step-endpoint-separation} realizes this distinction with a smooth function: one polarization increases its energy while Schwarz rearrangement strictly decreases it. The circle correction proves the final comparison even when individual sorting steps need not decrease the energy. Corollary~\ref{cor:onedimensionalendpoint} completes the rigidity statement at the doubled diameter.

The distinction between $\Rpos$ and $\Rmon$ is equally structural.
For nonnegative functions, control of the active kernel differences
recovers comparison beyond the sign change, while the signed tail
remains part of the functional. The modulus reduction uses the smaller
radius $\Rpos$. The resulting theory concerns the full fixed-positive-order
energy and is distinct both from its positive-part version and from
the renormalized order-zero energy.

Theorem~\ref{thm:schwarzdeficit} explains the scale dependence: each nonsymmetric profile has one critical dilation. Corollary~\ref{cor:schwarzvariational} obtains the universal threshold by minimizing the critical diameter over profiles normalized to diameter one. Theorem~\ref{thm:ordermoment} then orders the profile moments and the normalized threshold as the fractional order varies. The uniform moment gap for common linear deformations excludes that family from unrestricted threshold-minimizing sequences (Corollary~\ref{cor:affineexcluded}).

\paragraph{Status of the sharp higher-dimensional threshold.}
The explicit value of the unrestricted Schwarz threshold for $N\ge2$ remains open. The universal bounds proved here are $\Rmon\le \Dsch\le2\Rmon$. The doubled diameter is optimal within the class of functions having ball superlevel sets, even when their centers depend on the level. Further estimates for common affine deformations concern restricted classes and must not be interpreted as determining the universal threshold. These distinctions show that the geometry of the entire nested family of superlevel sets matters. The signed-moment
criterion in Proposition~\ref{prop:doubledcriterion} specifies the
remaining estimate; the local core obstruction in
Proposition~\ref{prop:coreobstruction} explains a limitation of the
translated-core construction below the doubled radius.

Theorem~\ref{thm:nonlinearcounter} identifies a further obstruction for functionals depending nonlinearly on the pointwise energy density. Even on supports strictly inside the valid comparison radius, strict increase and convexity in the density variable do not ensure rearrangement decrease. Theorem~\ref{thm:affinedensity} supplies a positive result under a convex-weight condition controlling the coupling between function values and energy densities. For integrands depending directly on distance and two function values, Theorem~\ref{thm:Phi} gives a sorting-gap criterion.

\appendix
\section{Deficit representation and the Schwarz threshold}\label{sec:schwarzthreshold}
This appendix begins with the exterior-ball potential calculation used in the ball-superlevel comparison. It then expresses the Schwarz energy gap through the cumulative Riesz rearrangement deficit. A weighted geometric mean of interaction distances determines the sign of the gap and the unique critical dilation of each nonsymmetric profile. Optimizing over profiles gives a variational characterization of the universal diameter threshold. The representation applies to bounded-support nonnegative functions in the stated energy space for all $sp>0$; evaluating its infimum in dimensions $N\ge2$ remains open.

{
\subsection{Exterior-ball potential monotonicity}\label{sec:ballpotentialappendix}
\begin{lemma}[Strict radial increase of the exterior-ball potential]\label{lem:ballpotential}
Let $q=N+sp$ and $0<b\leq\Rmon$. The function
\[
 T_b(x)=\int_{|y|>b}k_s(|x-y|)\,dy,\qquad |x|<b,
\]
is finite, smooth and radial on $B_b$, and strictly increases with $|x|$ there, including when $b=\Rmon$. It has the locally convergent expansion \eqref{eq:ballpotentialseries} below.
\end{lemma}
\begin{proof}
Finiteness and smoothness on compact subsets follow from positive separation from the boundary singularity and integrability of every differentiated tail. Rotational invariance makes the potential radial.
Write $\sigma_{N-1}=|\mathbb S^{N-1}|$ and, for $0\le r<b$, consider
\[
 A_q(r)=\int_{|y|>b}|re_1-y|^{-q}\,dy
       =\sum_{j=0}^{\infty}a_j(q)b^{N-q-2j}r^{2j},
 \qquad
 a_j(q)=\frac{\sigma_{N-1}}{q-N}
 \frac{(q/2)_j((q-N)/2)_j}{(N/2)_j\,j!}.
\]
Here $(a)_0=1$ and $(a)_j=a(a+1)\cdots(a+j-1)$ for $j\ge1$.
This expansion follows by spherical averaging and radial integration.
More explicitly, the coefficient of $r^{2j}$ in the spherical average
of $|y-re_1|^{-q}$ is obtained from
\[
 \Delta^j|y|^{-q}
 =\left[\prod_{\ell=0}^{j-1}(q+2\ell)(q-N+2+2\ell)\right]
 |y|^{-q-2j}
\]
by division by $4^j j!(N/2)_j$. Integration over $|y|>b$ gives the
displayed $a_j(q)$. The series converges locally uniformly for $r<b$,
also after differentiating in $r$ or in the exponent $q>N$.
For $j\ge1$, cancellation of the factor $q-N$ gives
\[
 \frac{a_j'(q)}{a_j(q)}
 =\frac1q+\sum_{\ell=1}^{j-1}
 \left(\frac1{q+2\ell}+\frac1{q-N+2\ell}\right).
\]
In differentiating with respect to this auxiliary exponent, $b_s$ and
$p$ are held fixed. Since $k_s(r)=(b_s+p\partial_q)r^{-q}$ and
$\log \Rmon=b_s/p+1/q$, it follows that
\begin{equation}\label{eq:ballpotentialseries}
 \begin{aligned}
 T_b(re_1)-T_b(0)
 ={}&p\sum_{j=1}^{\infty}a_j(q)b^{N-q-2j}r^{2j}\\
 &\times\left[\log\frac{\Rmon}{b}+
 \sum_{\ell=1}^{j-1}\left(\frac1{q+2\ell}
                    +\frac1{q-N+2\ell}\right)\right].
 \end{aligned}
\end{equation}
All coefficients are nonnegative for $b\le \Rmon$.

When $b<\Rmon$, the quadratic coefficient is positive. When $b=\Rmon$, it vanishes but every coefficient with $j\geq2$ is positive. Consequently the series strictly increases between any two distinct radii in $[0,b)$.
\end{proof}
}

\subsection{The cumulative deficit and its critical dilation}
\begin{theorem}[Cumulative deficit representation]\label{thm:schwarzdeficit}
Let $u\ge0$ belong to $\X$, let
$d=\diam(\esssupp u)\in(0,\infty)$, and set $q=N+sp$. Define
\begin{align}
 h_p(a,b)&=a^p+b^p-|a-b|^p,\notag\\
 M_u(r)&=\iint_{|x-y|<r}
 \bigl[h_p(u^*(x),u^*(y))-h_p(u(x),u(y))\bigr]\,dx\,dy.
 \label{eq:cumulativedeficit}
\end{align}
Then $u^*\in\X$, $M_u(r)\ge0$ for every $r>0$, and $M_u(r)=0$ for
$r\ge d$. Moreover,
\begin{equation}\label{eq:schwarzdeficitformula}
 \displaystyle
 \Lsp(u)-\Lsp(u^*)
 =\frac{C_s q}{2}\int_0^d
 M_u(r)r^{-q-1}\log\frac{\Rmon}{r}\,dr.
\end{equation}
The integral is absolutely convergent. If $u$ is not almost everywhere a
translate of $u^*$, put
\begin{equation}\label{eq:deficitprobability}
 \begin{gathered}
 W_u=\int_0^d M_u(r)r^{-q-1}\,dr,\qquad
 d\nu_u(r)=W_u^{-1}M_u(r)r^{-q-1}\,dr,\\
 g_u=\exp\!\left(\int_0^d\log r\,d\nu_u(r)\right).
 \end{gathered}
\end{equation}
Then $0<W_u<\infty$, $0<g_u<d$, and
\begin{equation}\label{eq:geometriccriterion}
 \Lsp(u)-\Lsp(u^*)
 =\frac{C_s q}{2}W_u\log\frac{\Rmon}{g_u}.
\end{equation}
Consequently the Schwarz inequality holds for this function precisely when
$g_u\le\Rmon$, and it is strictly reversed precisely when $g_u>\Rmon$.
\end{theorem}
\begin{proof}
\emph{Step 1: positivity and support of the cumulative deficit.}
For all $a,b\ge0$,
\begin{equation}\label{eq:hpdoublelayer}
 h_p(a,b)=p(p-1)\int_0^a\int_0^b
 |\alpha-\beta|^{p-2}\,d\beta\,d\alpha.
\end{equation}
For $p=2$ the integrand is $1$; for $1<p<2$ its diagonal singularity
is locally integrable. The identity follows by integration of the mixed
second derivative and the boundary values $h_p(a,0)=h_p(0,b)=0$.
Apply Theorem~\ref{thm:rieszinput}, whose sets have finite measure here because $u$ has bounded support, with the kernel
$\mathbf1_{B_r}(x-y)$ to the two superlevel sets
$\{u>\alpha\}$ and $\{u>\beta\}$, and integrate the resulting nonnegative
deficits against the weight in \eqref{eq:hpdoublelayer}. Tonelli's theorem
proves $M_u(r)\ge0$. This is also a special case of the supermodular
Riesz inequalities in \cite{BH}. All integrals defining $M_u$ are finite:
$0\le h_p(a,b)\le a^p+b^p$, $h_p(a,0)=0$, and the supports have finite
measure.

Theorem~\ref{thm:isodiametricinput}, applied to the bounded essential support of $u$, gives
$\diam(\esssupp u^*)\le d$. For $r\ge d$ both double integrals in
\eqref{eq:cumulativedeficit} include all pairs of nonzero values.
Their equality follows from equimeasurability, so $M_u(r)=0$ there.
Also apply Lemma~\ref{lem:positivekernel} to
$K(r)=r^{-q}(1+|\log r|)$. This kernel is nonnegative and decreasing
because $q>1$, its tail is integrable, and its energy is finite by
\eqref{eq:absolute}. The lemma gives finite $K$-energy of $u^*$, hence
$u^*\in\X$ by the definition of the logarithmic seminorm.

\emph{Step 2: the singular integration by parts.}
For $v=u,u^*$ write
$A_v(z)=\int_{\R^N}|v(x+z)-v(x)|^p\,dx$. Cancellation of the
$L^p$ terms, using $\|u^*\|_p=\|u\|_p$, gives
\[
 M_u(r)=\int_{|z|<r}[A_u(z)-A_{u^*}(z)]\,dz.
\]
In particular $M_u$ is locally absolutely continuous. For sufficiently
small $\varepsilon>0$, $k_s$ is positive and decreasing on
$(0,\varepsilon]$, whence
\[
 |k_s(\varepsilon)M_u(\varepsilon)|
 \le\int_{|z|<\varepsilon}k_s(|z|)
 [A_u(z)+A_{u^*}(z)]\,dz\longrightarrow0.
\]
This uses absolute convergence of both energies near the origin.
To check the other integrals used below, bound $|M_u(r)|$ by
$\int_{|z|<r}(A_u+A_{u^*})(z)\,dz$ and use Tonelli's theorem.
The elementary estimate, for $0<t<d$,
\[
 \int_t^d(1+|\log r|)r^{-q-1}\,dr
 \le C_d\bigl[1+t^{-q}(1+|\log t|)\bigr]
\]
then shows that
$\int_0^d M_u(r)(1+|\log r|)r^{-q-1}\,dr<\infty$.
Integration by parts, first away from zero and then letting
$\varepsilon\downarrow0$, therefore yields
\[
 \Lsp(u)-\Lsp(u^*)
 =\frac{C_s}{2p}\int_0^d k_s(r)M_u'(r)\,dr
 =-\frac{C_s}{2p}\int_0^d k_s'(r)M_u(r)\,dr.
\]
The boundary term at $d$ is zero. Since
\[
 k_s'(r)=pq\,r^{-q-1}\log\frac r{\Rmon},
\]
this proves \eqref{eq:schwarzdeficitformula}.

\emph{Step 3: normalization and strict positivity.}
The same calculation for $r^{-q}$ gives
\[
 \iint\frac{|u(x)-u(y)|^p-|u^*(x)-u^*(y)|^p}{|x-y|^q}\,dx\,dy
 =qW_u.
\]
Both individual fractional energies are finite by \eqref{eq:absolute}.
Corollary~\ref{cor:fractionalinput}, with $t=s$, therefore applies to $u$.
Its Schwarz equality criterion shows that $W_u>0$ whenever $u$ is not a
translate of $u^*$. The preceding integrability estimates
show that $\log r$ is $\nu_u$-integrable. Since $\nu_u$ is a probability
measure on $(0,d)$, its geometric mean satisfies $0<g_u<d$.
Substitution into \eqref{eq:schwarzdeficitformula} proves
\eqref{eq:geometriccriterion} and the sign assertions.
\end{proof}

\paragraph{A unique critical dilation.}
For a nonsymmetric profile $u$ as in Theorem~\ref{thm:schwarzdeficit},
put $u_\lambda(x)=u(x/\lambda)$, $\lambda>0$. The scaling laws
$W_{u_\lambda}=\lambda^{N-sp}W_u$ and $g_{u_\lambda}=\lambda g_u$ give
\begin{equation}\label{eq:profilecrossover}
 \Lsp(u_\lambda)-\Lsp(u_\lambda^*)
 =\frac{C_sq}{2}\lambda^{N-sp}W_u
 \log\frac{\Rmon}{\lambda g_u}.
\end{equation}
Since $W_u>0$, the gap is strictly positive for $\lambda<\Rmon/g_u$,
zero at $\lambda=\Rmon/g_u$, and strictly negative above that value.
The equality at this one scale occurs for a nonsymmetric function.
Its support diameter is $d\Rmon/g_u>\Rmon$, in agreement with the
small-support rigidity theorem. This is an exact description of the
comparison along every such dilation family.

\begin{corollary}[Exact variational characterization of the Schwarz threshold]
\label{cor:schwarzvariational}
Let $\mathcal A_{N,s,p}$ consist of the nonnegative functions $v\in\X$
with $\diam(\esssupp v)=1$ that are not translates of $v^*$. For these
functions define
\[
 \Theta(v)=\int_0^1\log\frac1r\,d\nu_v(r)>0,
 \qquad
 \theta_{N,s,p}=\inf_{v\in\mathcal A_{N,s,p}}\Theta(v).
\]
Then
\begin{equation}\label{eq:exactschwarzthreshold}
 \displaystyle
 \Dsch=\Rmon\exp(\theta_{N,s,p}),\qquad
 0\le\theta_{N,s,p}\le\log2.
\end{equation}
In dimension one,
\begin{equation}\label{eq:thetaonedimensional}
 \theta_{1,s,p}=\log2.
\end{equation}
In particular $\theta_{N,s,p}$ is independent of the choice of the
positive differentiable normalization $C_t$. Exact Schwarz sharpness at
$\Rmon$ is equivalent to the existence of $v_j\in\mathcal A_{N,s,p}$
such that
\begin{equation}\label{eq:missingschwarzconcentration}
 \int_0^1\log\frac1r\,d\nu_{v_j}(r)\longrightarrow0.
\end{equation}
\end{corollary}
\begin{proof}
For $v_\lambda(x)=v(x/\lambda)$, rearrangement commutes with dilation and
\[
 M_{v_\lambda}(r)=\lambda^{2N}M_v(r/\lambda),\quad
 W_{v_\lambda}=\lambda^{N-sp}W_v,\quad
 g_{v_\lambda}=\lambda g_v=\lambda e^{-\Theta(v)}.
\]
Its support diameter is $\lambda$. Formula
\eqref{eq:geometriccriterion} says that this profile satisfies the
Schwarz inequality exactly for
$\lambda\le\Rmon e^{\Theta(v)}$.
Every bounded-support nonnegative function that is not a translate of its
rearrangement is a dilation of a member of $\mathcal A_{N,s,p}$.
Hence every diameter bound at most
$\inf_v\Rmon e^{\Theta(v)}$ guarantees the inequality, and each larger
bound admits a strictly reversed comparison by choosing a profile with
smaller critical diameter and then a dilation strictly between the two
bounds. This proves the equality in \eqref{eq:exactschwarzthreshold},
including validity of the comparison at that endpoint.
Corollary~\ref{cor:schwarzthreshold} supplies the bounds on $\theta_{N,s,p}$.
Its definition contains only $N,s,p$ and the normalized profiles, proving
normalization independence. The final equivalence is the definition of a
zero infimum of the nonnegative quantities $\Theta(v)$.
For $N=1$, Theorem~\ref{thm:schwarzonedimensional} and
\eqref{eq:exactschwarzthreshold} give $e^{\theta_{1,s,p}}=2$, proving
\eqref{eq:thetaonedimensional}. Thus $\Theta(v)\ge\log2$ for every
$v\in\mathcal A_{1,s,p}$, and \eqref{eq:missingschwarzconcentration} is
impossible in dimension one.
\end{proof}

\subsection{Fractional deficit bounds}
\begin{theorem}[Comparison of logarithmic and fractional deficits]
\label{thm:quantdeficit}
For bounded-support nonnegative $u\in\X$, set
\[
 d=\diam(\esssupp u)>0,\qquad
 \Delta F_s(u)={F_s(u)-F_s(u^*)},\qquad
 \Delta L_s(u)=\Lsp(u)-\Lsp(u^*).
\]
If $u$ is not a translate of $u^*$ and $v(x)=u(dx)$, then
\begin{equation}\label{eq:exactdeficitratio}
 \Delta L_s(u)=p\Delta F_s(u)
 \left(\log\frac{\Rmon}{d}+\Theta(v)\right).
\end{equation}
In particular, in every dimension,
\begin{equation}\label{eq:basicquantdeficit}
 \Delta L_s(u)\ge p\log\frac{\Rmon}{d}\,\Delta F_s(u).
\end{equation}
For $N=1$, and in every dimension for the ball-superlevel class of
Theorem~\ref{rem:ballsuperlevelthreshold}, the stronger estimate is
\begin{equation}\label{eq:doublequantdeficit}
 \Delta L_s(u)\ge p\log\frac{2\Rmon}{d}\,\Delta F_s(u).
\end{equation}
Both estimates are strict for every nonsymmetric $u$ in their respective classes, at every $d>0$. Symmetric translates have both deficits zero, so the stated non-strict inequalities also include them.
\end{theorem}
\begin{proof}
The proof of Theorem~\ref{thm:schwarzdeficit} gives
$\Delta F_s(u)=C_sqW_u/(2p)$, while
\eqref{eq:geometriccriterion} gives
$\Delta L_s(u)=C_sqW_u\log(\Rmon/g_u)/2$.
Dilation gives $g_u=de^{-\Theta(v)}$, proving
\eqref{eq:exactdeficitratio}. Since $\Theta(v)\ge0$, this proves
\eqref{eq:basicquantdeficit}. In dimension one, use
\eqref{eq:thetaonedimensional}.
For the ball-superlevel class, dilation preserves the class.
Apply Theorem~\ref{rem:ballsuperlevelthreshold} to the dilation of $u$
having diameter $2\Rmon$. The exact scaling identity gives
\[
 0\le\lambda^{N-sp}
 \bigl(\Delta L_s(u)-p\log\lambda\,\Delta F_s(u)\bigr),
 \qquad\lambda=2\Rmon/d,
\]
which is \eqref{eq:doublequantdeficit}. The basic estimate is strict for
nonsymmetric $u$ because $\Theta(v)>0$ and $\Delta F_s(u)>0$.
For the doubled estimate, the dilation above is still nonsymmetric.
Its logarithmic deficit is strictly positive by the endpoint rigidity
in Corollary~\ref{cor:onedimensionalendpoint} when $N=1$, or by
Theorem~\ref{rem:ballsuperlevelthreshold} for ball superlevel sets.
This proves strictness also in \eqref{eq:doublequantdeficit}.
The argument holds for every $d>0$, although the right-hand side becomes
negative when $d>2\Rmon$.
\end{proof}

\begin{theorem}[A strict two-order deficit bound at the doubled diameter]
\label{thm:twoorderdeficit}
Fix $0<\sigma<s<1$. Let $u\ge0$ belong to $X^{s,\log,p}(\R^N)$, have bounded essential support of diameter $d>0$, and not be a translate of $u^*$. Assume either $N=1$, or that $u$ has ball superlevel sets as in Theorem~\ref{rem:ballsuperlevelthreshold}. For $t\in\{\sigma,s\}$ define the unnormalized fractional deficits
\[
 \mathfrak D_t(u)=\iint_{\R^{2N}}
 \frac{|u(x)-u(y)|^p-|u^*(x)-u^*(y)|^p}
 {|x-y|^{N+tp}}\,dx\,dy>0,
\]
and put
\begin{equation}\label{eq:twoorderfactor}
 Z_{s,\sigma}(u)
 =\left(\frac d2\right)^{p(s-\sigma)}
   \frac{N+\sigma p}{N+sp}
   \frac{\mathfrak D_s(u)}{\mathfrak D_\sigma(u)}.
\end{equation}
Then $Z_{s,\sigma}(u)>1$, and with $\Rmon$ evaluated at order $s$,
\begin{equation}\label{eq:twoorderdeficit}
 \Delta L_s(u)>
 p\log\frac{2\Rmon}{d}\,\Delta F_s(u)
 +\frac{\Delta F_s(u)}{s-\sigma}\log Z_{s,\sigma}(u).
\end{equation}
In particular, at the sharp doubled diameter,
\begin{equation}\label{eq:twoorderendpoint}
 d=2\Rmon\quad\Longrightarrow\quad
 \Delta L_s(u)>
 \frac{\Delta F_s(u)}{s-\sigma}\log Z_{s,\sigma}(u)>0.
\end{equation}
The factor $Z_{s,\sigma}$ is invariant under spatial dilations and independent of the energy normalization. Symmetric translates have both deficits zero and are excluded from its definition.
\end{theorem}
\begin{proof}
Membership at order $s$ gives membership and finite fractional energies at $\sigma$. The fractional rigidity statement in Corollary~\ref{cor:fractionalinput} gives positivity of both $\mathfrak D_t$. Normalize the support diameter by setting $v(x)=u(dx)$, and write $X(r)=\log(1/r)$ on $(0,1)$. With the cumulative deficit $M_v$ from Theorem~\ref{thm:schwarzdeficit}, set
\[
 W_t=\int_0^1 M_v(r)r^{-N-tp-1}\,dr,
 \qquad d\nu_t(r)=W_t^{-1}M_v(r)r^{-N-tp-1}\,dr.
\]
These are finite positive normalizing constants and probability measures. For $a=p(s-\sigma)>0$, the order tilt gives
\[
 d\nu_s=\frac{e^{aX}}{A}\,d\nu_\sigma,
 \qquad A=\int e^{aX}\,d\nu_\sigma=\frac{W_s}{W_\sigma}.
\]
The fractional cumulative-deficit identity and spatial scaling yield
\[
 \mathfrak D_t(v)=(N+tp)W_t,\qquad
 \mathfrak D_t(v)=d^{tp-N}\mathfrak D_t(u),\qquad
 A=2^aZ_{s,\sigma}(u).
\]
The comparison at the doubled diameter, applied at order $\sigma$ to a dilation of $v$, gives $\Theta_\sigma(v)\ge\log2$. This uses Theorem~\ref{thm:schwarzonedimensional} in dimension one and Theorem~\ref{rem:ballsuperlevelthreshold} for ball superlevel sets; their moment consequence is independent of normalization. To apply those theorems, one may choose any auxiliary positive normalization continuously differentiable near $\sigma$, for example $\widetilde C_t\equiv1$. The required dilation uses its own monotonicity radius, which cancels from the resulting moment bound. Thus the argument also covers $\sigma=1/2$, even when the piecewise normalization \eqref{eq:normalization} is discontinuous there. The measure $\nu_\sigma$ has a density and is not concentrated at one radius. Strict Jensen inequality therefore gives
\[
 \log A>a\int X\,d\nu_\sigma
 =a\Theta_\sigma(v)\ge a\log2,
\]
so $Z_{s,\sigma}(u)>1$.

For the measure at order $s$, the same tilt identity gives
$\int e^{-aX}\,d\nu_s=A^{-1}$. Another strict Jensen inequality yields
\[
 -\log A=\log\int e^{-aX}\,d\nu_s
 >-a\int X\,d\nu_s=-a\Theta_s(v).
\]
All first moments are finite by Theorem~\ref{thm:schwarzdeficit}. Hence
\[
 \Theta_s(v)>\log2+\frac{\log Z_{s,\sigma}(u)}{p(s-\sigma)}.
\]
Substitution in the exact deficit ratio \eqref{eq:exactdeficitratio} proves \eqref{eq:twoorderdeficit}, and $\Delta F_s(u)>0$ gives \eqref{eq:twoorderendpoint}. The powers of the spatial dilation cancel in \eqref{eq:twoorderfactor}, and its definition uses only unnormalized fractional deficits. This proves both invariances.
\end{proof}

\begin{remark}[Quantitative strictness without a uniform endpoint gap]
The positive correction in \eqref{eq:twoorderendpoint} depends on the profile and two ordinary fractional deficits. It improves the zero lower bound obtained by substituting $d=2\Rmon$ into \eqref{eq:doublequantdeficit}. No profile-independent positive endpoint gap is asserted. In higher dimensions this estimate retains the ball-superlevel hypothesis and does not determine the unrestricted Schwarz threshold.
\end{remark}

\subsection{Order monotonicity and rigidity at the sharp endpoint}
When the order varies, write $\Theta_s(v)$ and $\nu_{s,v}$ for the moment
and deficit probability measure. The cumulative function $M_v$ depends
on $p$ and $v$, but not on $s$.

\begin{theorem}[Strict profile monotonicity and monotonicity of the threshold]
\label{thm:ordermoment}
Fix $N\ge1$ and $p>1$. If $0<s_1<s_2<1$ and $v$ is a diameter-one
nonsymmetric profile in $X^{s_2,\log,p}(\R^N)$, then
\begin{equation}\label{eq:strictprofileorder}
 \Theta_{s_2}(v)>\Theta_{s_1}(v).
\end{equation}
Consequently
\begin{equation}\label{eq:thresholdorder}
 \theta_{N,s_1,p}\le\theta_{N,s_2,p},\qquad
 \frac{\Dsch(s_1)}{\Rmon(s_1)}
 \le\frac{\Dsch(s_2)}{\Rmon(s_2)}.
\end{equation}
For a fixed smooth nonsymmetric diameter-one profile,
\begin{equation}\label{eq:momentvariance}
 \frac{d}{ds}\Theta_s(v)
 =p\operatorname{Var}_{\nu_{s,v}}\!\left(\log\frac1r\right)>0.
\end{equation}
No monotonicity of the unnormalized diameter $\Dsch(s)$ is asserted.
The moment statements are independent of the normalization. The threshold
ratios are understood at differentiability points of the chosen
normalization, as required throughout the paper.
\end{theorem}
\begin{proof}
Membership at $s_2$ implies membership at $s_1$, since for $0<r<1$ the
lower-order logarithmic seminorm has the smaller weight. Set
$X(r)=\log(1/r)$ and $a=p(s_2-s_1)>0$. The measures satisfy
\[
 d\nu_{s_2,v}
 =\frac{e^{aX}}{\int e^{aX}\,d\nu_{s_1,v}}\,d\nu_{s_1,v}.
\]
The normalizing integral and its first tilted moment are finite by
Theorem~\ref{thm:schwarzdeficit} at order $s_2$.
For $\mu=\nu_{s_1,v}$, write
\[
 \begin{split}
 &\int Xe^{aX}\,d\mu
 -\left(\int X\,d\mu\right)\left(\int e^{aX}\,d\mu\right)\\
 &\quad=\frac12\iint
 (X(r)-X(t))(e^{aX(r)}-e^{aX(t)})\,d\mu(r)\,d\mu(t)>0.
 \end{split}
\]
The integrand is nonnegative. It is positive whenever $r\ne t$, and
$\mu$ has a density with respect to Lebesgue measure, with nonzero total
mass. Thus $X$ is not $\mu$-almost everywhere constant. Integrability of
the displayed terms justifies the identity, or one may first truncate
and pass to the limit. Dividing by $\int e^{aX}\,d\mu$ proves
\eqref{eq:strictprofileorder}.

The admissible classes satisfy
$\mathcal A_{N,s_2,p}\subset\mathcal A_{N,s_1,p}$, and hence
\[
 \inf_{v\in\mathcal A_{N,s_2,p}}\Theta_{s_2}(v)
 \ge\inf_{v\in\mathcal A_{N,s_2,p}}\Theta_{s_1}(v)
 \ge\inf_{v\in\mathcal A_{N,s_1,p}}\Theta_{s_1}(v).
\]
This proves \eqref{eq:thresholdorder} using
\eqref{eq:exactschwarzthreshold}. For a smooth fixed profile, higher
logarithmic moments are finite locally in $0<s<1$. {To justify this without assuming smoothness of the rearrangement, choose a slightly larger order below one. The fractional comparison bounds the rearranged fractional seminorm by that of the smooth profile, and the extra power controls the additional logarithmic factors at the smaller order.} Differentiating the
normalized exponential tilt gives \eqref{eq:momentvariance}; its
variance is positive by the same nonconcentration argument. Strictness
for each profile need not survive the infimum over profiles.
\end{proof}

\begin{corollary}[Complete one-dimensional rigidity at the sharp diameter]
\label{cor:onedimensionalendpoint}
For $N=1$ and every nonsymmetric diameter-one profile
$v\in X^{s,\log,p}(\R)$,
\begin{equation}\label{eq:strictlogtwo}
 \Theta_s(v)>\log2.
\end{equation}
Thus the infimum $\theta_{1,s,p}=\log2$ is not attained. For every nonzero
nonnegative $u\in X^{s,\log,p}(\R)$ with
$\diam(\esssupp u)\le2\Rmon$, equality in the logarithmic Schwarz
inequality holds if and only if $u=u^*(\cdot-a)$ almost everywhere for
one $a\in\R$.
\end{corollary}
\begin{proof}
Choose $0<\sigma<\min\{s,1/2\}$. Membership at $s$ implies membership at
$\sigma$. The comparison part of Theorem~\ref{thm:schwarzonedimensional}
and Corollary~\ref{cor:schwarzvariational} at order $\sigma$ give
$\Theta_\sigma(v)\ge\log2$. The strict profile inequality
\eqref{eq:strictprofileorder} gives
$\Theta_s(v)>\Theta_\sigma(v)\ge\log2$. No equality classification is
used in this argument. If $u$ is nonsymmetric,
Corollary~\ref{cor:fractionalinput} gives $\Delta F_s(u)>0$ and
\eqref{eq:exactdeficitratio} is strictly positive even when
$d=2\Rmon$. Conversely, translations preserve the energy.
\end{proof}

\begin{remark}[Sharp comparison versus attainment]
The strict inequality \eqref{eq:strictlogtwo} is compatible with the exact
threshold $2\Rmon$: the counterexamples of
Proposition~\ref{prop:schwarzcore} give diameter-one profiles whose moments
approach $\log2$. The minimizing profiles vary. The same reasoning applies
in every dimension to the ball-superlevel class, whose infimum is
$\log2$ by its sharp dilation threshold. In particular its moment
infimum is also not attained. More generally, if
$\theta_{N,\sigma,p}=\theta_{N,s,p}$ for some $\sigma<s$, strict profile
monotonicity excludes nonsymmetric equality at the unrestricted sharp
diameter at order $s$.
\end{remark}

{
\subsection{Translation moments for common affine deformations}
\begin{lemma}[A logarithmic moment bound from translations]
\label{lem:radialtranslationmoment}
Let $f\in C_c^\infty(\R^N)$ be nonzero, nonnegative, radial and
nonincreasing, supported in $\overline B_{1/2}$. For $0<\alpha<p$, set
\[
 E_\alpha(f)=\iint_{\R^{2N}}
 \frac{|f(x)-f(y)|^p}{|x-y|^{N+\alpha}}\,dx\,dy.
\]
Define the dimension-independent constant
\begin{equation}\label{eq:translationmarginseries}
 \delta_p=\sum_{n=3}^{\infty}n^{-p}
       \left(\log\frac n{n-1}-\frac{2}{n(n-1)}\right)>0.
\end{equation}
Then
\begin{equation}\label{eq:radiallogmomentmargin}
 \partial_\alpha\log\bigl(\alpha E_\alpha(f)\bigr)
 >\log2+\delta_p.
\end{equation}
\end{lemma}
\begin{proof}
Fix a unit vector $e$ and put
\[
 a(r)=\int_{\R^N}|f(x+re)-f(x)|^p\,dx,
 \qquad A=a(1)=2\|f\|_p^p.
\]
The last equality follows because the supports are disjoint at translation
one, up to null sets. Radiality makes $a$ independent of $e$. Moreover,
$a$ is nondecreasing. Indeed,
$2\|f\|_p^p-a(r)=\int h_p(f(x+re),f(x))\,dx$, and the positive
double-layer representation of $h_p$ expresses this quantity as a
positive mixture of intersection volumes of two balls, each
nonincreasing in their separation. We have $a(0)=0$ and $a(r)=A$ for
$r\ge1$.

\emph{A first-moment estimate.}
Writing $x=z+te$, take $I=(-1/2,1/2)$; every slice of $f$ vanishes
outside $I$. By symmetry and Fubini,
\[
 \begin{split}
 \int_0^1\bigl(A-a(r)\bigr)\,dr
 &=\frac12\int_{e^\perp}\int_{I^2}
 h_p\bigl(f(z+te),f(z+ve)\bigr)\,dt\,dv\,dz\\
 &\le\|f\|_p^p=\frac A2,
 \end{split}
\]
because $h_p(a,b)\le a^p+b^p$ and $|I|=1$. Thus
$\int_0^1a(r)\,dr\ge A/2$. Integration by parts gives
\[
 d\lambda(r)=A^{-1}\,da(r),\qquad
 \int_0^1 r\,d\lambda(r)
 =1-A^{-1}\int_0^1a(r)\,dr\le\frac12.
\]
Here $\lambda$ is a probability measure on $(0,1]$.

\emph{All integer translation scales.}
For each integer $n\geq2$, the triangle inequality for $n$ consecutive translations of length $1/n$ gives
\[
 a(1)^{1/p}\leq n\,a(1/n)^{1/p},\qquad
 \lambda((0,1/n])\geq n^{-p}.
\]
Tilt $\lambda$ by the decreasing weight $r^{-\alpha}$:
\[
 d\mu_\alpha(r)=\frac{r^{-\alpha}\,da(r)}{\int_0^1t^{-\alpha}\,da(t)}.
\]
A decreasing tilt decreases the mean of $r$ and increases the probability of each initial interval. The covariance identity for monotone functions therefore gives
\begin{equation}\label{eq:translationtiltconstraints}
 \int r\,d\mu_\alpha\leq\frac12,\qquad
 F_\alpha(r):=\mu_\alpha((0,r])\geq n^{-p}
 \quad\left(\frac1n\leq r\leq\frac1{n-1},\ n\geq3\right).
\end{equation}
These integrals are finite because $a(r)=O(r^p)$ at zero and $\alpha<p$. Since $a$ is continuous and increases from $0$ to $A$, the probability measure $\lambda$ is not a point mass. For $\alpha>0$ the strictly decreasing tilt therefore strictly decreases its first moment: $\int r\,d\mu_\alpha<\int r\,d\lambda\leq1/2$.

\emph{Integrating the tangent remainder over every scale.}
For $0<r\leq1$,
\[
 -\log r=\log2+1-2r+\phi(r),\qquad
 \phi(r)=-\log(2r)-1+2r\geq0.
\]
On $(0,1/2]$, $\phi(r)=\int_r^{1/2}(1/t-2)\,dt$. Tonelli's theorem and \eqref{eq:translationtiltconstraints} imply
\[
 \begin{aligned}
 \int\phi\,d\mu_\alpha
 &\geq\int_0^{1/2}(1/r-2)F_\alpha(r)\,dr\\
 &\geq\sum_{n=3}^{\infty}n^{-p}
       \int_{1/n}^{1/(n-1)}(1/r-2)\,dr=\delta_p.
 \end{aligned}
\]
Combining this with $\int r\,d\mu_\alpha\leq1/2$ proves
\[
 \int_0^1\log(1/r)\,d\mu_\alpha(r)>\log2+\delta_p.
\]
Finally, polar coordinates and Stieltjes integration by parts give
\[
 \alpha E_\alpha(f)
 =|\mathbb S^{N-1}|\int_0^1r^{-\alpha}\,da(r),
 \qquad
 \partial_\alpha\log\bigl(\alpha E_\alpha(f)\bigr)
 =\int_0^1\log(1/r)\,d\mu_\alpha(r).
\]
The boundary term at zero vanishes since $a(r)=O(r^p)$, and the logarithmic
derivative is justified by the same estimate on compact subintervals of
$0<\alpha<p$. This proves \eqref{eq:radiallogmomentmargin}.
\end{proof}

\begin{remark}[Evaluation and size of the uniform margin]\label{rem:translationmarginvalue}
Every summand in \eqref{eq:translationmarginseries} is positive. The terms for $n=3,4$ alone give
\[
 \begin{aligned}
 \delta_p&\geq4^{-p}(\log2-1/2)
 +(3^{-p}-4^{-p})(\log(3/2)-1/3)\\
 &>4^{-p}(\log2-1/2).
 \end{aligned}
\]
For a partial sum $S_M$ through $n=M\geq3$, a directly usable error bound is
\[
 0<\delta_p-S_M\leq\frac{1+1/M}{pM^p}.
\]
Indeed the $n$th term is at most $n^{-p}/(n-1)$, and the integral test bounds the remaining series. The resulting values, rounded to six decimal places, are
\begin{center}
\begin{tabular}{crr}
\toprule
$p$ & $\delta_p$ & $2e^{\delta_p}$\\
\midrule
2 & 0.033420 & 2.067968\\
3 & 0.007044 & 2.014139\\
4 & 0.001758 & 2.003518\\
\bottomrule
\end{tabular}
\end{center}
The number $2e^{\delta_p}\Rmon$ is a proved sufficient diameter for the common affine class, not an evaluation of the unrestricted Schwarz threshold. No optimality of this margin within the affine class is claimed.
\end{remark}

}

\subsection{The precise estimate needed for the higher-dimensional threshold}
\begin{proposition}[A signed-moment criterion for the doubled threshold]
\label{prop:doubledcriterion}
For fixed $N,s,p$, the equality $\Dsch=2\Rmon$ is equivalent
to the following estimate for every $v\in\mathcal A_{N,s,p}$, with
$q=N+sp$:
\begin{equation}\label{eq:doubledbalance}
 \int_0^{1/2} M_v(r)r^{-q-1}\log\frac1{2r}\,dr
 \ge
 \int_{1/2}^1 M_v(r)r^{-q-1}\log(2r)\,dr.
\end{equation}
The estimate holds for every one-dimensional profile and for the
ball-superlevel class in every dimension. For arbitrary profiles in
$N\ge2$, it is not proved here.
\end{proposition}
\begin{proof}
The difference between the two sides is
$W_v(\Theta_s(v)-\log2)$, with $W_v>0$. Hence the estimate for every
profile is equivalent to $\theta_{N,s,p}\ge\log2$.
The opposite inequality follows from the translated-core
counterexamples. Formula \eqref{eq:exactschwarzthreshold} proves the
equivalence. The stated subclasses follow from their sharp comparison
theorems.
\end{proof}

The cumulative Riesz inequality gives $M_v(r)\ge0$ at each radius. In the exact deficit formula, this yields the sufficient bound $\Rmon$, but it does not compare the weighted contributions from $(0,1/2)$ and $(1/2,1)$ in \eqref{eq:doubledbalance}. The strict tilt argument orders the
moments of each fixed profile as $s$ varies; it also does not give the
missing uniform bound $\Theta_s(v)\ge\log2$ in higher dimensions.
The stronger estimate in Theorem~\ref{rem:affinethreshold} does cross
$\log2$, but its factorization relies on one common linear map applied
to a radial profile. Corollary~\ref{cor:affineexcluded} in fact rules out
that entire family as unrestricted minimizing sequences. Extending its
translation-moment estimate directly to the cumulative deficit would
require the latter to be a positive translation measure of the same
form; the present argument supplies no such representation.
Similarly, \eqref{eq:ballpotentialseries} controls nested levels when the
outer level is a ball, whereas arbitrary superlevel sets have no such
radial exterior potential. Applying that series to nonball outer levels
would therefore be unjustified. These are the precise limitations of
the two extensions considered here. The unconditional unrestricted
higher-dimensional conclusion remains
$\Rmon\le \Dsch\le2\Rmon$.

\begin{remark}[The remaining sharpness question for $N\ge2$]
In this remark assume $N\ge2$. For every fixed nonsymmetric profile, $\Theta(v)>0$: its first sign
reversal under dilation occurs strictly above $\Rmon$. A proof of exact
sharpness must therefore use a sequence of changing profiles, not dilations
of a single example. Condition \eqref{eq:missingschwarzconcentration}
would force the weighted deficits to concentrate near the largest available
pair distance. Indeed, for every $0<\delta<1$,
\[
 \nu_{v_j}((0,1-\delta])
 \le\frac{\Theta(v_j)}{\log(1/(1-\delta))}\longrightarrow0.
\]
Concentration in this probability sense alone would not suffice without
control of the logarithmic moment near zero. Neither the two-bump
polarization construction nor the transverse Steiner construction proves
\eqref{eq:missingschwarzconcentration}: Schwarz rearrangement also changes
the internal interactions of their components. Thus
\eqref{eq:exactschwarzthreshold} identifies the remaining geometric question
but does not establish $\theta_{N,s,p}=0$ or a positive lower bound for it.
The rigorous unconditional conclusion remains
$\Rmon\le \Dsch\le2\Rmon$.
\end{remark}

\begin{proposition}[A local obstruction to reducing the core construction]
\label{prop:coreobstruction}
In the notation of Proposition~\ref{prop:schwarzcore}, suppose instead that
$0<r_0<r_1<\Rmon$, and choose a sufficiently small nonzero radial
nonincreasing smooth core $w$. Then for every sufficiently small nonzero
translation $a$,
\[
 u_a^*=u_0,\qquad \Lsp(u_a^*)<\Lsp(u_a).
\]
Thus translating a tiny central core in such a background gives the usual
Schwarz comparison, rather than a counterexample.
\end{proposition}
\begin{proof}
Here $k_s'(\rho)<0$ throughout $[r_0,r_1]$. Consequently
\eqref{eq:tailhessian} gives $D^2T_\rho(0)>0$, uniformly in that compact
interval. By continuity there are $\eta,\kappa>0$ such that
$D^2T_\rho(x)\ge\kappa I$ for $|x|<\eta$ and
$r_0\le\rho\le r_1$. The positive weights $g_\alpha'$ in the proof of
Proposition~\ref{prop:schwarzcore} now give
$D^2F(a)\ge\kappa M_w I$ for a sufficiently small core and displacement.
Since $F$ is even, $\nabla F(0)=0$; hence
\[
 \Lsp(u_a)-\Lsp(u_0)
 \ge\frac{C_s\kappa M_w}{2p}|a|^2>0.
\]
The superlevel-set argument for $u_a^*=u_0$ is unchanged.
\end{proof}

\section{A bounded-domain consequence}
\begin{corollary}[P\'olya--Szeg\H{o} inequality on bounded domains]\label{cor:boundedPS}
Let $\Omega\subset\R^N$ be bounded and open, with
\[
 \diam\Omega\le\Rmon.
\]
Let $u:\Omega\to[0,\infty)$ have zero extension
\[
 \widetilde u(x)=
 \begin{cases}u(x),&x\in\Omega,\\0,&x\notin\Omega,\end{cases}
 \qquad \widetilde u\in\X.
\]
Let $\Omega^*$ be the centered ball with $|\Omega^*|=|\Omega|$, and write $u^*=(\widetilde u)^*$. Then
\[
 u^*\in X_0^{s,\log,p}(\Omega^*),\qquad
 \Lsp(u^*)\le\Lsp(\widetilde u).
\]
For nonzero $u$, equality holds if and only if
$\widetilde u=u^*(\cdot-a)$ almost everywhere for some $a\in\R^N$.
The same energy comparison holds for polarization and Steiner symmetrization, with transformed domains $\Omega^H$ and $\Omega^{\star_e}$, respectively, and the equality cases in Theorems~\ref{thm:polar} and~\ref{thm:steiner}.
\end{corollary}
\begin{proof}
Apply Theorems~\ref{thm:polar}, \ref{thm:steiner}, and \ref{thm:schwarz}
to the zero extension $\widetilde u$. Its support has diameter at most
$\diam\Omega$. Their finiteness and equality statements give the energy
claims. Inclusion preservation of the set rearrangements places each
transformed support in the corresponding transformed domain.
\end{proof}

\begin{remark}[The doubled one-dimensional domain bound]
When $N=1$, the Schwarz and Steiner conclusions of
Corollary~\ref{cor:boundedPS}, including rigidity, hold under the larger
bound $\diam\Omega\le2\Rmon$. This follows from
Theorem~\ref{thm:schwarzonedimensional} applied to the zero extension.
The polarization threshold remains $\Rmon$.
\end{remark}

\begin{remark}[Exterior interactions and the geometric hypothesis]
The corollary concerns the full zero exterior energy. In terms of the function on $\Omega$, it is
\[
 \begin{aligned}
 \Lsp(\widetilde u)
 ={}&\frac{C_s}{2p}\iint_{\Omega\times\Omega}
 |u(x)-u(y)|^p k_s(|x-y|)\,dx\,dy\\
 &+\frac{C_s}{p}\int_\Omega |u(x)|^p
 \left(\int_{\Omega^c}k_s(|x-y|)\,dy\right)dx.
 \end{aligned}
\]
All terms are absolutely integrable under the stated energy-space assumption, by \eqref{eq:absolute}. The exterior term is retained; an energy integrated only over $\Omega\times\Omega$ requires a separate argument. Membership in $L^p(\Omega)$ alone does not ensure that the zero extension belongs to $\X$.

The sufficient geometric hypothesis is a diameter bound. A bound on $|\Omega|$ alone does not imply it. For a ball $\Omega=B_r$, the hypothesis is equivalent to
\[
 |\Omega|\le\omega_N\left(\frac{\Rmon}{2}\right)^N.
\]
More generally, if $\diam\Omega\le A|\Omega|^{1/N}$ for a prescribed shape constant $A$, then $|\Omega|\le(\Rmon/A)^N$ is sufficient. The ambient domain may also be larger if the support of the zero extension itself has diameter at most $\Rmon$.
\end{remark}

\section{Additional rearrangement results}\label{sec:appendix}
\subsection{Generalized P\'olya--Szeg\H{o} inequalities}
\subsubsection{Strictly increasing transformations}
\begin{corollary}[Strictly increasing transforms]\label{cor:transform}
Let $A:[0,\infty)\to[0,\infty)$ be continuous and strictly increasing, with $A(0)=0$. Suppose $u\ge0$ vanishes outside a set of diameter at most $\Rmon$, and $A(u)\in\X$. If $V$ is Borel measurable and $V(u)\in L^1(\R^N)$, set
\[
 J(u)=\Lsp(A(u))+\int_{\R^N}V(u)\,dx.
\]
Then $J(Tu)\le J(u)$ for $T$ equal to a polarization, a Steiner symmetrization, or Schwarz rearrangement. The equality cases are exactly those of the corresponding theorem, expressed for $u$ itself.
\end{corollary}
\begin{proof}
Apply Theorem~\ref{thm:polar} for polarization, Theorem~\ref{thm:steiner}
for Steiner symmetrization, and Theorem~\ref{thm:schwarz} for Schwarz
rearrangement, each to $A(u)$ rather than $u$. Their hypotheses hold:
$A(u)\in\X$ by assumption, $A(u)\ge0$, and $A(u)$ has the same positivity
set as $u$ because $A(0)=0$ and $A$ is strictly increasing. Equation \eqref{eq:commute} and equimeasurability give the inequality. Injectivity of $A$ transfers every equality identity back to $u$. All rearrangements are defined because the support is bounded and, for $\alpha>0$, $|\{u>\alpha\}|\le A(\alpha)^{-p}\|A(u)\|_p^p$. Potential integrability on the full space in particular requires $V(0)=0$ in this bounded-support setting.
\end{proof}

\subsubsection{General two-variable interactions}
The next result treats interactions that cannot necessarily be reduced to a single change of the function values. Its algebraic assumption is the submodularity appropriate to an energy inequality; compare the supermodular interaction formulation in \cite{BH}.

\begin{theorem}[Generalized interaction inequality]\label{thm:generalG}
Let $G:[0,\infty)^2\to[0,\infty)$ be continuous and symmetric. Suppose that for some $M_G>0$,
\begin{equation}\label{eq:Ggrowth}
 0\le G(a,b)\le M_G|a-b|^p\qquad(a,b\ge0),
\end{equation}
and that for $a>b\ge0$ and $d>c\ge0$,
\begin{equation}\label{eq:Gmonge}
 G(a,c)+G(b,d)\ge G(a,d)+G(b,c).
\end{equation}
For nonnegative $u\in\X$ supported in a set of diameter at most $\Rmon$, define
\[
 \mathcal E_G(u)=\frac{C_s}{2p}\iint_{\R^{2N}}
 G(u(x),u(y))k_s(|x-y|)\,dx\,dy.
\]
Then the integral is absolutely convergent, and
\begin{equation}\label{eq:allG}
 \mathcal E_G(u^H)\le\mathcal E_G(u),\qquad
 \mathcal E_G(u^{\star_e})\le\mathcal E_G(u),\qquad
 \mathcal E_G(u^*)\le\mathcal E_G(u).
\end{equation}
If \eqref{eq:Gmonge} is strict for every such pair of strict orderings, the equality cases for nonzero $u$ are respectively
\[
 u^H=u\ \text{or}\ u^H=u\circ\sigma_H;\qquad
 u=u^{\star_e}(\cdot-ae);\qquad
 u=u^*(\cdot-a),
\]
with $a\in\R$ in the Steiner case and $a\in\R^N$ in the Schwarz case. Without strictness, these configurations still give equality, but additional equality cases may occur.
\end{theorem}
\begin{proof}
The bound \eqref{eq:Ggrowth} and \eqref{eq:absolute} give absolute
convergence for $u$ and its transforms. Membership of those transforms in
$\X$ follows, respectively, from Theorems~\ref{thm:polar},
\ref{thm:steiner}, and \ref{thm:schwarz}, whose common hypotheses
$u\ge0$, $u\in\X$, and $\diam S\le\Rmon$ are exactly those assumed here. Define
\[
 Q_G=G(a,c)+G(b,d)-G(a_+,c_+)-G(a_-,c_-),
\]
using the four values and their sorted pairs from Lemma~\ref{lem:Q}. It vanishes for equally ordered pairs and is nonnegative for oppositely ordered pairs by \eqref{eq:Gmonge}. Under strictness it is positive exactly when $(a-b)(c-d)<0$.

For an integrable radial kernel the four-quadrant calculation is purely algebraic and gives the analogue of \eqref{eq:defectK}. Apply it to $k_\varepsilon$ and pass to the limit using \eqref{eq:Ggrowth} and the finite power energies. Equivalently, on the active set one may use monotone convergence of the nonnegative kernel differences. This yields
\[
 \mathcal E_G(u)-\mathcal E_G(u^H)
 =\frac{C_s}{2p}\int_{H\times H}[k_s(r)-k_s(\rho)]Q_G\,dx\,dy.
\]
If $Q_G>0$, the orderings are opposite, so the same nonnegativity argument as in \eqref{eq:activedistance} gives $r\le\rho\le\diam S\le\Rmon$. The integrand is nonnegative. Under strictness, positive-measure sets of both orderings make it strictly positive, proving the polarization equality characterization.

For Steiner and Schwarz symmetrization, apply Theorem~\ref{thm:polarapprox} to the nonnegative $L^p$ function $u$ with the respective target $T=\star_e$ or $T=*$. The generalized energies decrease. To justify the limit, split $k_s=P+K$ as in Lemma~\ref{lem:lsc}. The $P$-energy with integrand $G$ is lower semicontinuous by continuity, nonnegativity, and Fatou's lemma. The $K$-energy is continuous under strong $L^p$ convergence. Indeed, along an almost-everywhere convergent subsequence, its integrands converge pointwise, and their absolute values are bounded by a constant times $|v_n(x)|^p+|v_n(y)|^p$. These majorants converge in $L^1(\R^{2N},|K(|x-y|)|\,dx\,dy)$ when $v_n\to v$ in $L^p$. The generalized dominated-convergence argument therefore gives convergence of the integrals. This proves \eqref{eq:allG}.

If the endpoint energies are equal and \eqref{eq:Gmonge} is strict, every intermediate polarization has equality. The compositions are thus Euclidean motions of the original function. The compactness argument in Theorem~\ref{thm:steiner}, or Theorem~\ref{thm:schwarz}, then gives exactly the stated translations. The converse follows from invariance.
\end{proof}

\begin{remark}[Examples and potentials]
The theorem covers
\[
 G(a,b)=\sum_{j=1}^m\lambda_j|A_j(a)-A_j(b)|^p,
 \qquad \lambda_j>0,
\]
when the $A_j$ are Lipschitz and strictly increasing. The growth bound follows from Lipschitz continuity, while strict submodularity follows by applying \eqref{eq:strictQ} to the transformed values. The $A_j$ need not be the same function.

A nonhomogeneous example is $G(a,b)=\Psi(|a-b|)$ with
\[
 \Psi(t)=(1+t^p)^{1/p}-1,\qquad
 0\le\Psi(t)\le t^p/p,
 \qquad
 \Psi''(t)=(p-1)t^{p-2}(1+t^p)^{1/p-2}>0\quad(t>0).
\]
The even function $r\mapsto\Psi(|r|)$ is strictly convex, so its difference interaction is strictly submodular. The theorem also applies to $v=A(u)$ whenever $A$ is continuous and strictly increasing, $A(0)=0$, and $A(u)\in\X$. Finally, adding $\int V(u)$ with $V(u)\in L^1$ changes neither the comparisons nor the equality cases, by equimeasurability. These statements concern the specified interaction class, not arbitrary nonlocal Lagrangians.
\end{remark}

\subsubsection{General distance-dependent integrands}\label{sec:generalPhi}
We next allow the distance and function values to enter the integrand jointly. Consider
\begin{equation}\label{eq:PhiEnergy}
 \mathcal I_\Phi(u)=\iint_{\R^{2N}}
 \Phi(|x-y|,u(x),u(y))\,dx\,dy.
\end{equation}
The scalar integrand $\Phi(r,a,b)$ may change sign and need not be a product of a function of $r$ with a function of $a,b$. The class remains invariant under Euclidean motions because the spatial variables enter through $|x-y|$. Arbitrary dependence on $x,y$ would require additional geometric hypotheses and is not included in this theorem. For the role and necessity of structural conditions in the classical generalized rearrangement setting, compare \cite{BH,HNecessity}.

\paragraph{The sorting-defect condition}
For $a>b\ge0$ and $d>c\ge0$, define the distance-dependent Monge gap
\begin{equation}\label{eq:MongePhi}
 \Delta_\Phi(r;a,b,c,d)
 =\Phi(r,a,c)+\Phi(r,b,d)-\Phi(r,a,d)-\Phi(r,b,c).
\end{equation}
The essential condition is that $r\mapsto\Delta_\Phi(r;a,b,c,d)$ be nonincreasing on the relevant distance interval. The gap itself need not be positive. This distinction is essential for signed kernels: even the logarithmic model need not be submodular in the values at every distance, but its Monge gap decreases on the admissible range.

\begin{theorem}[Distance-dependent generalized P\'olya--Szeg\H{o} inequality]\label{thm:Phi}
Fix $\ell>0$, and put
\[
 w_s(r)=(1+|\log r|)r^{-N-sp}.
\]
Assume that $\Phi:(0,\infty)\times[0,\infty)^2\to\R$ has the following properties:
\begin{enumerate}[label=\textup{(\roman*)},leftmargin=*]
\item $\Phi$ is measurable in $r$, continuous in $(a,b)$ for every $r>0$, and symmetric in $a,b$.
\item For some $M>0$ and all $r>0$, $a,b\ge0$,
\begin{equation}\label{eq:Phigrowth}
 |\Phi(r,a,b)|\le M w_s(r)|a-b|^p.
\end{equation}
\item There is a measurable $h\ge0$ with $h(|\cdot|)\in L^1(\R^N)$ such that
\begin{equation}\label{eq:Philower}
 \Phi(r,a,b)\ge-h(r)|a-b|^p.
\end{equation}
\item For every $a>b\ge0$ and $d>c\ge0$, and every $0<r<\rho\le \ell$,
\begin{equation}\label{eq:Phimon}
 \Delta_\Phi(r;a,b,c,d)\ge\Delta_\Phi(\rho;a,b,c,d).
\end{equation}
\end{enumerate}
Let $u\ge0$ belong to $\X$ and vanish outside a measurable set $S$ with $\diam S\le \ell$. Then all the following integrals converge absolutely and
\begin{equation}\label{eq:PhiAll}
 \mathcal I_\Phi(u^H)\le\mathcal I_\Phi(u),\qquad
 \mathcal I_\Phi(u^{\star_e})\le\mathcal I_\Phi(u),\qquad
 \mathcal I_\Phi(u^*)\le\mathcal I_\Phi(u).
\end{equation}
If \eqref{eq:Phimon} is strict for every choice of strict value orderings and every $r<\rho$, equality for nonzero $u$ occurs exactly in the configurations
\[
 u^H=u\ \text{or}\ u^H=u\circ\sigma_H;\qquad
 u=u^{\star_e}(\cdot-ae);\qquad
 u=u^*(\cdot-a),
\]
respectively, with one constant translation parameter in each case.
\end{theorem}

\begin{proof}
\emph{Finiteness and membership.} The positive kernel
\[
 W(r)=r^{-N-sp}\log(1/r)\,\mathbf1_{(0,1)}(r)
\]
is radially nonincreasing. Apply Lemma~\ref{lem:positivekernel} to this $W$. Its energy equals the
defining logarithmic seminorm, which is finite because $u\in\X$; it is
nonnegative, decreasing, finite at each positive radius, and has integrable
tails. The lemma proves that all three transforms have finite $W$-energy,
hence belong to $\X$, without a support-diameter restriction for this
auxiliary positive-kernel step. Equation \eqref{eq:absolute}, followed by \eqref{eq:Phigrowth}, proves absolute convergence of \eqref{eq:PhiEnergy} for $u$ and its transforms.

\emph{Polarization.} Use the notation $a,b,c,d,r,\rho$ of Theorem~\ref{thm:polar}. If the two value pairs are ordered in the same direction, their grouped contribution is unchanged by sorting. If $a>b$ and $d>c$, the difference between the original and polarized grouped contributions is
\begin{equation}\label{eq:Phidefect}
 \Delta_\Phi(r;a,b,c,d)-\Delta_\Phi(\rho;a,b,c,d).
\end{equation}
The reverse strict ordering is reduced to this case by interchanging $x$ and $y$, using symmetry of $\Phi$ in its value variables. Nonnegativity of $u$ implies that the original occupied cross pair lies in $S$, so $r\le\rho\le\diam S\le \ell$. Hence \eqref{eq:Phidefect} is nonnegative by \eqref{eq:Phimon}. Splitting into the four reflected quadrants and integrating proves the first inequality in \eqref{eq:PhiAll}. The identity is legitimate directly: its grouped expression consists of the original and polarized energy terms, each absolutely integrable. It is not necessary to integrate the two Monge gaps separately, which might produce divergent intermediate expressions.

If both strict-ordering sets $A$ and $B$ from Theorem~\ref{thm:polar} have positive measure, then $r<\rho\le \ell$ almost everywhere on $A\times B$. Strictness in \eqref{eq:Phimon} makes the energy defect strictly positive. Thus equality forces $|A|=0$ or $|B|=0$: respectively,
$u^H=u\circ\sigma_H$ or $u^H=u$ almost everywhere. Reflection invariance
proves the converse.

\emph{Lower semicontinuity.} Define
\[
 \widehat\Phi(r,a,b)=\Phi(r,a,b)+h(r)|a-b|^p\ge0.
\]
Its energy is lower semicontinuous under strong $L^p$ convergence by continuity in the values and Fatou's lemma along an almost-everywhere convergent subsequence. The correction $E_h(v)$ is continuous under strong $L^p$ convergence by the bounded difference-map argument in Lemma~\ref{lem:lsc}. Therefore $\mathcal I_\Phi=\mathcal I_{\widehat\Phi}-E_h$ is lower semicontinuous, with the usual extended-value interpretation of its nonnegative part.

\emph{Steiner and Schwarz symmetrization.} Approximate the chosen symmetrization by successive polarizations. Their support diameters are at most $\ell$ by Lemma~\ref{lem:diameter}, so every step decreases $\mathcal I_\Phi$. Lower semicontinuity proves the remaining inequalities in \eqref{eq:PhiAll}. Under strictness, equality of the endpoint energies forces equality at every step. Each step is therefore an identity or a reflection of the whole function. The Euclidean-motion compactness arguments in Theorems~\ref{thm:steiner} and \ref{thm:schwarz} give the asserted translations. Their energy invariance proves sufficiency.
\end{proof}

\paragraph{A differential criterion}
\begin{corollary}\label{cor:Phidifferential}
Fix $\ell>0$ and assume properties \textup{(i)}--\textup{(iii)} of
Theorem~\ref{thm:Phi}: measurability in distance, continuity and symmetry in
values, the growth bound \eqref{eq:Phigrowth}, and the integrable lower
bound \eqref{eq:Philower}. Suppose in addition that $\Phi$ is three times
continuously differentiable on $(0,\ell)\times(0,\infty)^2$, and that the
rectangle identities for its mixed derivatives extend to rectangles meeting
a zero value by limits from positive values. Assume also that each gap $r\mapsto\Delta_\Phi(r;a,b,c,d)$ is left-continuous at $\ell$. A sufficient condition for \eqref{eq:Phimon} is
\begin{equation}\label{eq:mixedPhi}
 \partial_r\partial_a\partial_b\Phi(r,\alpha,\beta)\ge0
 \quad(0<r<\ell,\ \alpha,\beta>0).
\end{equation}
If this mixed derivative is strictly positive almost everywhere on $(0,\ell)\times(0,\infty)^2$, then \eqref{eq:Phimon} is strict for every strict value rectangle and $0<r<\rho\le \ell$, and the rigidity conclusion of Theorem~\ref{thm:Phi} applies.
\end{corollary}
\begin{proof}
For $a>b$ and $d>c$,
\[
 \Delta_\Phi(r;a,b,c,d)
 =-\int_b^a\int_c^d\partial_a\partial_b\Phi(r,\alpha,\beta)\,d\beta\,d\alpha.
\]
Differentiating in $r$ gives a nonpositive derivative under \eqref{eq:mixedPhi}. Equivalently, integrate the mixed derivative over $(r,\rho)\times(b,a)\times(c,d)$; a strictly positive integral makes the gap strictly decrease. Conditions on the derivatives are only sufficient; Theorem~\ref{thm:Phi} itself uses the finite-difference condition and permits nonsmooth integrands.
\end{proof}

\paragraph{Examples and the logarithmic model}
\begin{corollary}[Sums of signed logarithmic interactions]\label{cor:Phisums}
Let $q=N+sp$ and, for $j=1,\ldots,m$, let
\[
 k_j(r)=(\beta_j-\gamma_j\log r)r^{-q},\qquad
 \gamma_j>0,\qquad \lambda_j>0.
\]
Let $G_j$ be continuous, symmetric, nonnegative, satisfy the power-growth condition \eqref{eq:Ggrowth}, and be submodular as in \eqref{eq:Gmonge}. Then
\[
 \Phi(r,a,b)=\sum_{j=1}^m\lambda_j k_j(r)G_j(a,b)
\]
satisfies Theorem~\ref{thm:Phi} whenever
\[
 \ell\le\min_{1\le j\le m}\exp(\beta_j/\gamma_j+1/q).
\]
If at least one $G_j$ is strictly submodular, the three equality classifications also hold.
\end{corollary}
\begin{proof}
The growth bound is immediate. Each negative kernel part begins at a positive radius and is integrable on $\R^N$; its multiple by the growth constant for $G_j$ supplies the function $h$ in \eqref{eq:Philower}. Differentiation shows $k_j$ is strictly decreasing up to $\exp(\beta_j/\gamma_j+1/q)$. Moreover,
\[
 \Delta_\Phi(r)=\sum_j\lambda_j k_j(r)\Delta_{G_j},
\]
where $\Delta_{G_j}\ge0$ for the strict opposite orderings. Thus the gap decreases on $(0,\ell]$ and decreases strictly if one $\Delta_{G_j}>0$.
\end{proof}

\begin{remark}
The original energy is recovered by $m=1$, $\beta_1=b_s$, $\gamma_1=p$, $\lambda_1=C_s/(2p)$, and $G_1(a,b)=|a-b|^p$. Taking distinct $G_j$ and distinct ratios $\beta_j/\gamma_j$ generally yields an integrand that does not factor into one distance kernel and one value interaction.

One may also add $K(r)G_0(a,b)$ with a nonnegative nonincreasing kernel $K$ dominated by a constant multiple of $w_s$, and a nonnegative submodular $G_0$ satisfying the same power-growth bound. Its Monge-gap contribution is nonincreasing, so it preserves the inequalities. A strictly decreasing contribution with a strictly submodular interaction supplies rigidity. Integrable potential terms depending only on $u$ are preserved by equimeasurability.
\end{remark}

\paragraph{Limits of the generalization}
Submodularity in the value variables at each fixed distance, by itself, does not control how an interaction changes when the distance is shortened. The operative condition for the general integrand is \eqref{eq:Phimon}, involving both the distances and the values. Nor does the theorem imply a result for arbitrary spatially dependent weights: such weights can destroy reflection or translation invariance and alter the equality cases. The power-growth and integrable-negative-part hypotheses ensure finiteness and passage to the symmetrized limit; they are sufficient assumptions, not claimed necessary ones. If the monotonicity condition holds for every $\ell>0$, the inequalities hold without a support-size restriction for boundedly supported functions. The signed logarithmic kernel does not have that global monotonicity; Propositions~\ref{prop:sharpdiameter}, \ref{prop:sharpsteiner}, and \ref{prop:schwarzcore} give the respective counterexamples.

\subsection{Functions of the total energy and of its pointwise density}\label{sec:density}
An expression such as $j(u,\text{energy})$ has two distinct meanings. The total energy is a scalar, while a pointwise energy density is a function of $x$. These interpretations require different arguments. Neither should be confused with applying $j$ to the operator $\mathfrak L_{s,p}u$ itself.

\subsubsection{Dependence on the total energy}
\begin{proposition}[Monotone dependence on the total scalar energy]\label{prop:totalj}
Let $N\ge1$, $1<p<\infty$, $0<s<1$, and let $T$ be a polarization,
one-direction Steiner symmetrization, or Schwarz rearrangement. Suppose
$u\ge0$ belongs to $\X$ and vanishes outside a set of diameter at most
$\Rmon$. Theorems~\ref{thm:polar}, \ref{thm:steiner}, and
\ref{thm:schwarz} then give the respective comparison. Set $E=\Lsp(u)$ and $E_T=\Lsp(Tu)\le E$. Let $\Omega$ have finite positive measure, let $u$ vanish outside $\Omega$, and denote its transformed domain by $T\Omega$. Suppose $j:[0,\infty)\times\R\to\R$ is Borel measurable and all displayed integrals are finite. Define
\[
 J_\Omega(u)=\int_\Omega j(u(x),\Lsp(u))\,dx.
\]
If $t\mapsto j(a,t)$ is nondecreasing for every $a\ge0$, then
\[
 J_{T\Omega}(Tu)\le J_\Omega(u).
\]
If it is nonincreasing for every $a$, the comparison is reversed. Convexity in $t$ is not needed. If the relevant monotonicity is strict for every $a$, equality is equivalent to $E_T=E$, and consequently has the equality configuration in Theorem~\ref{thm:polar}, \ref{thm:steiner}, or \ref{thm:schwarz}, according to the chosen $T$. The same conclusions hold for $j(\|u\|_p,\Lsp(u))$.
\end{proposition}
\begin{proof}
The distributions of $Tu$ on $T\Omega$ and of $u$ on $\Omega$, including the measures of their zero sets, agree. Thus
\[
 J_{T\Omega}(Tu)=\int_\Omega j(u(x),E_T)\,dx.
\]
Compare this with $\int_\Omega j(u(x),E)\,dx$ using $E_T\le E$. Strict monotonicity makes the comparison strict if $E_T<E$, because $|\Omega|>0$. The norm-dependent formulation uses $\|Tu\|_p=\|u\|_p$.
\end{proof}

\subsubsection{Definition of the pointwise density}
For $u\in\X$, define, at almost every $x$,
\begin{equation}\label{eq:localdensity}
 e_u(x)=\frac{C_s}{2p}\int_{\R^N}|u(x)-u(y)|^p k_s(|x-y|)\,dy.
\end{equation}
Absolute convergence of the double integral and Fubini's theorem give
\[
 e_u\in L^1(\R^N),\qquad \int_{\R^N}e_u(x)\,dx=\Lsp(u).
\]
The density may be signed. For the functional
\begin{equation}\label{eq:Jdensity}
 J_j(u)=\int_{\R^N}j(u(x),e_u(x))\,dx,
\end{equation}
the second variable of $j$ must therefore be defined on the relevant subset of $\R$. A fractional power of $e_u$ is not automatically a real-valued quantity.

The inequality for $\int e_u$ supplies neither a pointwise comparison between $e_{Tu}$ and $e_u$ nor the joint-distribution comparison involving $Tu$ and $e_{Tu}$ needed in \eqref{eq:Jdensity}. Generalized local-gradient inequalities, as studied in \cite{Hajaiej10}, do not by themselves establish a theorem for this nonlocal density. A genuine additional hypothesis controlling the coupling is required. The counterexample in Appendix~\ref{sec:nonlinearcounter} makes this obstruction explicit even when the support is small.

\subsubsection{A sufficient weighted class and both comparison directions}
\begin{theorem}[Affine density dependence]\label{thm:affinedensity}
Let $u\ge0$ belong to $\X\cap L^\infty(\R^N)$, with support diameter at most $\Rmon$. Put $M=\|u\|_\infty$. Assume $A\in C^1([0,M])$ is nonnegative and convex, and let $V:[0,M]\to\R$ be Borel measurable with $V(u)\in L^1(\R^N)$. For
\[
 j_+(a,t)=V(a)+A(a)t,\qquad
 j_-(a,t)=V(a)-A(a)t,
\]
the corresponding full-space integrals are finite and, for all three symmetrizations,
\begin{equation}\label{eq:affineboth}
 J_{j_+}(Tu)\le J_{j_+}(u),\qquad
 J_{j_-}(Tu)\ge J_{j_-}(u).
\end{equation}
If $A>0$ on $[0,M]$ and $u\ne0$, equality in either comparison has precisely the corresponding polarization, Steiner, or Schwarz configuration.
\end{theorem}
\begin{proof}
Symmetry in $x,y$ gives
\begin{equation}\label{eq:weightedsym}
 \int A(u(x))e_u(x)\,dx
 =\frac{C_s}{2p}\iint G_A(u(x),u(y))k_s(|x-y|)\,dx\,dy,
\end{equation}
where
\[
 G_A(a,b)=\frac{A(a)+A(b)}2|a-b|^p.
\]
On the bounded value range this has the power-growth bound required in Theorem~\ref{thm:generalG}. For $a>b$, its mixed derivative is
\begin{equation}\label{eq:weightedmixed}
 \partial_a\partial_bG_A(a,b)
 =-\frac p2(a-b)^{p-2}
 \left[(p-1)(A(a)+A(b))+(a-b)(A'(a)-A'(b))\right].
\end{equation}
The bracket is nonnegative because $A\ge0$ and $A'$ is nondecreasing. Integrating this mixed derivative on a value rectangle gives submodularity. For $1<p<2$, the diagonal singularity is locally integrable, so the same rectangle identity holds. If $A>0$, the mixed derivative is strictly negative away from the diagonal, giving strict submodularity. One may restrict the interaction theorem to $[0,M]$, since all transforms preserve the essential supremum. Apply Theorem~\ref{thm:generalG} to \eqref{eq:weightedsym} on the invariant
value interval $[0,M]$. Continuity, symmetry, nonnegativity, power growth,
and submodularity of $G_A$ have just been checked; $u\ge0$, $u\in\X$,
and $\diam(\esssupp u)\le\Rmon$ are assumptions of the present theorem.
Equimeasurability gives equality of the integrals of $V$. Changing the sign of the weighted term reverses the inequality and leaves equality unchanged.
\end{proof}

\begin{remark}
Thus an increase under symmetrization is obtained for $j_-$, while the usual P\'olya--Szeg\H{o} decrease is obtained for $j_+$. Convexity of $A$ is a sufficient coupling condition, not a claimed necessary condition. The theorem does not assert the result for arbitrary nonnegative weights $A$, even though $t\mapsto A(a)t$ is nondecreasing for every fixed $a$.
\end{remark}

\subsection{Failure for nonlinear functions of the pointwise density}\label{sec:nonlinearcounter}
The following counterexample shows that smoothness, strict monotonicity, and strict convexity with respect to the energy-density variable do not, on their own, imply a generalized P\'olya--Szeg\H{o} inequality. The counterexample has small support, so it is independent of the large-distance obstruction for the full logarithmic kernel.

\begin{theorem}[Failure for nonlinear pointwise-density functionals]\label{thm:nonlinearcounter}
Fix $N\ge1$, $1<p<\infty$, and $0<s<1$ with $sp<1$. There exist a nonnegative $u\in C_c^\infty(\R^N)$ and a smooth integrand $j:[0,\infty)\times\R\to\R$ such that
\[
 \diam(\esssupp u)\le\Rpos/2<\Rmon,
 \qquad
 \partial_tj(a,t)>0,\quad \partial_{tt}j(a,t)>0
 \quad(a\ge0,\ t\in\R),
\]
but nevertheless
\begin{equation}\label{eq:nonlinearfailure}
 \displaystyle J_j(u^*)>J_j(u),\qquad
 J_j(v)=\int_{\R^N}j(v(x),e_v(x))\,dx.
\end{equation}
Both integrals are finite. The integrand is genuinely nonlinear in its second variable. No smoothness of $u^*$ is needed or asserted.
\end{theorem}

\begin{proof}
\emph{Step 1: geometry and admissibility.} Set $r=\Rpos/24$ and define
\[
 H=B_r(4re_1),\qquad L=B_r(-4re_1),\qquad
 H_*=B_r(0),\qquad
 L_*=B_{2^{1/N}r}(0)\setminus B_r(0).
\]
All four sets have the relevant volumes $|H|=|L|=|H_*|=|L_*|=:v$. For a height $M>1$ to be chosen, put
\begin{equation}\label{eq:twolvl}
 u=M\mathbf1_H+\mathbf1_L,\qquad
 u^*=M\mathbf1_{H_*}+\mathbf1_{L_*}.
\end{equation}
The rearrangement formula follows by comparing superlevel volumes below and above the level $1$. The intermediate piecewise constant function has support diameter $10r=5\Rpos/12<\Rpos/2$.

These functions belong to $\X$. Indeed, for each ball or annulus $E$ used here,
\[
 |E\mathbin\triangle(E-z)|\le C_E|z|\quad(|z|<1),
\]
and the seminorm of $\mathbf1_E$ is bounded near zero by a constant times
\[
 \int_0^1 t^{-sp}\log(1/t)\,dt<\infty.
\]
The $L^p$ norm is finite and the remaining distances are controlled by integrable kernels. Finite linear combinations are also in $\X$.

\emph{Step 2: a strictly larger interaction after rearrangement.} Define
\[
 I(E,F)=\int_E\int_F k_s(|x-y|)\,dy\,dx,
 \qquad
 P_k(E)=\int_E\int_{E^c}k_s(|x-y|)\,dy\,dx.
\]
All these integrals for the chosen sets converge absolutely by the preceding admissibility argument; the touching boundary between $H_*$ and $L_*$ causes no divergence when $sp<1$. Distances between $H,L$ lie in $[6r,10r]$, while distances between $H_*,L_*$ are at most $(1+2^{1/N})r\le3r$. All these distances are in the decreasing region. Hence
\begin{equation}\label{eq:Qdensity}
 \mathcal Q:=I(L_*,H_*)-I(L,H)
 \ge v^2[k_s(3r)-k_s(6r)]>0.
\end{equation}
Let $\mathcal B=P_k(L_*)-P_k(L)$, a finite real number, and choose
\begin{equation}\label{eq:Mdensity}
 M=1+\left(2+\frac{|\mathcal B|}{\mathcal Q}\right)^{1/p}.
\end{equation}
Then $(M-1)^p-1=1+|\mathcal B|/\mathcal Q$.

\emph{Step 3: select the intermediate value.} Define the smooth nonnegative polynomial
\[
 A(a)=\frac{a^2(a-M)^2}{(M-1)^2},
 \qquad A(0)=A(M)=0,\quad A(1)=1.
\]
It weights only the level $1$ among the three values of $u$. In particular,
\[
 J_0(u):=\int A(u(x))e_u(x)\,dx=\int_L e_u(x)\,dx.
\]
For $x\in L$, the value differences are zero on $L$, $(M-1)^p$ on $H$, and $1$ on the zero set. Thus
\[
 J_0(u)=\frac{C_s}{2p}
 \left[P_k(L)+((M-1)^p-1)I(L,H)\right],
\]
and the same formula holds for $u^*$ with the starred sets. It follows that
\begin{align}
 \Delta&:=J_0(u^*)-J_0(u)\notag\\
 &=\frac{C_s}{2p}
 \left[\mathcal B+((M-1)^p-1)\mathcal Q\right]\notag\\
 &=\frac{C_s}{2p}(\mathcal B+|\mathcal B|+\mathcal Q)
 \ge\frac{C_s}{2p}\mathcal Q>0.\label{eq:Deltadensity}
\end{align}
Thus failure already occurs for a nonnegative weight that does not satisfy the sufficient convex-coupling hypothesis of Theorem~\ref{thm:affinedensity}.

\emph{Step 4: a strictly convex nonlinear perturbation.} Put
\[
 \phi(t)=\log(1+e^t)-\log2.
\]
This function satisfies $\phi(0)=0$, $|\phi(t)|\le|t|$, and
\[
 \phi'(t)=\frac{e^t}{1+e^t}>0,\qquad
 \phi''(t)=\frac{e^t}{(1+e^t)^2}>0.
\]
Since both densities are integrable, the number
\[
 K=\int_{\R^N}\phi(e_{u^*}(x))\,dx
 -\int_{\R^N}\phi(e_u(x))\,dx
\]
is finite. Define
\begin{equation}\label{eq:jcounter}
 \eta=\frac{\Delta}{2(1+|K|)}>0,\qquad
 j(a,t)=A(a)t+\eta\phi(t).
\end{equation}
These definitions are noncircular: the sets, $M$, $u$, $A$, $\Delta$, and $K$ are fixed before $\eta$ and $j$ are chosen. For all $a\ge0$, $t\in\R$,
\[
 \partial_tj(a,t)=A(a)+\eta\frac{e^t}{1+e^t}>0,
 \qquad
 \partial_{tt}j(a,t)=\eta\frac{e^t}{(1+e^t)^2}>0.
\]
At the values attained by $u$ or $u^*$, $A$ is bounded and $|j(a,t)|\le(A(a)+\eta)|t|$, so the functionals are finite. Finally,
\begin{align*}
 J_j(u^*)-J_j(u)
 &=\Delta+\eta K\\
 &\ge\Delta-\frac{\Delta|K|}{2(1+|K|)}
 >\frac\Delta2>0.
\end{align*}
This proves \eqref{eq:nonlinearfailure} for the intermediate function.

\emph{Step 5: a smooth test function with the same strict failure.}
Let $\rho_\varepsilon$ be a nonnegative smooth radial mollifier supported in $B_\varepsilon$, and set $U_\varepsilon=\rho_\varepsilon*u$. Then $0\le U_\varepsilon\le M$, $U_\varepsilon\to u$ in $L^p$, and their total variations are uniformly bounded. They have support diameter at most $5\Rpos/12+2\varepsilon\le\Rpos/2$ for sufficiently small $\varepsilon$.
The rearrangement contraction gives $U_\varepsilon^*\to u^*$ in $L^p$. These rearrangements are radial decreasing, bounded by $M$, and supported in a common ball $B_B$. Their total variations are also uniformly bounded: the coarea formula gives
\[
 \operatorname{TV}(U_\varepsilon^*)
 =\int_0^M P(\{U_\varepsilon^*>t\})\,dt
 \le M|\mathbb S^{N-1}|B^{N-1}.
\]
For $N=1$, each nonempty bounded interval has perimeter two.

For either sequence $f_\varepsilon=U_\varepsilon,U_\varepsilon^*$, with limit $f=u,u^*$, the BV translation estimate gives
\[
 \int|f_\varepsilon(x+h)-f_\varepsilon(x)|^p\,dx
 \le (2M)^{p-1}|h|\operatorname{TV}(f_\varepsilon).
\]
The same bound holds for $f$. Thus the absolute-kernel energy on $|h|<\delta$ is uniformly bounded by
$C\int_0^\delta r^{-sp}(1+|\log r|)\,dr\to0$, using $sp<1$. On $|h|\ge\delta$ the kernel is integrable; strong $L^p$ convergence and the bounded difference map of Lemma~\ref{lem:lsc} give convergence of the $p$th powers of the differences in the weighted $L^1$ space. Combining the two ranges and using Fubini yields
\[
 \|e_{U_\varepsilon}-e_u\|_1\to0,\qquad
 \|e_{U_\varepsilon^*}-e_{u^*}\|_1\to0.
\]
The polynomial $A$ is bounded on $[0,M]$, and $\phi$ in Step 4 is globally one-Lipschitz. Along a common almost-everywhere convergent subsequence, dominated convergence and these density limits imply
\[
 J_j(U_\varepsilon)\to J_j(u),\qquad
 J_j(U_\varepsilon^*)\to J_j(u^*).
\]
For example, split the weighted difference into
$A(U_\varepsilon)(e_{U_\varepsilon}-e_u)$ and
$(A(U_\varepsilon)-A(u))e_u$; the first tends to zero in integral and the second is dominated by a constant times $|e_u|$. Lipschitz continuity handles $\phi$.
The positive gap in Step 4 therefore persists for sufficiently small $\varepsilon$. Take $U_\varepsilon$ as the asserted $u$, keeping the same smooth $j$. Its strict density monotonicity and convexity are unchanged, and the bounded weight and integrable densities ensure both functionals are finite.
\end{proof}

\begin{remark}[The quadratic case discussed above]
For $p=2$ the construction applies for $0<s<1/2$ and \eqref{eq:Mdensity} reads $M=1+\sqrt{2+|\mathcal B|/\mathcal Q}$. The coefficient in Step 3 becomes $M^2-2M$, and \eqref{eq:Deltadensity} has prefactor $C_s/4$. Thus the quadratic example is a special case of the same construction for every $p$ with $sp<1$.
\end{remark}

\begin{remark}[Sign, regularity, and the positive-part density]
For the full signed density, $j$ is defined on all real second arguments and need not be nonnegative there. The assertion is that it is smooth, strictly increasing, and strictly convex in that argument. The failure involves coupling to the first variable and does not contradict Theorem~\ref{thm:affinedensity}, whose weight is required to be convex. Nor does it show that every nonlinear dependence on the density must fail.

The same construction also works for the nonnegative density obtained by replacing $k_s$ with $k_s^+=\max\{k_s,0\}$: recompute $\mathcal B$, $\mathcal Q$, $M$, and the subsequent constants for that kernel. Its small-distance ordering still gives \eqref{eq:Qdensity}, and every subsequent identity is unchanged. In that version the density arguments are nonnegative and $j(a,t)\ge0$ for $a,t\ge0$. Thus the obstruction is not solely a consequence of negative density values. This statement does not alter the unrestricted rearrangement inequality for the positive-part \emph{linear} energy itself.
\end{remark}

\subsubsection{Scope of the density counterexample}
The construction shows that strict increase and convexity in the pointwise energy-density variable are insufficient for a generalized Schwarz P\'olya--Szeg\H{o} inequality, even on small supports. Its test function is smooth and compactly supported, all energies are finite, and $sp<1$ ensures that the uniform BV estimate is integrable in the smoothing step. The construction leaves the regime $sp\ge1$ open. A broader theorem must control the joint behavior of $u$ and $e_u$; Theorem~\ref{thm:affinedensity} provides one sufficient convex-weight condition. For comparison, Proposition~\ref{prop:totalj} concerns monotone composition with the \emph{total scalar energy}, a different functional.

\subsection{Signed functions and zero exterior domains}
\begin{proposition}[Reduction to absolute values]\label{prop:modulus}
Let real-valued $u\in\X$ vanish outside $S$, with $\diam S\le\Rpos$. Then
\[
 \Lsp(|u|)\le\Lsp(u).
\]
For $u\ne0$, equality holds precisely when $u$ has a constant sign almost everywhere on its nonzero set.
\end{proposition}
\begin{proof}
The modulus belongs to $\X$ by $||a|-|b||\le|a-b|$. The energy difference has integrand
\[
 \bigl(|u(x)-u(y)|^p-\bigl||u(x)|-|u(y)|\bigr|^p\bigr)k_s(|x-y|).
\]
It vanishes if either point is outside $S$. On $S\times S$ its factors are nonnegative, and the kernel is strictly positive except possibly at the endpoint distance $\Rpos$, a set of product measure zero. If the positive and negative sets both have positive measure, the first factor is strictly positive on their product, giving strict inequality. Constant-sign functions give equality.
\end{proof}

For signed functions we apply Steiner or Schwarz rearrangement to $|u|$; we do not assert the nonnegative polarization theorem for a direct sorting of signed values.

\begin{corollary}[Rearrangement of the modulus for signed functions]\label{cor:signed}
Let $N\ge1$, $1<p<\infty$, $0<s<1$, and let real-valued $u\in\X$
vanish almost everywhere outside a measurable set $S$ with
$\diam S\le\Rpos$. Then
\[
 \Lsp((|u|)^{\star_e})\le\Lsp(u),\qquad
 \Lsp((|u|)^*)\le\Lsp(u).
\]
For nonzero $u$, equality in the first comparison holds precisely when
\[
 u(z+te)=\eta (|u|)^{\star_e}(z+(t-a)e),\qquad \eta\in\{-1,1\},\ a\in\R.
\]
For the second comparison it holds precisely when $u=\eta (|u|)^*(\cdot-a)$ for $a\in\R^N$.
\end{corollary}
\begin{proof}
Combine Proposition~\ref{prop:modulus} with Theorems~\ref{thm:steiner} and \ref{thm:schwarz}. Equality in the chain requires equality at both steps.
\end{proof}

\begin{corollary}[Signed zero exterior functions]\label{cor:domains}
Let $\Omega\subset\R^N$ be bounded and open with $\diam\Omega\le\Rpos$.
For every real-valued $u\in X_0^{s,\log,p}(\Omega)$,
\[
 (|u|)^{\star_e}\in X_0^{s,\log,p}(\Omega^{\star_e}),\qquad
 (|u|)^*\in X_0^{s,\log,p}(\Omega^*),
\]
and their energies do not exceed $\Lsp(u)$. Equality is characterized
by Corollary~\ref{cor:signed}.
\end{corollary}
\begin{proof}
Use Corollary~\ref{cor:signed} for the zero extension and inclusion
preservation of the set rearrangements. The nonnegative domain
statements are already contained in Corollary~\ref{cor:boundedPS}.
\end{proof}

\begin{remark}
Equality for one function does not force the ambient domain to be symmetric if that function vanishes on part of the domain. If $u>0$ almost everywhere in $\Omega$, Steiner equality does imply that $\Omega$ agrees almost everywhere with one translate of $\Omega^{\star_e}$ in direction $e$; Schwarz equality gives a translated ball. No regularity of the domain boundary is needed for these almost-everywhere statements.
\end{remark}

\begin{remark}[The positive-part energy]
The positive-part kernel $k_s^+(r)=\max\{k_s(r),0\}$ is nonnegative and radially nonincreasing, and therefore satisfies the three comparisons in Lemma~\ref{lem:positivekernel} for every nonnegative bounded-support $u\in\X$, without a diameter restriction. Its equality theory, however, must be treated separately: a finite interaction radius can allow components or variations separated beyond the interaction range to become undetectable. The full signed kernel, the positive-part energy, and the endpoint renormalized logarithmic energy are distinct functionals.
\end{remark}

\end{document}